\documentclass[12pt]{amsart}

\usepackage[leqno]{amsmath}
\usepackage{amssymb}
\usepackage{amsxtra}
\usepackage{amscd}
\usepackage{hyperref}
\usepackage{graphicx, color}
\usepackage{dcpic}

\newtheorem{thm}{Theorem}
\newtheorem{prop}[thm]{Proposition}

\newtheorem{defn}[thm]{Definition}

\newtheorem{remark}[thm]{Remark}

\newtheorem{lemma}[thm]{Lemma}
\newtheorem{cor}[thm]{Corollary}

\newcommand{\w}[1]{\widetilde {#1}}
\newcommand{\ww}{\widetilde W}
\newcommand{\by}{\partial^{\boxtimes}}
\newcommand{\bo}{\boxtimes}
\newcommand{\bys}{\partial_{\bullet}^{\boxtimes}}
\newcommand{\bos}{\boxtimes_{\bullet}}
\newcommand{\vr}{\varphi}
\newcommand{\ze}{\zeta}

\newcommand{\al}{\alpha}
\newcommand{\be}{\beta}
\newcommand{\ga}{\gamma}
\newcommand{\de}{\delta}

\newcommand{\field}[1]{\mathbb{F}_{#1}}
\newcommand{\si}{\sigma}

\newcommand{\R}{\mathbb{R}}

\newcommand{\ring}[1]{\mathbb{P}_{#1}}
\newcommand{\F}{\mathbb{F}}
\newcommand{\I}{\ensuremath{\mathbb{I}}}
\newcommand{\N}{\mathbb{N}}
\newcommand{\res}[1]{\ensuremath{\mathrm{\textsc{res}}(#1)}}
\newcommand{\half}{\mathbb{H}}

\newcommand{\Z}{\mathbb{Z}}

\newcommand{\disp}[1]{\displaystyle{#1}}

\newcommand{\pairing}[2]{\big\langle\,#1\,,\,#2\,\big\rangle}
\newcommand{\m}[1]{m_{#1,\bullet}}
\newcommand{\mm}[1]{m'_{#1,\bullet}}

\newcommand{\lk}[1]{\mathcal{#1}}

\newcommand{\diagram}[4][1.0]{\begin{center}\begin{figure}\includegraphics[scale=#1]{#3}\caption{#4}\label{fig:#2}\end{figure}\end{center}}

\newcommand{\graph}[1]{\Gamma_{L}}
\newcommand{\circles}[1]{\ensuremath{\mathrm{\textsc{bc}}(#1)}}
\newcommand{\fcir}[1]{\ensuremath{\mathrm{\textsc{fcir}}(#1)}}
\newcommand{\cir}[1]{\ensuremath{\mathrm{\textsc{cir}}(#1)}}
\newcommand{\arcs}[1]{\ensuremath{\mathrm{\textsc{arc}}(#1)}}

\newcommand{\inlinediag}[2][0.33]{\includegraphics[scale=#1]{#2}}

\newcommand{\cross}[1]{\ensuremath{\mathrm{\textsc{cr}}(#1)}}

\newcommand{\dec}[1]{\ensuremath{\mathrm{\textsc{dec}}(#1)}}

\newcommand{\match}[1]{\ensuremath{\mathrm{\textsc{Match}}(#1)}}
\newcommand{\state}[1]{\ensuremath{\mathrm{\textsc{State}}(#1)}}
\newcommand{\stn}[1]{\ensuremath{\mathrm{\textsc{st}}_n(#1)}}

\newcommand{\resolution}[1]{\ensuremath{\mathrm{\textsc{res}}(#1)}}
\newcommand{\complex}[1]{\ensuremath{[\![#1]\!]}}

\newcommand{\weigh}[2] {(\mu_{\mathcal{B}\Gamma_{n}} \otimes \ensuremath{\mathbb{\I}}_d)(\ensuremath{\mathbb{\I}} \otimes {#1}){#2} }
\newcommand{\weighv}[2] {(\ensuremath{\mathbb{\I}} \otimes {#1}){#2} }
\newcommand{\lefty}[1]{\overleftarrow{#1}}
\newcommand{\righty}[1]{\overrightarrow{#1}}

\newcommand{\rightDelta}[1]{\overrightarrow{\delta_{#1}}}
 \newcommand{\rightvde}[1]{\overrightarrow{\delta _{\mathcal{#1}}}}
\newcommand{\rightTde}[1]{\overrightarrow{\delta _{#1,\bullet}}}

\newcommand{\bridgeGraph}[1]{\Gamma_{#1}}

\newcommand{\dcomplex}[2]{\ensuremath{[\![\, \overrightarrow#1 \,\rangle\!\!\!\rangle} \longrightarrow \mathcal{B}\Gamma_{n} \otimes_{\lk I_n} \ensuremath{[\![\, \overrightarrow#2 \,\rangle\!\!\!\rangle}}
\newcommand{\acomplex}[2]{\ensuremath{\langle\!\!\!\langle\,\overleftarrow #1\,]\!]} \otimes_{\lk I_n}  \mathcal{B}\Gamma_{n}  \longrightarrow\ensuremath{\langle\!\!\!\langle\,\overleftarrow #2\,]\!]}}
\newcommand{\awcomplex}[2]{\ensuremath{\langle\!\!\!\langle\,\overleftarrow #1\,]\!]} \otimes_{\lk I_2}  
\mathcal{B}\Gamma_{2}  \longrightarrow\ensuremath{\langle\!\!\!\langle\,\overleftarrow #2\,]\!]}}

\newcommand{\cleaves}[1]{\mathcal{C\!L}^*_{#1}}

\newcommand{\leftBridges}[1]{\ensuremath{\overleftarrow{\mathrm{\textsc{Br}}}(#1)}}
\newcommand{\rightBridges}[1]{\ensuremath{\overrightarrow{\mathrm{\textsc{Br}}}(#1)}}
\newcommand{\bridges}[1]{\ensuremath{\mathrm{\textsc{Bridge}}(#1)}}

\newcommand{\rightcomplex}[1]{\ensuremath{[\![\, \overrightarrow{#1} \,\rangle\!\!\!\rangle} }
\newcommand{\leftcomplex}[1]{\ensuremath{\langle\!\!\!\langle\,\overleftarrow #1\,]\!]}}

\title{TWISTING BORDERED KHOVANOV HOMOLOGY}
\author{Nguyen D. Duong}

\begin{document}
\maketitle
\begin {abstract}
We describe a bordered version of totally twisted Khovanov homology. We first twist Roberts's type $D$ structure by adding a ``vertical" type $D$ structure which generalizes the vertical map in twisted tangle homology. One of the distinct advantages of our type $D$ structure is that it is homotopy equivalent to a type $D$ structure supported on ``spanning tree" generators. We also describe how to twist Roberts's type $A$ structure for a left tangle in such a way that pairing our type $A$ and type $D$ structures will result in the totally twisted Khovanov homology.     
\end {abstract}
\section {Introduction} \label {sect:1}

\noindent Modeled from the construction of the bordered Heegaard-Floer homology given by R. Lipshitz, P. Ozsvath, D. Thurston \cite {LOT}, L. Roberts, in \cite {LR1} and \cite {LR2}, described how to obtain the Khovanov homology of a link from a diagram divided into two parts, left and right tangles. To those tangles, he associated different types of tangle invariants in such a way that gluing those two invariants recovers the Khovanov homology. There are several ways to extend Khovanov homology to tangles by M. Asaeda, J. Przytycki and A. Sikora in \cite {APS}, A. Lauda and H. Pfeiffer in \cite {PLau}, D. Bar-Natan in \cite {Bar2} or M. Khovanov in \cite {Khta}. In \cite {NR}, by using the idea of the totally twisted Khovanov homology in \cite {LR3}, the author and L. Roberts took the construction of M. Asaeda, J. Przytycki and A. Sikora and twisted them to obtain a tangle invariant (called twisted tangle homology). Furthermore, reduced Khovanov homology can be described neatly and explicitly in terms of spanning trees (at least for the case of knots, proved by T. Jaeger \cite {TJ}). In an effort to combine those theories and to study the reduced Khovanov homology, we will twist Roberts's tangle invariants to get the twisted invariants for left and right tangle components of a link. By using an algebraic tool in \cite {LOT}, the result of gluing those two tangle invariants is the totally twisted Khovanov homology of links.
\ \\
\ \\
Before we describe the chain complex, let's recall the setting of \cite {LR1}, \cite {LR2}. Let $T$ be a link diagram of a link $\lk T$ in $\R^2$, which is transverse to the $y$-axis. The $y$-axis divides $T$ into two parts: a left tangle $\lefty {T}$ and a right one $\righty {T}$. Two left (right) tangles will be equivalent if they are related by an ambient isotopy preserving the $y$-axis pointwise. 
\ \\
Labeling each arc $f$ of ${\righty T}$  with a formal variable $x_{f}$, we form the polynomial ring $\ring{\righty T}=$ $\Z_2[x_f|f\in \arcs{\righty T}]$ where $\arcs{\righty T}$ denotes the set of arcs of $\righty {T}$. The field of fractions of $\ring{\righty T}$ will be denoted $\field{\righty {T}}$. The same thing can be done for $\lefty {T}$ to get $\field{\lefty {T}}$.
\ \\
\noindent We then associate to each right (left) tangle diagram $\righty {T}$ (respectively $\lefty {T}$) a bigraded vector space $\rightcomplex {T}$ over $\field{\righty {T}}$ (respectively $\leftcomplex {T}$ over $\field{\lefty {T}}$), generated by a collection of  states $(r,s)$, described as follows:\\
\begin {enumerate}
\item $r$ is a pair $(\lefty {m},\rho)$ (respectively $(\rho,\righty {m})$) where $\lefty {m}$ ($\righty {m}$) is a left (right) planar matching (a left (right) planar matching is a collection of $n$ embedded arcs in the left (right) half plane whose boundaries are the intersection points of $T$ and the $y$-axis) and $\rho$ is a resolution of $\righty {T}$ ($\lefty {T}$).
\item $s$ is a decoration on each circle of $\lefty {m}\cup \rho$ (respectively $\rho \cup \righty {m}$) by either $+$ or $-$.
\end {enumerate}
In \cite {LR1}, L. Roberts defines a type $D$ structure $\rightDelta {T}$ on $\rightcomplex {T}$ over a differential graded algebra $(\mathcal{B}\Gamma_{n},\lk I_n)$ such that the homotopy class as a type $D$ structure of $(\rightcomplex {T},\rightDelta {T})$ is an invariant of $\righty {\lk {T}}$. Inspired by the idea of the twisted tangle homology \cite {NR}, we will define a $(0,-2)$ ``vertical" type $D$ structure $\rightvde {V}$ on $\rightcomplex {T}$. With respect to the bigrading of $\rightcomplex {T}$ and $\mathcal{B}\Gamma_{n}$, $\rightDelta {T}$ is $(1,0)$ map while $\rightvde {V}$ is $(0,-2)$ map. By collapsing the bigrading using the formula $\ze=i-j/2$, both $\rightDelta {T}$ and $\rightvde {V}$ become order $1$ maps. Note that this is the usual $\delta$-grading on Khovanov homology but we will not use $\delta$ here since it overlaps the notation of the type $D$ structure $\rightDelta {T}$. We will also prove that $\rightvde {V}$ commutes with $\rightDelta {T}$ in a sense that $\rightTde {T}:=\rightDelta {T}+\rightvde {V}$ is a type $D$ structure on $\rightcomplex {T}$ (see proposition \ref {prop:11}).  Furthermore, using a trick to move weights (the formal variables labeled to the arcs) developed by T. Jaeger with the aid of stable homotopy equivalence, we will prove the following theorem in section \ref 
{sect:IVR} :
\begin {thm} \label {thm:main}
The homotopy class of the type $D$ structure $(\rightcomplex {T}, \rightTde {T})$ is (stably) invariant under the first, second and third Reidemeister moves applied to $\righty {T}$
\end {thm}
\noindent Additionally, in section \ref {sect:deformation}, using the cancellation lemma, we will get a homotopy equivalent type $D$ structure of $(\rightcomplex {T}, \rightTde {T})$ supported on states which do not contain any free circles (a free circle is the one which does not intersect with the $y$-axis). The collection of such states will be denoted $\stn {\righty {T}}$. We also denote:\\
$$
\rightcomplex {CT}:=span_{\F_{\righty {T}}}\{(r,s)\in \stn {\righty {T}}\}
$$      
We now describe a type $D$ structure $\rightDelta {T,n}$ on $\rightcomplex {CT}$, which is a deformation retraction of $(\rightcomplex {T}, \rightTde {T})$. The map
$$
\rightDelta {T,n}: \rightcomplex {CT} \longrightarrow \mathcal{B}\Gamma_{n}\otimes_{\lk I_n}\rightcomplex {CT} [-1] 
$$
is defined by specifying its image on each generator $(r,s) \in \stn {\righty {T}}$ as: 
\ \\
$\rightDelta {T,n}(r, s) =\sum_{(r',s')} \big \langle(r,s), (r',s') \big \rangle I_{\partial(r,s)} \otimes (r',s')    + \disp{\sum_{\gamma \in \bridges{r}}  B(\gamma)}$ 
$$
\hspace*{7.8cm}+\disp{\sum_{C \in \cir{\partial(r,s)}, s(C) = +} (\lefty{e_{C}}+\righty {w_C} \righty {e_C}) \otimes (r,s_{C}) }
$$
The last two terms will be defined in section \ref  {sect:3}. We now need to clarify the notation of the first term. The coefficient $\pairing{(r,s)}{(r',s')}$ is calculated from the weights assigned to the arcs. First, $\pairing{(r,s)}{(r',s')} = 0$ unless they satisfy the two following conditions:
\begin {enumerate}
\item $r$ is different from $r'$ at two crossings $c_{1}$ and $c_{2}$ which are resolved with the $0$-resolution in $r$ and the $1$-resolution in $r'$.
\item $\partial (r,s)=\partial (r',s')$ where $\partial(r,s)$, $\partial (r',s')$ are the cleaved links associated to $(r,s)$ and $(r',s')$ respectively (we will recall the definition of cleaved links in section \ref {sect:clv} and the definition of $\partial (r,s)$  in section \ref {sect: twisted tangle}).   
\end {enumerate}

\noindent If $(r,s)$ and $(r',s')$ satisfy the above conditions, let $r_{01}$ ($r_{10}$) be the resolution where $c_{1}$ is resolved with the $0$-resolution (respectively $1$-resolution) and $c_{2}$ resolved with the $1$-resolution (respectively $0$-resolution). Furthermore, if $C$ is a free circle, let $\righty {w_C} = \sum_{f \in \arcs{C}} x_{f}$. Then:\\ 
\begin {equation}
\pairing{(r,s)}{(r',s')} = \left\{\begin{array}{cl} 
\displaystyle{1/{\righty {w_{C_{10}}}} + 1/{\righty {w_{C_{01}}}} } & r_{10}\mathrm{\ and\ }r_{01}\mathrm{\ have\ free\ circles\ }C_{10}\mathrm{\ and\ }C_{01}\mathrm{\ respectively} \\
\displaystyle{1/{\righty {w_{C_{10}}}} } & r_{10}\mathrm{\ has\ a\ free\ circle\ }C_{10}\mathrm{\ but\ }r_{01}\mathrm{\ does\ not}\\
\displaystyle{1/{\righty {w_{C_{01}}}}}  & r_{01}\mathrm{\ has\ a\ free\ circle\ }C_{01}\mathrm{\ but\ }r_{10}\mathrm{\ does\ not}\\
0 & \mathrm{\ neither\ }r_{10}\mathrm{\ nor\ }r_{01}\mathrm{\ has\ a\ free\ circle}\\
\end{array}\right.
\end {equation}
Then, as the consequence of theorem \ref {thm:main}, we have the following corollary:
\begin {cor}
The homotopy class of type $D$-structure $(\rightcomplex {CT},\rightDelta {T,n})$ is an invariant of $\righty {\lk T}$.
\end {cor}
\ \\
In \cite {LR2}, associated to a left tangle diagram $\lefty {T}$ of a tangle $\lefty {\lk T}$, there exists a type $A$ structure $(\leftcomplex {T},m_1,m_2)$ over $(\mathcal{B}\Gamma_{n},\lk I_n)$ where $m_1$ is a modified version of the differential defined by M. Asaeda, J. Przytycki and A. Sikora  and $m_2$ is a right action $\acomplex {T}{T}$. In section \ref {sect:Adef}, we will define a new type $A$ structure $(\leftcomplex {T},\m {1},\m {2})$ where $\m {1}$ is obtained from $m_1$ by adding a modified version of the vertical differential in the twisted tangle homology. The only difference between the actions $m_2$ and $\m {2}$ is the action of left decoration elements on $\leftcomplex {T}$. In section \ref  {sect: Ainvariant}, we will sketch a proof which shows that $(\leftcomplex {T},\m {1},\m {2})$ is an invariant of $\lefty {\lk T}$ in the category of $A_{\infty}$ modules:
\begin {thm}\label {thm:A}
Let $\lefty {\lk T}$ be a left tangle with a diagram $\lefty {T}$. The homotopy class of $(\leftcomplex {T},\m {1},\m {2})$ is an invariant of the tangle $\lefty {\lk T}$.
\end {thm}   
\noindent Using the gluing theory described in \cite {LOT}, we can pair the type $D$ structure $(\rightcomplex {T},\rightTde  {T})$ and the type $A$ structure $(\leftcomplex {T},\m {1},\m {2})$ to form a chain complex whose underlying module is $\leftcomplex {T}\otimes_{\lk I_n} \rightcomplex {T}$ and differential is defined by the following formula:\\
$$
\bys(x\otimes y)=\m {1}(x)\otimes y + (\m {2}\otimes\I)(x\otimes \rightTde {T}(y)) 
$$
\noindent Let $T$ be the link diagram obtained by gluing $\lefty {T}$ and $\righty {T}$ along their end points. In section \ref {sect:gluing}, we will prove one of the main theorems of this paper:\\
\begin {thm}
$(\leftcomplex {T}\otimes_{\lk I_n} \rightcomplex {T},\bys )$ is chain isomorphic to $(\complex {T},\widetilde {\partial})$ where $(\complex {T},\widetilde {\partial})$ denotes the totally twisted Khovanov chain complex associated to $T$.
\end {thm} 
\ \\
\subsection {The structure of this paper} In section \ref {sect:clv}, we recall the definition of the cleaved algebra $\mathcal{B}\Gamma_{n}$ in \cite {LR1} with two minor modifications: the ground ring $\Z_2$ instead of $\Z$ and the requirement on decorations of marked circles. In section \ref {sect: twisted tangle}, we construct the expanded complex $\rightcomplex {T}$ and state our main result for the twisted tangle homology. It will be followed by the construction of the ``vertical" type $D$ structure $(\rightcomplex {T},\rightvde {V})$ in section \ref {sect:3}. In section \ref {sect:dproof}, we will prove that $\rightTde {T}:=\rightDelta {T}+\rightvde {V}$ is a type $D$ structure on $\rightcomplex {T}$. We, next, will prove that $(\rightcomplex {CT},\rightDelta {n,T})$ is the deformation retraction of $(\rightcomplex {T}, \rightTde {T})$ in section \ref {sect:deformation}. The proof of theorem \ref {thm:main} will be described in section \ref {sect:IVR} by using a trick to move weights, described in section \ref {sect:movingweight}. The whole section \ref {sect:stable} is devoted to establish the definitions of the stable homotopy equivalence of type $A$ and type $D$ structures. In section \ref {sect:Adef} and \ref {sect:Aprp}, we define the type $A$ structure on $\leftcomplex {T}$. We will sketch the proof of theorem \ref {thm:A} in section \ref  {sect: Ainvariant}. The relationship between the totally twisted Khovanov homology and the chain complex obtained by pairing our twisted type $A$ and type $D$ structures will be described in section \ref {sect:gluing}. In the last section \ref {sect: example}, we will give examples calculating our type $D$ and type $A$ structures for several knots and links.
\ \\
\subsection {Acknowledgment } I would like to thank Lawrence Roberts for a lot of helpful conversations, corrections to the earlier draft of this paper and for letting me to use his latex code to recall the definitions and properties of the cleaved algebra in section \ref {sect:clv}.

\newpage
\section{Recall the definition of The algebra from cleaved links in \cite {LR1}} \label {sect:clv}
\noindent Let $P_n$ be the set of $2n$ points $p_1=(0,1)$, $p_2=(0,2), $..., $p_{2n}=(0,2n)$ on the $y$-axis. In \cite {LR1}, Roberts associates $P_n$ with a bigraded differential algebra $\mathcal{B}\Gamma_{n}$ (over $\Z$) equipped with a $(1,0)$-differential $d_{\Gamma_{n}}$ which satisfies a Leibniz identity.  \\
\ \\
\noindent In this section, we will recall the definitions and properties of $(\mathcal{B}\Gamma_{n},d_{\Gamma_{n}})$. As we mentioned earlier, we will replace the ground ring $\Z$ by $\Z_2$ in our definitions and thus, we do not need to worry about the sign issues.  \\
\ \\
\noindent  We denote the closed half-planes $(-\infty,0]\times \R$ and $\R\times [0,\infty)$ by $\lefty {\half}$ and $\righty {\half}$ respectively. Since $\mathcal{B}\Gamma_{n}$  is described by generators and relations, we first recall the definitions of its generators and then rewrite the set of the relations, ignoring the signs. 
\begin {defn} \cite {LR1}
A $n$-decorated, cleaved link $(L,\sigma)$ is an embedding of disjoint circles in $\R^2$ such that:
\begin {enumerate}
\item Each circle contains at least 2 points in $P_n$,
\item $L\cap \{0\times (-\infty,\infty)\}=P_n$, 
\item $\sigma$ is a function which assigns either $+$ or $-$ to each circle of $L$, called decoration.
\end {enumerate}
\end {defn}
\noindent Let $\cir{L}$ be the set of circle components of $L$. We call $p_{2n}$ the marked point and the circle of $L$ which passes through the marked point is called the marked circle. We denote the marked circle by $L(*)$. We then denote $\cleaves {n}$ the set of equivalence classes of $n$-decorated, cleaved links whose decorations on the marked circle is $-$.
\begin {defn} \cite {LR1}
A right $($left$)$ planar matching $M$ of $P_n$ in the right $($left$)$ half plane $\righty {\half}$ $(\lefty {\half})$ is a proper embedding of $n$-arcs $\alpha_i:[0,1]\rightarrow \righty {\half}$ $(\lefty {\half})$ such that $\alpha_i(0)$ and $\alpha_i(1)$ belong to $P_n$. 
\end {defn}
\noindent $L$ can be obtained canonically by gluing a right planar matching $\righty {L}$ and a left one $\lefty {L}$ along $P_n$. We denote the equivalence classes of right and left planar matchings by $\righty {\match {n}}$ and $\lefty {\match {n}}$ respectively.\\
\ \\
\noindent We next describe the definition of a bridge of a cleaved link.
\begin {defn}\cite {LR1}
A bridge of a cleaved link $L$ is an embedding of $\gamma:[0,1]\rightarrow \R^2\backslash\{0 \times (-\infty,\infty)\}$ such that:
\begin {enumerate}
\item $\gamma (0)$ and $\gamma (1)$ are on distinct arcs of $\righty {L}$ $($or $\lefty {L}$$)$
\item $\gamma (0,1)\cap L=\emptyset$
\end {enumerate}
\end {defn}
\noindent Depending upon the location of $\gamma$, $\gamma$ is called either left bridge or right bridge. The equivalence classes of bridges, left bridges and right bridges of $L$ are denoted by $\bridges {L}$, $\leftBridges{L}$ and $\rightBridges{L}$ respectively.  \\
\ \\
\noindent A generator $e$ of $\mathcal{B}\Gamma_{n}$ can be thought as an oriented edge in a directed graph whose source and target vertices, denoted by $s(e)$ and $t(e)$, are decorated cleaved links (see \cite {LR1} for more detail). Corresponding to each $(L,\sigma) \in \cleaves {n} $, there is an idempotent $I_{(L,\sigma)}\in \mathcal{B}\Gamma_{n}$. $\mathcal{B}\Gamma_{n}$ is freely generated over $\Z_2$ by the idempotents and the following elements, subject to the relations described below: \\
\begin{enumerate}
\item For each circle $C \in \cir{L}$ where $\sigma(C)=+$, we have a ``dual" decorated cleaved link $(L,\sigma_C)$ where $\sigma_C (C)=-$ and $\sigma_C(D)=\sigma(D)$ for each $D \in \cir{L}\backslash\{C\}$. There are two elements $\righty{e_{C}}$ and $\lefty{e_{C}}$, called {\em right and left decoration elements}, whose sources are $(L,\sigma)$ and targets are $(L,\sigma_C)$. $C$ is called the support of $\righty{e_{C}}$ and $\lefty{e_{C}}$.
\item Let $\gamma \in \bridges{L}$, then there is a {\em bridge element} $e_{(\gamma; \sigma , \sigma_{\gamma})}$ whose source is $(L,\sigma)$ and target is $(L_{\gamma},\sigma_{\gamma})$ where $L_{\gamma}$ is obtained from $L$ by surgering along $\gamma$ and $\sigma_{\gamma}$ is any decoration compatible with $\sigma$ and the Khovanov sign rules. Additionally, $e_{(\gamma; \sigma , \sigma_{\gamma})}$ is called a left (right) bridge element if $\gamma\in \leftBridges {L}$ ($\rightBridges {L})$ and will be denoted by $\lefty{e_{\gamma}}$ ($\righty{e_{\gamma}}$) if the context is clear. Note that $L_{\gamma}$ has a special bridge $\gamma^{\dagger}$ which is the image of the co-core of the surgery.
\end{enumerate}
\ \\
\noindent With these generators and idempotents, we have:
\begin{prop} \cite{LR1}
$\mathcal{B}\bridgeGraph{n}$ is finite dimensional
\end{prop}

\noindent Furthermore, $\mathcal{B}\bridgeGraph{n}$ can be bigraded as in \cite{LR1}. In this paper, we collapse the bigrading by using $\ze(i,j)=i-j/2$ to give a new grading on $\mathcal{B}\bridgeGraph{n}$. On the generating elements, the new grading is specified by setting:
$$
\begin{array}{llcl}
\mathrm{I}_{(L,\sigma)} & \longrightarrow (0,0)&  \longrightarrow  0{\hspace {0.43cm} }\\
\righty{e_{C}} & \longrightarrow (0,-1) & \longrightarrow 1/2\\
\lefty{e_{C}} & \longrightarrow  (1,1) & \longrightarrow  1/2 \\
\righty{e_{\gamma}} &\longrightarrow  (0,-1/2) & \longrightarrow  1/4 \\
\lefty{e_{\gamma}} &\longrightarrow  (1,1/2) & \longrightarrow  3/4\\
\end{array}
$$
\noindent This assignment provides the grading to every other element by extending the grading on generators homomorphically.
\ \\
\ \\
\noindent Based on the above set of generators, there is a set of commutativity relation of generators, divided into the following cases:
\ \\
\ \\
$\bf Disjoint$ $\bf support$ $\bf and$ $\bf squared$ $\bf bridge$ $\bf relations:$ \label {rel:dsasb}
We describe the set of relations in this case by using the following model:
\begin{equation}
e_{\alpha}e_{\beta'} = e_{\beta}e_{\alpha'} 
\end{equation} 
\noindent We require that $e_{\alpha}$ and $e_{\alpha'}$ are the same type of elements (decoration or bridge) and they also have same locations (left or right). The same requirements are applied for the pair $e_{\beta}$ and $e_{\beta'}$.
\ \\
\ \\
Let $(L,\sigma)\in \cleaves{n}$ such that $I_{(L,\sigma)}$ is the source of both $e_{\alpha}$ and  $e_{\beta}$, we have the following cases:
\begin{enumerate}
\item \label {rel:1} If $C$ and $D$ are two distinct $+$ circles of $(L,\sigma)$, there are two ways to obtain $(L,\sigma_{C,D})$ from $(L,\sigma)$ by changing the decoration on either $C$ or $D$ from $+$ to $-$ first and then changing the decoration on the remaining $+$ circle. The recorded algebra elements for two paths will form a relation.
\item \label {rel:2}If $e_{\alpha}=e_{(\gamma,\sigma,\sigma')}$ for a bridge $\gamma$ in $(L,\sigma)$ and $e_{\beta}$ is a decoration element for $C \in \cir{L}$, with $C$  not in the support of $\gamma$, due to the disjoint support, there will exist $e_{\alpha'} = e_{(\gamma,\sigma_{C},\sigma_{C}')}$ and $e_{\beta'}$ which is a decoration element whose source is $(L_{\gamma},s')$ and target is $(L_{\gamma},s'_C)$.
\item \label {rel:3} If $e_{\alpha}= e_{(\gamma,\sigma,\sigma')}$ and $e_{\beta}=e_{(\eta,\sigma,\sigma'')}$ are bridge elements for distinct bridges $\gamma$ and $\eta$ in $(L,\sigma)$, with $\eta \in B_{d}(L,\gamma)$ and $e_{\beta'}= e_{(\eta,\sigma',\sigma''')}$ and $e_{\alpha'}=e_{(\gamma,\sigma'',\sigma''')}$ for some decoration $\sigma'''$ on $L_{\gamma,\eta}$, we obtain a commutativity relation. We note that $B_{d}(L,\gamma)$ stands for the set of bridges of $(L,\gamma)$ neither of whose ends is on an arc with $\gamma$. 
\item If $e_{\alpha}= e_{(\righty{\gamma},\sigma,\sigma')}$ and $e_{\beta}=e_{(\righty{\eta},\sigma,\sigma'')}$ are bridge elements for distinct {\em right} bridges $\gamma$ and $\eta$ in $(L,\sigma)$, and $e_{\beta'}= e_{(\righty{\delta},\sigma',\sigma''')}$ and $e_{\alpha'}=e_{(\righty{\omega},\sigma'',\sigma''')}$ 
, such that $L_{\gamma, \delta} = L_{\eta,\omega}$, and some compatible decoration $\sigma'''$, it will form a commutativity relation. 
\item If $e_{\alpha}= e_{(\lefty{\gamma},\sigma,\sigma')}$ and $e_{\beta}=e_{(\lefty{\eta},\sigma,\sigma'')}$ are bridge elements for distinct {\em left} bridges in $(L,\sigma)$, with $\lefty{\eta} \in B_{o}(L,\lefty{\gamma})$, and $e_{\beta'}= e_{(\lefty{\delta},\sigma',\sigma''')}$, $e_{\alpha'}=e_{(\lefty{\omega},\sigma'',\sigma''')}$ with $L_{\gamma, \delta} = L_{\eta,\omega}$, and some compatible decoration $\sigma'''$, we form a commutativity relation. We denote  $B_{o}(L,\lefty{\gamma})$ the set of bridges of $(L,\lefty{\gamma})$ whose one end is on the same arc at $\lefty {\gamma}$ and lying on opposite side of the arc as $\lefty {\gamma}$.
\end{enumerate}

\ \\
\noindent {\bf Other bridge relations:} \cite {LR1} Suppose $\gamma \in \leftBridges{L}$ and $\eta \in B_{\pitchfork}(L_{\gamma},\gamma^{\dagger})$ where $B_{\pitchfork}(L_{\gamma},\gamma^{\dagger})$ consists of the classes of bridges all of whose representatives intersect $\gamma^{\dagger}$, then 
$$
e_{(\gamma,\sigma,\sigma')}e_{(\eta,\sigma',\sigma'')} = 0
$$
whenever $\sigma'$ and $\sigma''$ are compatible decorations.\\
\ \\
\ \\
\noindent The second possibility is that $\lefty{\alpha} \in B_{s}(L,\lefty{\beta})$ where $B_{s}(L,\lefty{\gamma})$ stands for the set of bridges whose one end is on the same arc at $\lefty {\gamma}$ and lying on same side of the arc as $\lefty {\gamma}$. There is a natural left bridge $\lefty {\gamma}$ by sliding $\lefty {\alpha}$ over $\lefty {\beta}$. In this case there are three paths from $s(\lefty {e_\gamma}\lefty {e_\beta})$ to $t(\lefty {e_\gamma}\lefty {e_\beta})$ and they will form a relation whenever the decorations are compatible:
$$
\lefty {e_\alpha}\lefty {e_\beta}+\lefty {e_\beta}\lefty {e_\delta}+\lefty {e_{\gamma}}\lefty {e_{\eta}}=0
$$
where $\lefty {\delta}$ and $\lefty {\eta}$ are the images of $\lefty {\alpha}$ in $L_{\lefty {\beta}}$ and in $L_{\lefty {\gamma}}$ respectively.
\\

\noindent Additionally, if  there is a circle $C \in \cir{L}$ with $\sigma(C)=+$, and there are elements $\righty{e}_{(\gamma,\sigma,\sigma')}$ and $\righty{e}_{(\gamma^{\dagger},\sigma',\sigma_{C})}$ for a bridge $\gamma \in \rightBridges{L}$ then
\begin{equation}\label{rel:special}
\righty{e}_{(\gamma,\sigma,\sigma')}\righty{e}_{(\gamma^{\dagger},\sigma',\sigma_{C})} = \righty{e_{C}}
\end{equation}
\ \\
\noindent{\bf Relations for decoration edges:} When the support of $e_{C}$ is not disjoint from the bridge $\gamma$ of $\righty{e}_{(\gamma,\sigma,\sigma_{\gamma})}$, the relations are different depending upon the location of $e_{C}$.

\begin{enumerate} 
\item {\bf The relations for $\righty{e_{C}}$:} \label {rel:7} If $C_1$ and $C_2$ are two $+$ circles in $L$ and $\gamma$ is a bridge which merges $C_1$ and $C_2$ to form a new circle $C$, we then obtain the following relation: 
\begin{equation} \label {rel:4}
\righty{e_{C_{1}}}m_{(\gamma,\sigma_{C_{1}},\sigma_{C})} = \righty{e_{C_{2}}}m_{(\gamma,\sigma_{C_{2}},\sigma_{C})} = m_{(\gamma,\sigma,\sigma_{\gamma})} \righty{e_{C}}
\end{equation}

\ \\
Similarly, if $C$ is a $+$ circle in $L$ and $\gamma$ is a bridge which divides $C$ into $C_1$ and $C_2$ in $\cir{L_{\gamma}}$, then we impose the relation:
\begin{equation}
\righty{e_{C}}f_{(\gamma,\sigma_{C},\sigma_{C,\gamma})} =  f_{(\gamma,\sigma,\sigma^{1}_{\gamma})}\righty{e_{C_{1}}} = f_{(\gamma,\sigma,\sigma^{2}_{\gamma})} \righty{e_{C_{2}}}
\end{equation}
where $\sigma^{i}_{\gamma}$ assigns $+$ to $C_{i}$ and $-$ to $C_{3-i}$. 

\item {\bf The relations for $\lefty{e_{C}}$:} Since there are two types of decoration elements, in the above relations, if we replace the right decoration elements by the left ones, we obtain the following relations:
\begin{equation}\label{rel:lefty1}
\lefty{e_{C_{1}}}m_{(\gamma,\sigma_{C_{1}},\sigma_{C})} + \lefty{e_{C_{2}}}m_{(\gamma,\sigma_{C_{2}},\sigma_{C})} + m_{(\gamma,\sigma,\sigma_{\gamma})} \lefty{e_{C}} = 0
\end{equation}
\begin{equation}\label{rel:lefty2}
\lefty{e_{C}}f_{(\gamma,\sigma_{C},\sigma_{C,\gamma})} +  f_{(\gamma,\sigma,\sigma^{1}_{\gamma})}\lefty{e_{C_{1}}} + f_{(\gamma,\sigma,\sigma^{2}_{\gamma})} \lefty{e_{C_{2}}} = 0
\end{equation}
\end{enumerate}
\ \\
The main result about $\mathcal{B}\Gamma_{n}$ in \cite{LR1} is:
\begin{prop}
Let $(L,\sigma) \in \cleaves{n}$ such that there is a circle $C \in \cir{L}$ with $\sigma(C) = +$. Let $\lefty{e_{C}}$ be
the decoration element corresponding to $C$. Let
\begin{equation}\label{eq:dga}
d_{\Gamma_{n}}(\lefty{e_{C}}) = \sum e_{(\gamma,\sigma,\sigma_{\gamma})}e_{(\gamma^{\dagger},\sigma_{\gamma},\sigma_{C})}
\end{equation}
where the sum is over all $\gamma \in \leftBridges{L}$ with $C$ as active circle, and all decorations $\sigma_{\gamma}$ which define
compatible elements. Let $d_{\Gamma_{n}}(e) = 0$ for every other generator $e$ (including idempotents). Then $d_{\Gamma_{n}}$ can be extended to an order 1 differential on graded algebra $\mathcal{B}\bridgeGraph{n}$ which satisfies the following Leibniz identity:
\begin{equation}\label{eqn:Leibniz}
d_{\Gamma_{n}}(\alpha\beta) = (d_{\Gamma_{n}}(\alpha)\big)\beta + \alpha\big(d_{\Gamma_{n}}(\beta)\big)
\end{equation}
\end{prop}
\noindent $(\mathcal{B}\bridgeGraph{n}, d_{\Gamma_{n}})$ denotes this differential, graded algebra over $\Z_2$.  \\
\ \\
$\bf {Note.}$ Because we require the decoration on a marked circle is $-$, our $\mathcal{B}\Gamma_{n}$, described above, is actually a subalgebra of $\mathcal{B}\Gamma_{n}$ defined in \cite {LR1}. \\
\section {Twisted tangle homology and its expansion} \label {sect: twisted tangle}
\subsection {Twisted tangle homology } \label {sect:3.1}
\noindent In \cite {NR}, the author and L. Roberts define a twisted Khovanov homology version for tangles embedded in thickened surfaces by twisting the reduced Khovanov tangle chain complex, defined by  M. Asaeda, J. Przytycki and A. Sikora in \cite{APS}. \\
\ \\
\noindent Let $\righty T \subset \righty {\half}= [0,\infty) \times \R\times\{0\}$ be a tangle diagram for an oriented tangle $\righty {\lk T} \subset [0,\infty) \times \R^2 \subset \R^3$. The set of crossings in $\righty {T}$ will be denoted $\cross{\righty T}$ and the set of arcs will be denoted $\arcs{\righty T}$. The number of positive crossings will be denoted $n_{+}(\righty T)$ and the number of negative crossings will be denoted $n_{-}(\righty T)$. We will often omit the reference to $\righty {T}$ when the choice of the diagram is clear.\\
\ \\
\noindent Following \cite{LR3}, \cite{NR}, we label each arc $f \in \arcs{\righty T}$  with a formal variable $x_{f}$ and form the polynomial ring $\ring{\righty T}=$ $\Z_2[x_f|f\in \arcs{\righty T}]$. The field of fractions of $\ring{\righty T}$ will be denoted $\field{\righty T}$.\\ 
\begin{defn}
\noindent For each subset $S \subset \cross{\righty T}$, the resolution $\rho_{S}$ of $\righty T$ is a collection of arcs and circles in $\righty {\half}$, considered up to isotopy, found by locally replacing each crossing $s \in \cross{\righty {T}}$ according to the rule \\
$$
\inlinediag{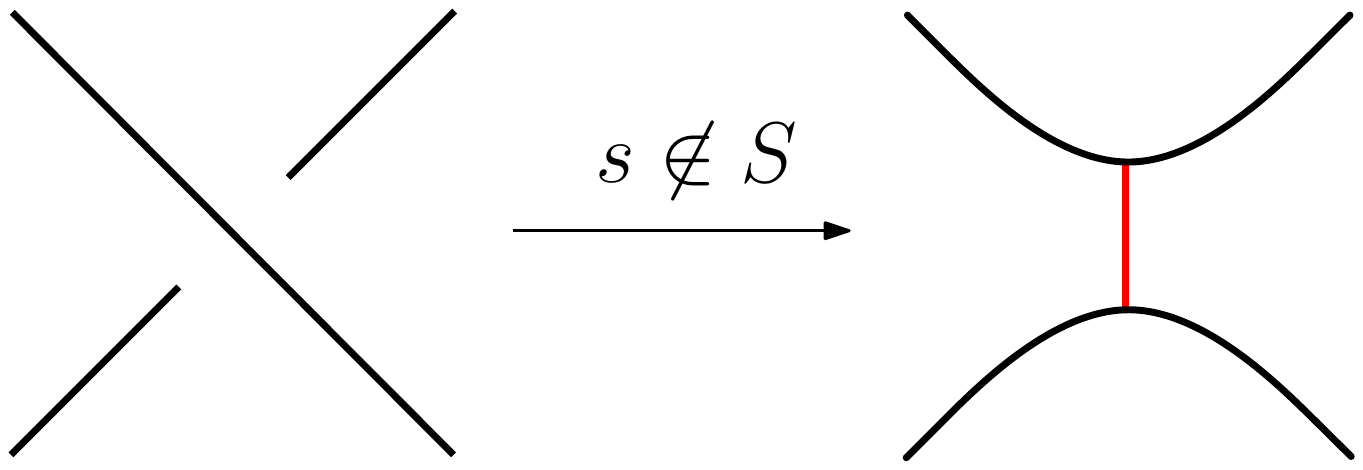} \hspace{.5in} \inlinediag{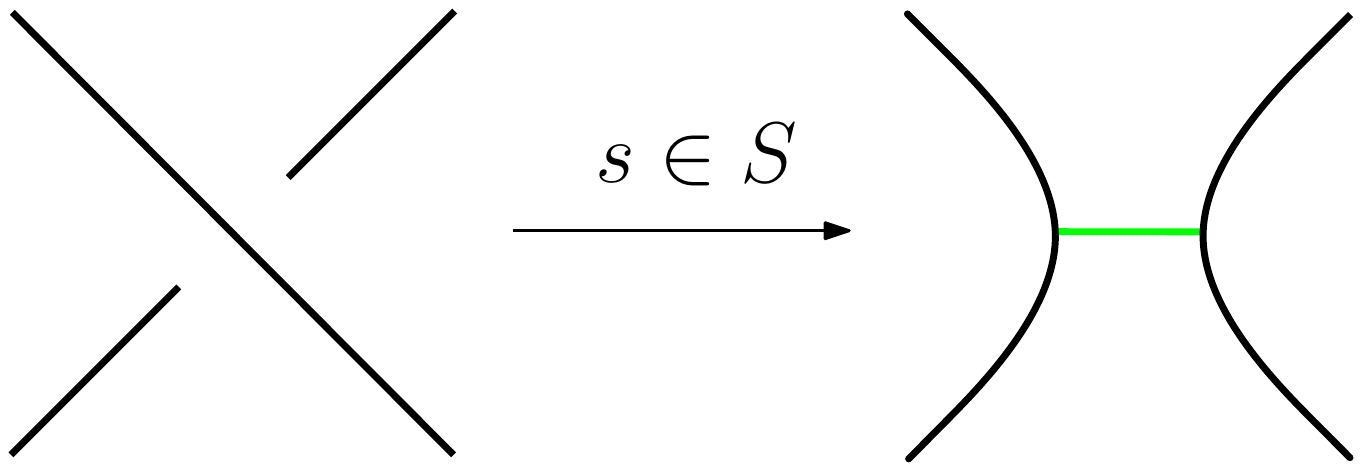} 
$$ 
The set of resolutions for $\righty {T}$ will be denoted $\res{\righty {T}}$. For each resolution $\rho$, denote by $h(\rho)$ the number of elements in the corresponding subset $S \subset \cross{\righty {T}}$. 
\end{defn}
\noindent Given a crossing $c \in \cross{\righty {T}}$ and a resolution $\rho=\rho_{S}$, we will also use the notation $\rho(c) = 0$ for $c \not\in S$ and $\rho(c)=1$ for $c \in S$. Thus, $\rho$ will stand both for the resolution diagram and for the indicator function for the set of crossings defining the resolution. The local arc $\beta$ which shows up when we resolve the crossing $c$ is called a resolution bridge. Depending on the value of $\rho$ at $c$, which is either $0$ or $1$, $\beta$ is active or inactive respectively. 
\ \\
\ \\
\noindent We denote $\rho_{\beta}=\rho \cup \{c\}$ if $\beta$ is active. Furthermore, let $\circles{\rho}$ be the set of circles and arcs in $\rho$ while $\righty {m_{\rho}}$ stands for the planar matching obtained from $\rho$ by deleting all of circles in $\circles{\rho}$. The set of free circles (the ones not crossing the $y$-axis) of $\rho$ will be denoted by $\fcir\rho$. \\
\ \\
\noindent We next assign a weight to each circle (or arc) in a resolution by adding the formal variables along each circle (or arc). \\

\begin{defn}
Let $\rho$ be a resolution for $\righty {T}$, then for each $C \in \circles{\rho}$, we define:
$$
\righty {w_C}=\sum_{f \in \arcs{C}} x_{f}
$$
\end {defn}

\noindent For each resolution $\rho(\righty{T})$, let $\fcir{\rho}=\{C_1,...,C_k\}$. To each $C_i$, we associate the complex $\mathcal{K}([C_i])$:\\
$$
0 \longrightarrow \field{\righty {T}}\,v_{+} \stackrel{\cdot \righty {w_{C_i}} }{\longrightarrow} \field{\righty {T}}\,v_{-} \longrightarrow 0
$$
where $v_{\pm}$ occur in bigradings $(0,\pm 1)$. The differential in $\mathcal{K}([C_i])$ will be denoted $\partial_{C_i}$
\ \\
Then we associate $\rho$ with the bigraded Koszul chain complex defined as
$$
\mathcal{K}(\rho) = \mathcal{K}([C_{1}], \ldots, [C_{k}]) = \mathcal{K}([C_{1}]) \otimes_{\field{\righty {T}}} \mathcal{K}([C_{2}]) \otimes_{\field{\righty {T}}} \cdots \otimes_{\field{\righty {T}}} \mathcal{K}([C_{k}])
$$
The differential in $\mathcal{K}(\rho)$ is $\partial_{\mathcal{K}(\rho)} = \sum_{i } \partial_{C_{i}}$. \\
\ \\
Because we need to shift the gradings, we use the following notation:

\begin{defn}
Let $M = \oplus_{\vec{v} \in \Z^{k}} M_{\vec{v}}$ be a $\Z^{k}$-graded $R$-module where $R$ is a ring, then $M\{\vec{w}\}$ is the $\Z^{k}$-graded module with $(M\{\vec{w}\})_{\vec{v}} \cong M_{\vec{v}-\vec{w}}$. 
\end{defn}
\noindent The vertical complex for the resolution $\rho$ of the tangle diagram $\righty {T}$ is defined as:
$$\mathcal{V}(\rho) = \mathcal{K}(\rho)\{[h(\rho),h(\rho)]\}$$ 
where the differential $\partial_{\mathcal{V}(\rho)}$ will change the bigrading by $(0,-2)$\\
\ \\
We now define a bigraded chain group:
$$
C_{\mathrm{APS}}(\righty{T}) = \bigoplus_{\rho \in \resolution{\righty{T}}} \mathcal{V}{(\rho)}
$$
$C_{\mathrm{APS}}(\righty{T})$ is the bigraded chain complex with a $(0,-2)$ differential:
$$\partial'_{\mathcal{V}}= \bigoplus_{\rho \in \res{\righty {T}}} \partial_{\mathcal{V}(\rho)}$$
Additionally, in \cite{APS}, M. Asaeda, J. Przytycki and A. Sikora define a $(1,0)$  differential $\partial_{APS}$ on this bigraded module, satisfying:
$$
\partial_{APS}\partial'_{\mathcal{V}}+\partial'_{\mathcal{V}}\partial_{APS}=0
$$
Therefore, by collapsing the bigrading using $\ze(i,j)=i-j/2$, we can make $\partial_{APS}+\partial'_{\mathcal{V}}$  a differential of $C_{\mathrm{APS}}(\righty{T})$. 
Furthermore, since both $\partial_{APS}$ and $\partial'_{\mathcal{V}}$ are constructed to preserve the right planar matching, we can  decompose: 
$$
(C_{\mathrm{APS}}, \partial_{APS}+\partial'_{\mathcal{V}})=\bigoplus_{m \in \righty {\match{n}}} (C_{\mathrm{APS}}(\righty{T},m), \partial_{APS,m}+\partial'_{\mathcal{V},m})
$$
We now can describe the main theorem of \cite {NR}:
\begin {thm}
For each ${m \in \righty {\match{n}}}$, the homology $H_{\ast}(C_{\mathrm{APS}}, \partial_{APS,m}+\partial'_{\mathcal{V},m})$, as a relative $\ze$-graded module, is an invariant of the isotopy class of $\righty {\lk T}$.   
\end {thm}
\noindent Following \cite{LR1}, we first expand the chain group:\\
$$
\rightcomplex{T}:=\bigoplus_{(L,\sigma)\in \cleaves {n}}C_{\mathrm{APS}}(\righty {T},\righty {L},\sigma)=\bigoplus_{(L,\sigma)\in \cleaves {n}}C_{\mathrm{APS}}(\righty {T},\righty {L})\{0,\frac {i(L,\sigma)}{2}\}
$$
where $ {i(L,\sigma)}$ is computed by subtracting the number of $-$ non-marked circles in $L$ from the number of $+$ circles in $L$. We also denote:
$$
d_{APS}=\bigoplus_{(L,\sigma)\in \cleaves {n}}\partial_{APS,\righty {L}}
$$
$$
\partial_{\mathcal{V}}=\bigoplus_{(L,\sigma)\in \cleaves {n}}\partial'_{\mathcal {V},\righty {L}}
$$
\ \\
\ \\
\noindent As we can see, a generator of $\rightcomplex {T}$ corresponds $1$-$1$ with a triple $(\lefty {m},\rho,s)$ where $\lefty {m}$ is a left planar matching, $\rho$ is a resolution of $\righty {T}$ and $s:\cir{\lefty {m}\# \rho}  \longrightarrow \{+,-\}$ such that $s(\lefty {m}\# \righty {m_\rho}(\ast))=-$. Here $\cir{\lefty {m}\# \rho}$ is the collection of all cleaved and free circles of $\lefty {m}\# \rho$. \\
\ \\
\noindent We denote $r=\lefty {m}\# \rho$ and we call $(r,s)$ a state. The collection of states of $\righty {T}$ will be denoted $\state {\righty {T}}$. $\partial (r,s)$ denotes a decorated cleaved link, obtained by deleting all of the free circles of $\lefty {m}\# \rho$ and the decoration of this cleaved link is induced from the decoration $s$ on $\cir{\lefty {m}\# \rho}$.\\ Therefore:
$$
\rightcomplex {T}= \bigoplus_{(L,\sigma)\in \cleaves {n}}\bigoplus_{\partial (r,s)=(L,\sigma)}\F_{\righty {T}}(r,s)
$$
\ \\
Additionally, each state $(r,s)=(\lefty {m}\# \rho,s)$ has a $\ze$-grading, computed from a bigrading $(h(r,s),q(r,s))$ where $h(r,s)=h(\rho)-n_-$ and $q(r,s)=h(\rho)+\frac {i(L,\sigma)}{2}+\#(\text {+ free circles})-\#(\text {- free circles})+n_+-2n_-$.   
\ \\
\ \\
$\bf Note$. The whole process can be applied exactly the same to give the construction of the expanded complex $\leftcomplex {T}$ associated to a left tangle $\lefty {T}$.
\ \\   
\ \\
In the next section, we will describe a ``vertical" type $D$ structure on this underlying module $\rightcomplex {T}$ with respect to the $\ze$-grading.
\section {A type D structure in the twisted tangle homology} \label {sect:3}
\noindent Using the idea of the twisted tangle homology, we will describe a way twisting of Roberts's type $D$ structure on $\rightcomplex {T}$ \cite {LR1}. First of all, we define a ``vertical" type $D$ structure on $\rightcomplex {T}$.\\
\ \\
Given $\rightcomplex {T}$ equipped with the $\ze$-grading, we define a left $\lk I_n$-module map:
$$
\rightvde {V}: \dcomplex {T}{T}[-1]
$$
by specifying the image of $\rightvde{V}$ on each generator $\xi=(r,s)$ of $\rightcomplex {T}$:
\begin {equation} \label {eqn:dvertical}
\rightvde {V}(r,s)= I_{\partial (r,s)}\otimes \partial_{\mathcal{V}}(r,s) +\sum\limits_{C \in \cir{\partial (r,s)},s(C)=+}\righty {w_C} \righty {e_C}\otimes (r,s_C)
\end {equation}
$\bf Note:$ 
$\rightcomplex {T}$ is a left $\lk I_n$-module, defined by the trivial action: $I_{(L,\sigma)}(r,s)=(r,s)$ if $\partial (r,s)=(L,\sigma)$ or $0$ else. In the definition of $\rightvde {V}$ and elsewhere in this paper, $\mathcal{B}\Gamma_{n}$ stands for $\mathcal{B}\Gamma_{n}\otimes_{\Z_2} \F_{\righty {T}}$. Note that $\mathcal{B}\Gamma_{n}\otimes_{\Z_2} \F_{\righty {T}}$ has a differential graded algebra structure over $\F_{\righty {T}}$, induced by the differential graded algebra structure over $\Z_2$ of $\mathcal{B}\Gamma_{n}$ (this fact will be mentioned in section \ref {sect:stable}). Also, the decoration $s_C$ on $r$ is obtained from $s$ by changing the decoration on a cleaved circle $C$ from $+$ to $-$ and using the same decorations for other circles. Therefore, if $\partial (r,s)= (L,\sigma)$ then $\partial(r,s_C)=(L,\sigma_C)$.
\ \\
\ \\
Furthermore, the bigradings of the terms $I_{\partial (r,s)}\otimes \partial_{\mathcal{V}}$ are decreased by $(0,2)$ because that the bigradings of idempotents are $(0,0)$ and $\partial_{\mathcal{V}}$ is $(0,-2)$ differential. Therefore, the $\ze$-grading is increased by $1$. Similarly, since the $\ze$-grading of $\righty {e_C}$ is $1/2$ and $\ze(r,s)-\ze(r,s_C)=1/2[i(L,\sigma)/2-i(L,\sigma_C)/2]=1/2$, the $\ze$-gradings of the terms $\righty {e_C} \otimes (r,s_C)$ are larger by $1$. As a result, $\rightvde {V}$ is $\ze$-grading preserving. 
\ \\
\ \\
Next, we prove the following proposition:
\begin {prop} \label {prop:10} $\rightvde {V}$ is a type $D$-structure on $\rightcomplex {T}$$:$

$$
\weigh {\rightvde{V}}{\rightvde{V}}+(d_{\Gamma_{n}} \otimes \I_d) \rightvde {V}=0
$$
\end {prop}
\noindent $\bf Proof.$ Since the image of $d_{\Gamma_{n}}$ on idempotents or right decoration elements equals $0$, it suffices to verify that:
$$
\weigh {\rightvde{V}}{\rightvde{V}}(r,s)=0
$$  
for each generator $\xi=(r,s)$ of $\rightcomplex {T}$\\
\ \\
Since the image of $\rightvde {V}(r,s)$ contains states $(r,s_C)$ (with coefficients in $\mathcal{B}\Gamma_{n}$), the image of $\weighv {\rightvde{V}}{\rightvde{V}}(r,s)$ are the states $(r,s_{C_{1,2}})$ where the decoration $s_{C_{1,2}}$ of $r$ is obtained from $s$ by changing the decoration on $C_1,C_2 \in \cir{r}$ from $+$ to $-$. Therefore, we have:
$$
\weigh {\rightvde{V}}{\rightvde{V}}(r,s)=\sum_{\substack{C_1,C_2\in \cir{r}\\ C_1 \neq C_2}} A_{(r,s_{C_{1,2}})}(r,s_{C_{1,2}})
$$
We will prove that each $ A_{(r,s_{C_{1,2}})}=0$. Since each circle in $\cir{r}$ is either cleaved or free, we have the following cases:
\begin {enumerate}
\item Both $C_1$ and $C_2$ are free circles, there are two ways to obtain $(r,s_{C_{1,2}})$ from $(r,s)$:

$$
(r,s) \xrightarrow {\righty {w_{C_1}}I_{\partial (r,s)}} (r,s_{C_1}) \xrightarrow {\righty {w_{C_2}}I_{\partial (r,s)}} (r,s_{C_{1,2}})
$$
$$
(r,s) \xrightarrow {\righty {w_{C_2}}I_{\partial (r,s)}} (r,s_{C_2}) \xrightarrow {\righty {w_{C_1}}I_{\partial (r,s)}} (r,s_{C_{1,2}})
$$
Therefore, $A_{(r,s_{C_{1,2}})}= (\righty {w_{C_1}}\righty {w_{C_2}}+\righty {w_{C_2}}\righty {w_{C_1}}) I_{\partial (r,s)}=0$. Note that the elements above the arrows indicate the algebra elements in $\mathcal{B}\Gamma_{n}$.
\item $C_1$ is free and $C_2$ is cleaved. We have two following paths:
$$
(r,s) \xrightarrow {\righty {w_{C_1}}I_{\partial (r,s)}} (r,s_{C_1}) \xrightarrow {\hspace* {0.25cm}\righty {w_{C_2}}\righty {e_{C_2}}\hspace* {0.25cm}} (r,s_{C_{1,2}})
$$
$$
(r,s) \xrightarrow {\hspace* {0.25cm}\righty {w_{C_2}}\righty {e_{C_2}}\hspace* {0.25cm}} (r,s_{C_2}) \xrightarrow {\righty {w_{C_1}}I_{\partial (r,s_{C_2})}} (r,s_{C_{1,2}})
$$
and therefore,
$$
A_{(r,s_{C_{1,2}})}= \righty {w_{C_2}}\righty {w_{C_1}}\righty {e_{C_2}}I_{\partial (r,s_{C_{1,2}})} +\righty {w_{C_1}}\righty {w_{C_2}} I_{\partial (r,s_{C_1})}\righty {e_{C_2}}=0
$$ 
\item Both $C_1$ and $C_2$ are cleaved circles. We have:
$$
(r,s) \xrightarrow {\hspace* {0.25cm}\righty {w_{C_1}}\righty {e_{C_1}}\hspace* {0.25cm}} (r,s_{C_1}) \xrightarrow {\hspace* {0.25cm}\righty {w_{C_2}}\righty {e_{C_2}}\hspace* {0.25cm}} (r,s_{C_{1,2}})
$$
$$
(r,s) \xrightarrow {\hspace* {0.25cm}\righty {w_{C_2}}\righty {e_{C_2}}\hspace* {0.25cm}} (r,s_{C_1}) \xrightarrow {\hspace* {0.25cm}\righty {w_{C_1}}\righty {e_{C_1}}\hspace* {0.25cm}} (r,s_{C_{1,2}})
$$
By the relation \ref {rel:1} defining $\mathcal{B}\Gamma_{n}$, we have $\righty {e_{C_2}}\righty {e_{C_1}}+\righty {e_{C_1}}\righty {e_{C_2}}=0$. Therefore, $A_{(r,s_{C_{1,2}})}=0$.  $\Diamond$
\end {enumerate} 
\noindent In \cite{LR1}, L. Roberts discovered a type $D$ structure $\rightDelta {T} $ on the bigraded module $\rightcomplex {T}$, which is bigrading preserving into $\mathcal{B}\Gamma_{n} \otimes_{\lk I_n}\rightcomplex {T}[(-1,0)]$  and thus, is also $\ze$-grading preserving into $\mathcal{B}\Gamma_{n} \otimes_{\lk I_n}\rightcomplex {T}[-1]$.\\
\ \\ 
\ \\
As the next step, we show that ${\rightvde{V}}$ commutes to $\rightDelta {T}$ in the following sense:
\begin {prop} \label{prop:11}
The type $D$ structures ${\rightvde{V}}$ and $\rightDelta {T}$ on $\rightcomplex {T}$ satisfy:
$$
\weigh {\rightvde{V}}{\rightDelta {T}}+\weigh {\rightDelta {T}}{\rightvde{V}}=0
$$   
\end {prop}

\noindent The proof of this proposition will be presented in the next section. Combined with the fact that $\rightvde{V}$ and ${\rightDelta {T}}$ are type $D$ structures on $\rightcomplex {T}$, proposition \ref {prop:11} ensures that $\rightTde {T}={\rightvde{V}}+{\rightDelta {T}}$ is a type $D$ structure on $\rightcomplex T$. In section \ref {sect:IVR}, we will prove one of the main theorems of this paper:\\
\begin {thm}
The homotopy class of the type $D$ structure $(\rightcomplex {T}, \rightTde {T})$ is invariant under the first, second and third Reidemeister moves applied to $\righty {T}$.
\end {thm}


\section{Proof of proposition \ref {prop:11}} \label {sect:dproof}
\noindent Before proving the proposition, let us shortly recall the definition of $\rightDelta {T}$ defined in \cite {LR1}. For each generator $(r,s)$ of $\rightcomplex {T}$, let

\noindent $\rightDelta {T}(r,s):= I_{\partial (r,s)}\otimes d_{APS}+\sum\limits_{\ga\in \bridges {r}} \lk B(\ga)$\\ 
$$
\hspace* {2.5cm}+ \sum\limits_{\ga \in \dec {r,s}}\righty {e_{C(\ga)}}\otimes (r_{\ga},s_{\ga})+\sum\limits_{C \in \cir{\partial (r,s)},s(C)=+}\lefty {e_C}\otimes (r,s_C)
$$
where:
\begin {enumerate}
\item $\bridges {r}$ is the union of left bridges of $\partial (r,s)$ and the subset of active resolution bridges of $r$, containing active resolution bridges $\gamma$ such that the right planar matching of $r_{\gamma}$ is different from the one of $r$. Furthermore, $\lk B(\ga)$ is the sum of $(r_{\gamma},s_{\gamma}^i)$ where $s_{\gamma}^i$ is computed from Khovanov Frobenius algebra, with the recorded coefficient in $\mathcal{B}\Gamma_{n}$ corresponding to the bridge element whose source is $\partial(r,s)$ and target is $\partial (r_{\gamma},s_{\gamma}^i)$.  
\item $\dec{r,s}$ is a collection of active resolution bridges $\gamma$ of $r$ such that either both feet of $\gamma$ belong to a $+$ cleaved circle $C$ of $r$, or one foot of $\gamma$ belongs to a $+$ cleaved circle $C$ of $r$ and another one belongs to a $-$ free circle of $r$. In both cases, $s_{\gamma}$ is computed from Khovanov Frobenius algebra with a condition that $s_{\gamma}(C_{\gamma})=-$ where $C_\gamma$ is a cleaved circle of $r_{\gamma}$, obtained from $C$ by surgering along $\gamma$.   
\end {enumerate}
$\bf  Proof$ $\bf of$ $\bf Proposition$ $\bf \ref {prop:11}.$ It suffices to prove that for each generator $\xi=(r,s)$ of $\rightcomplex {T}$, 
$$
\weigh {\rightvde{V}}{\rightDelta {T}}(r,s)+\weigh {\rightDelta {T}}{\rightvde{V}}(r,s)=0
$$
We rewrite the left hand side as
$$
\weigh {\rightvde{V}}{\rightDelta {T}}(r,s)+\weigh {\rightDelta {T}}{\rightvde{V}}(r,s)=\sum\limits_{(r',s')} A(r',s') \otimes (r',s')
$$
 where $A(r',s')$ is computed by taking the sum of products $e_{\al}e_{\be}$, modeled by:
$$
(r,s) \xrightarrow {e_{\al}}  (r_1,s_1) \xrightarrow {e_{\be}} (r',s')
$$  
Note that ${e_{\al}}$ and $e_{\be}$ are the elements in $\mathcal{B}\Gamma_{n}$ corresponding to $\partial(r,s) \rightarrow \partial(r_1,s_1)$ and $\partial(r_1,s_1) \rightarrow \partial(r',s')$, respectively. Additionally, one of them is either an idempotent or a right decoration element (this term comes from ${\rightvde{V}}$) and the other is either an idempotent, or a bridge element, or a right decoration element, or a left decoration element (this term comes from $\rightDelta {T}$).    
\ \\
\ \\
\noindent Our goal is to prove that $A(r',s')=0$. We have the following two cases:\\
$\bf Case$ $\bf I.$ The term coming from ${\rightvde{V}}$ is a right decoration element $\righty {e_C}$. In this case, we have the following subcases:
\begin {enumerate}
\item The term coming from $\rightDelta {T}$ is a left decoration element $\lefty {e_D}$. Then $(r',s')$ is obtained from $(r,s)$ by changing the decorations on two cleaved circles $C$, $D$ of $r$ from $+$ to $-$. Therefore, there are exactly two paths from $(r,s)$ to $(r',s')$:
 $$
(r,s) \xrightarrow {\righty {w_C} \righty {e_C}}(r,s_C) \xrightarrow {\hspace {1.35ex}\lefty {e_D}\hspace {1.35ex}} (r',s')
$$ 
 $$
(r,s) \xrightarrow {\hspace {1.35ex} \lefty {e_D}\hspace {1.35ex}}(r,s_D) \xrightarrow {\righty {w_C}\righty {e_C}} (r',s')
$$
We have: $A(r',s')=\righty {w_C}[ \righty {e_C}\lefty {e_D}+ \lefty {e_D}\righty {e_C}]=0$ by relation \ref {rel:1} in the definition of $\mathcal{B}\Gamma_{n}$.
\item The term coming from $\rightDelta {T}$ is either an idempotent, or a right decoration element, or a bridge element $e_{\gamma}$ where $\gamma \in \bridges {r}$, such that the support of $\gamma$ is disjoint from $C$. In this case, there are again exactly two paths from $(r,s)$ to $(r',s')$. Similar to the previous case and the proof of proposition \ref {prop:10}, using the relation \ref {rel:1} or \ref {rel:2} in the set of disjoint support and squared bridge relations, we can compute that $A(r',s')=0$
\item The term from $\rightDelta {T}$ is a bridge element $e_{\gamma}$ where $\gamma \in \bridges{r}$ and $C$ is in the support of $\gamma$. We only describe the merging case since the case of division can be handled by a similar calculation. Let $C_1$ be a cleaved circle in $\partial (r,s)$, merged with $C$ by merging on $\ga$ to form another cleaved circle $C_2$. For the path from $(r,s)$ to $(r',s')$ to exist, $s(C_1)=+$. Then, there are three paths in $A (r',s')$:
$$
(r,s) \xrightarrow {\hspace* {0.3cm}\righty {w_C }\righty {e_C}\hspace* {0.3cm}} (r,s_C) \xrightarrow {\hspace* {0.55cm}e_\ga\hspace* {0.55cm}}(r',s')
$$
$$
(r,s) \xrightarrow {\hspace* {0.2cm}\righty {w_{C_1} }\righty {e_{C_1}}\hspace* {0.2cm}} (r,s_{C_1}) \xrightarrow {\hspace* {0.55cm}e_{\ga_1}\hspace* {0.55cm}}(r',s')
$$
$$
(r,s) \xrightarrow {\hspace* {0.55cm} e_{\ga_2}\hspace* {0.55cm}}(r_\ga,s_\ga) \xrightarrow {\hspace* {0.2cm}\righty {w_{C_2} }\righty {e_{C_2}}\hspace* {0.2cm}}(r',s')
$$ 
 As the result, $A(r',s')= \righty {w_C }\righty {e_C} e_\ga +\righty {w_{C_1} }\righty {e_{C_1}}e_{\ga_1}+e_{\ga_2}\righty {w_{C_2} }\righty {e_{C_2}}$. Since $\righty {w_{C_2}}=\righty {w_{C_1}}+\righty {w_C}$, we can rewrite: $A(r',s')=\righty {w_C}(\righty {e_C}e_\ga+e_{\ga_2} \righty {e_{C_2}})+\righty {w_{C_1}}(\righty {e_{C_1}}e_{\ga_1}+e_{\ga_2} \righty {e_{C_2}})=0$, by relation \ref {rel:4} in the set of relations for decoration edges.
\end {enumerate}
$\bf Case$ $\bf II.$ The term coming from ${\rightvde{V}}$ is an idempotent obtained by changing the decoration on a free circle $C$ from $+$ to $-$. We have the following subcases:
\begin {enumerate}
\item The term from $\rightDelta {T}$ is an idempotent. In this case, we know that the product of the weights of two paths from $(r,s)$ to $(r',s')$ will be canceled out by \cite {NR}. Therefore, $A(r',s')=0$.
\item The term from $\rightDelta {T}$ is a bridge element $e_{\gamma}$ where $\gamma \in \bridges{r}$. Since the support of $\gamma$ is disjoint from $C$, we only have two paths from $(r,s)$ to $(r',s')$ and thus $A(r',s')=\righty {w_C} I_{\partial {(r,s)}}e_\gamma+ e_\gamma \righty {w_C}  I_{\partial {(r,s)}}=0$. 
\item The term from $\rightDelta {T}$ is a right decoration element, obtained by surgery of $(r,s)$ along $\gamma \in \dec{r,s}$. It is possible that the support of $\gamma$ is disjoint or not disjoint from $C$. When the support of $\gamma$ is disjoint from $C$, the proof of $A(r,s)=0$ is similar to the case II.2. When the support of $\gamma$ is not disjoint from $C$, the situation is more interesting. Since the case of division can be handled similarly, we here only present the proof of the case when $\gamma$ merges a $+$ cleaved circle $D \in \partial {(r,s_C)}$ with $C$ to give a cleaved circle $C_1 \in \partial {(r',s')}$. In this case, there are exactly three paths from $(r,s)$ to $(r',s')$:
$$
(r,s) \xrightarrow {\hspace* {0.25cm}I_{\partial {(r,s)}} \hspace* {0.25cm}} (r_2,s_2) \xrightarrow {\hspace* {0.2cm}\righty {w_{C_1}} \righty {e_{C_1}}\hspace* {0.2cm}} (r',s')
$$  
$$
(r,s) \xrightarrow {\hspace* {0.1cm}\righty {w_C} I_{\partial {(r,s)}}\hspace* {0.1cm}} (r,s_C) \xrightarrow {\hspace* {0.5cm}\righty {e_D}\hspace* {0.5cm}} (r',s')
$$
$$ 
(r,s) \xrightarrow {\hspace* {0.25cm}\righty {w_D} \righty {e_D}\hspace* {0.25cm}} (r,s_D) \xrightarrow  {\hspace* {0.2cm}I_{\partial {(r,s_D)}}\hspace* {0.2cm}}  (r',s')
$$
Thus, $A(r',s')=(\righty {w_{C_1}}+\righty {w_{D}}+\righty {w_{C}}) \righty {e_D}=0$ because $\righty {e_D}\equiv \righty {e_{C_1}}$ and $\righty {w_{C_1}}=\righty {w_{C}}+\righty {w_{D}}$.     
\end {enumerate}

\section{The deformation retraction of the type $D$ structure} \label {sect:deformation}
\noindent Let $\stn {\righty T}$ be the collection of states of $\righty {T}$, consisting of those states that do not have any free circle in their resolutions. In this section, we will prove that the type $D$ structure $(\rightcomplex {CT},\rightDelta {n,T})$, described in section \ref  {sect:1}, is homotopy equivalent to  $(\rightcomplex {T}, \rightTde T)$ as type $D$ structures.\\
\ \\
First of all, we state the type $D$ cancellation lemma whose proof can be found in \cite {LR1}.
\ \\
\ \\
Let $(A,d)$ be a differential graded algebra over a ground ring $R$ (characteristic $2$). Let $N$ be a graded module over $R$. Suppose over $R$, $N$ can be generated by a basis $\{x_1,...,x_n\}$. Let $a_{ij} \in A$ such that:
\begin{equation} \label{eqn:structure}
d(a_{ik}) + \sum_{j=1}^{n} a_{ij}\,a_{jk} = 0  \hspace{0.5in} i,k \in \{1, \ldots, n\}
\end{equation}
and $gr(a_{ij})=|x_i|-|x_j|+1$. Then $\delta : N \rightarrow (A\otimes_R N)[-1]$, defined by:
$$
\delta(x_{i}) = \sum_{j=1}^{n} a_{ij} \otimes x_{j}
$$
is a type $D$ structure on $N$.
\begin {lemma}\cite {LR1} \label {lemma:cancel}
Let $\delta$ be a $D$ structure on $N$. Suppose there is a basis $B$ whose structure coefficients satisfy $a_{ii} = 0$ and $a_{12} = 1_{A}$. Let $\overline{N} = \mathrm{span}_{R}\{\overline{x}_{3}, \ldots, \overline{x}_{n}\}$. Then 
$$
\overline{\delta}(\overline{x}_{i}) = \sum_{j \geq 3}(a_{ij} - a_{i2}\,a_{1j}) \otimes \overline{x}_{j}
$$
is a $D$ structure on $\overline{N}$. Furthermore, the maps
$$
\begin{array}{lcl}
\iota : \overline{N} \rightarrow A \otimes N & \hspace{0.75in} & \iota(\overline{x}_{i}) = 1_{A} \otimes x_{i} -  a_{i2} \otimes x_{1} \\
\ &\ \\
\pi :  N \rightarrow A \otimes \overline{N} & \hspace{0.75in} & \pi(x_{i}) = \left\{
\begin{array}{ll}
0 & i = 1\\
\sum_{j \geq 3}a_{1j} \otimes \overline{x}_{j} & i = 2\\ 
1_{A} \otimes \overline{x}_{i} & i \geq 3
\end{array}
\right.
\end{array}
$$
realize $\overline{N}$ as a deformation retraction of $N$ with $\iota \ast \pi \simeq_{H} \I_{N}$ using the {\em homotopy} $H: N \rightarrow A \otimes N[-1]$
$$
H(x_{i}) = \left\{\begin{array}{ll} 1_{A} \otimes x_{1} & i = 2\\ 0 & i \neq 2 \end{array}\right.
$$ 
\end{lemma}
\begin {prop}
$(\rightcomplex {CT},\rightDelta {n,T})$, defined in section \ref {sect:1}, is a type $D$ structure over $(\mathcal{B}\Gamma_{n},\lk I_n)$. Additionally, $(\rightcomplex {T}, \rightTde {T})$, defined in section \ref {sect:3}, is homotopy equivalent to $(\rightcomplex {CT},\rightDelta {n,T})$.
\end {prop}
\noindent $\bf Proof.$ Let $(r,s)$ be a state of $\rightcomplex {T}$ containing a free circle $C$. Then the decoration $s(C)$ is either $+$ or $-$. Corresponding to $(r,s)$, there is a state $(r,s_1)$ of $\rightcomplex {T}$ obtained by changing the decoration on $C$ from $\pm$ to $\mp$. We call this pair of states a mutual pair. We will use the cancellation lemma \ref {lemma:cancel} to cancel out the mutual pairs $(r,s)$ and $(r,s_1)$. Indeed, $\rightTde {T} (r,s)= \righty {w_C} I_{\partial (r,s)}\otimes (r,s_1) +Q= \righty {w_C} 1_{\mathcal{B}\Gamma_{n}}\otimes (r,s_1) +Q$ where $Q$ is a linear combination supported on states not equal to either $(r,s)$ or $(r,s_1)$. Since $\righty {w_C} 1_{\mathcal{B}\Gamma_{n}}$ is invertible, we can cancel this pair to get a deformation retraction of $(\rightcomplex {T},\rightTde {T})$ supported on $\state {\righty T} \setminus \{(r,s),(r,s_1)\}$. We also note that the new perturbation does not alter the relationship between another mutual pair. As a result, we can cancel out all of the mutual pairs of states and what we have left is a homotopy equivalent type $D$ structure $\delta_n$ supported on $\stn {\righty T}$. We also need to verify that $\delta_n$ is the same as $\rightDelta {n,T}$. Let $(r,s), (r',s')$ be two states in $ \stn {\righty T}$ such that $\langle \delta_n(r,s),(r',s') \rangle \neq 0$. Therefore, under the action of $\rightTde {T}$, we have the following sequence of transitions of states in $\state {\righty T}$:\\
$$
(r,s)=(r_0,s_0^+) \rightarrow (r_1,s_1^-)\rightarrow (r_1,s_1^+) \rightarrow.... \rightarrow (r_k,s_k^+) \rightarrow (r_{k+1},s_{k+1}^-)=(r',s')
$$     
where each transition $(r_i,s_i^-) \rightarrow (r_i,s_i^+)$ comes from a mutual pair and let $C_i$ be the free circle where $s_i^-(C_i)=-$ and $s_i^+(C_i)=+$. The other transitions come from the fact that $(r_{i+1},s^-_{i+1})$ is in the image of $\rightTde {T}(r_i,s^+_i)$. By denoting the number of $+$ and $-$ free circles for each state $(r,s)$ of $\righty {T}$ by $J(r,s)=(J_+(r,s),J_-(r,s))$, we see that  $J(r_i,s_i^+)-J(r_i,s_i^-)=(1,-1)$ for each $i \in \{1,...,k\}$. Additionally, we evaluate $J_i=J(r_{i+1},s_{i+1}^-) - J(r_i,s_i^+)$ as the following cases:
\begin {enumerate}
\item  If the corresponding coefficient from $\partial (r_i,s_i^+) \rightarrow \partial (r_{i+1},s_{i+1}^-) $ is either $\lefty {e_C}$ or bridge element then $J_i=(0,0)$. 
\item  If the corresponding coefficient is $\righty {e_C}$ then $J_i$ belongs to $\{(1,0),(0,-1),(0,0)\}$. 
\item If this transition comes from $I\otimes d_{APS}$, $J_i$ belongs to $\{(-1,0), (0,1)\}$.
\end {enumerate}   
Since $(r,s), (r',s') \in \stn {\righty {T}}$, we have $J(r,s)=J(r',s')=(0,0)$. Furthermore, we have:
$$
J(r',s')-J(r,s)=\sum \limits_{i=1}^k\big[ J(r_i,s_i^+)-J(r_i,s_i^-)\big ]+\sum \limits_{i=0}^{k}J_i=(0,0)
$$
Therefore, $\sum \limits_{i=0}^{k}J_i=(-k,k)$. Looking through all of possible cases of $J_i$, $k$ has to be either $0$ or $1$. If $k=0$, $(r',s')$ is obtained from $(r,s)$ by either changing a decoration on a cleaved circle from $+$ to $-$ or performing surgery along $\gamma \in \bridges{r}$. In this case:
$$
\langle\delta_n(r,s),(r',s')\rangle=\langle\rightTde {T} (r,s),(r',s')\rangle
$$
If $k=1$, we need to have $\partial (r,s)=\partial (r_1,s_1^-)=\partial (r_1,s_1^+)=\partial (r',s')$. As a result, we know how to calculate $\delta_n$. We call a sequence of transitions an $\lk I$-transition if it is of the following form: $(r,s)\rightarrow(r_1,s_1^-)\rightarrow (r_1,s_1^+) \rightarrow (r',s') $ where the states in this sequence of transitions have the same ``cleaved link" boundaries. By using the formula in lemma \ref {lemma:cancel}, we have:
$$
\langle \delta_n(r,s),(r',s') \rangle=\sum \limits _{\lk I-transition} (1/\righty {w_C})I_{\partial(r,s)}
$$  
Also, we see that $(r_1,s_1^-)$ is obtained from $(r,s)$ by surgery along an active resolution bridge $\gamma_1$ of $r$. Since $(r,s) \in \stn {\righty {T}}$ and a new free circle $C$ is created in $(r_1,s_1^-)$, $\gamma_1$ has to divide a cleaved circle $C_1$ of $r$ into $C$ and another cleaved circle $C_2$ of $r_1$. Furthermore, $(r',s')$ is obtained from $(r_1,s_1^+)$ by  surgering along  an active resolution bridge $\gamma_2$ of $r_1$. Since $(r',s') \in \stn {\righty {T}}$, $\gamma_2$ merges $C$ to a cleaved circle $D$ of $r_1$. If $C_2\equiv D$, there are two $\lk I$-transitions from $(r,s)$ to $(r',s')$. If $C_2\neq D$, there is only one $\lk I$-transition from $(r,s)$ to $(r',s')$. Therefore, $\delta_n$ is $\rightDelta {T,n}$. $\Diamond$    
\section{Invariance of the type $D$ structure under the weight moves} \label{sect:movingweight}
\noindent Following \cite{NR} and \cite{TJ}, we will prove that the homotopy type of $(\rightcomplex {T}, \rightTde T)$ is invariant under the following weight moves. Let $\righty T$ and $\righty T'$ be weighted tangle diagrams of a tangle $\righty {\lk T}$ with weighted arcs before and after the weight $w$ is moved along the crossing $c$ as the following figures:\\
$$
\inlinediag{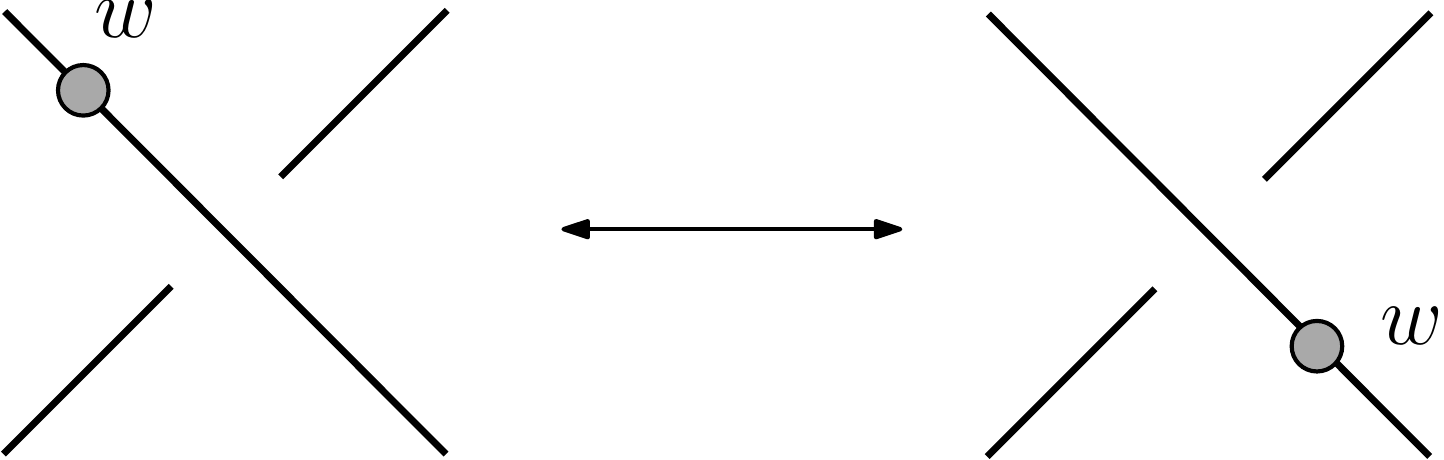} \hspace{0.5in} \inlinediag{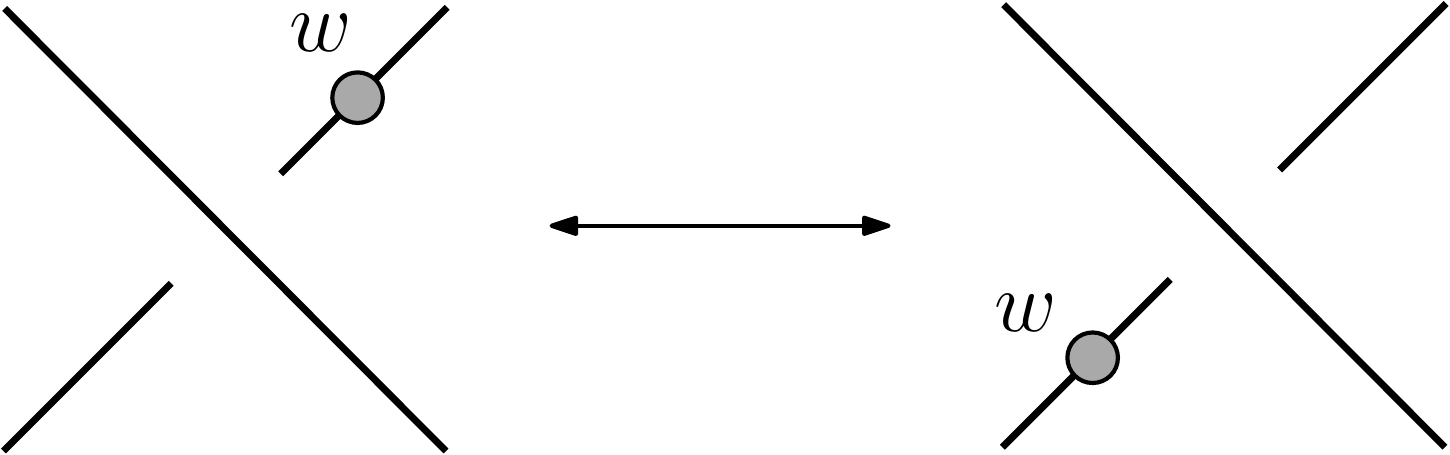} 
$$

\ \\
\noindent  We need to show that $(\rightcomplex {T}, \rightTde T)$ is homotopy equivalent to $(\rightcomplex {T'}, \rightTde {T'})$.\\
\ \\
 Let $\overline {D_c}: \dcomplex {T}{T'}$ be the $\field {\righty {T}}$-linear map defined as follows. Let $\xi =(r,s)$ be a generator of $\rightcomplex T$. $\overline {D_c}(r,s)$ is defined to be either $0$ if $r(c)=0$ or the sum of element(s) $e \otimes (r',s')$ if $r(c)=1$ where:
\begin {enumerate}
\item $r'$ is the resolution obtained from $r$ by surgering along the inactive bridge $\gamma$ at $c$ of $r$. 
\item $s'(D)=s(D)$ for all circles $D$ not abutting $c$ and the signs on circles abutting $c$ are computed by using $Khovanov$ sign rules. 
\item $e$ is the element in $\mathcal{B}\Gamma_{n}$ whose source is ${\partial (r,s)}$ and target is ${\partial (r',s')}$. We note that since $r'$ is obtained from $r$ by surgering a bridge on the right side, $e$ is either a right bridge element, or a right decoration element, or an idempotent.
\end {enumerate}
\noindent Additionally, if $e$ is a bridge element in $\mathcal{B}\Gamma_{n}$, we have:\\
$$\ze(r',s')=h(r',s')-q(r',s')/2=[h(r,s)-1]-[q(r,s)-3/2]/2=\ze(r,s)-1/4$$
Therefore, in this case, $\ze(r',s')+\ze(e)=\ze (r,s)$. By a similar computation, if $e$ is either an idempotent or a right decoration element, then $\ze(r',s')+\ze(e)=\ze (r,s)$. As a result, $\overline {D_c}$ is a $\ze$-grading preserving map.
\ \\
\begin {prop}\label{prop:weigh}The map $\Psi: \dcomplex {T}{T'}$ where $\Psi(r,s)= I_{\partial (r,s)} \otimes (r,s) + w.\overline {D_c}(r,s)$ is a type $D$ homomorphism.
\end {prop}  
\noindent $\bf Proof$. To ease the notation, we let $\de$, $\de'$ stand for  $\rightTde {T}$, $\rightTde {T'}$ respectively.\\
\ \\
It suffices to prove:
\begin{equation}\label{eqn:ident}
\weigh {\de'}\Psi+\weigh \Psi\de+(d_{\Gamma_{n}} \otimes \I_d)\Psi=0
\end {equation}
when applied to each $(r,s) \in \state {\righty {T}}$.  
\ \\
\ \\
Since the image of $\Psi$ does not have any term of the form $\lefty{e_C} \otimes(r_1,s_1)$, the last term of (\ref{eqn:ident}) will be $0$. Also, the map $\de$ can be written as the sum: $\de=D_T+D_{T,c}+E_T+E_{T,c}$ where $D_T$ and $D_{T,c}$ are the terms in the image of $\de$ obtained by surgering either one left bridge or one active bridge at $c_1 \ne c$, and $c$ respectively. $E_T$ ($E_{T,c}$ respectively) is the term in the image of $\de$ obtained by changing the decoration on a circle not abutting $c$ (abutting $c$) from $+$ to $-$.  We can write down a similar sum for $\de'$: $\de'=D_T'+D_{T',c}+E_{T'}+E_{T',c}$.  
\ \\
\ \\
As we can canonically identify $\rightcomplex T$ as $\rightcomplex {T'}$, let $\de_1$ be a $D$ structure on $\rightcomplex T$, which is precisely the same as $\de'$ on $\rightcomplex {T'}$. Note that: $\de+\de_1= E_{T,c}+E_{T',c}$. 
\ \\
\ \\
Let $\I_T: \dcomplex {T}{T'}$ be the map, defined by $\I_T(r,s)= I_{\partial (r,s)} \otimes (r,s)$. The left hand side of (\ref{eqn:ident}) can be written as the sum:
\ \\

$\weigh {\de'}{\I_T}+\weigh {\I_T}{\de_1}+ \weigh {\de'}{w\cdot\overline {D_c}}$
\ \\
\ \\
\indent $+\weigh {w \cdot\overline{D_c}} {\de}+\weigh {\I_T} {(E_{T,c}+E_{T',c})}$ 

\ \\
The sum of the first two terms will be $0$ since both are equal $\de'$. Rearranging the other terms, we need to prove:
\begin {equation} \label{eqn:ident1}
\weigh {\de'}{w\cdot\overline {D_c}}+\weigh {w\cdot\overline {D_c}}{\de}=\weigh {\I_T}{(E_{T,c}+E_{T',c})}
\end {equation}
By using the decomposition of $\de$ and $\de'$, we will in turn prove the following equations:
\begin {equation} \label{eqn:ident2}
\weigh {D_{T'}}{w \cdot\overline {D_c}}+\weigh {w\cdot \overline {D_c}}{D_T}=0
\end {equation}
\begin {equation} \label{eqn:ident3}
\weigh {E_{T'}}{w\cdot\overline {D_c}}+\weigh {w \cdot\overline {D_c}}{E_T}=0
\end {equation}
\begin {equation} \label{eqn:ident4}
\weigh {{E}_{T',c}}{w\cdot \overline {D_c}}+\weigh {w \cdot\overline {D_c}}{E_{T,c}}=0
\end {equation}
\begin {equation} \label{eqn:ident5}
\weigh {D_{T',c}}{w\cdot\overline {D_c}}+\weigh {w \cdot\overline {D_c}}{D_{T,c}}=\weigh {\I_T}{(E_{T,c}+E_{T',c})}
\end {equation}
The proof of (\ref{eqn:ident2}) comes from the proof that $\rightDelta {T}$ is a type $D$ structure on $\rightcomplex {T}$ (for the mirror crossing at $c$) and (\ref{eqn:ident3}) is one of case of proposition \ref {prop:11} (for the mirror crossing at $c$). For (\ref {eqn:ident4}), we remark that both terms will be $0$ if $r(c)=0$ because $E_{T,c}$ does not change the value of $r(c)$. $\weigh {E_{T',c}}{w\cdot \overline {D_c}}$ will be $0$ unless there is a $+$ circle abutting $c$ in the image of $\overline {D_c}$. As a result, $\weigh {E_{T',c}}{w\cdot \overline {D_c}}$ will be $0$ if one of the circle(s) abutting $c$ contains the marked point. In this case, it is also true that $\weigh {w\cdot\overline {D_c}}{E_{T,c}}=0$. Since $\overline {D_c}$, $E_{T,c}$ and $E_{T',c}$ only change the decorations on the circles abutting c, it suffices to prove (\ref {eqn:ident4}) by checking all of the possibilities for the circles abutting $c$. We can also assume that the circles abutting $c$ do not contain the marked point. Therefore, we have the following cases to check: 
\begin {enumerate}
\item All of circles abutting $c$ are free circles, 
\item There is at least one free circle and one cleaved circle
\item All of circles abutting $c$ are cleaved circles. 
\end {enumerate}
\noindent Case $(1)$ can be done similarly to case $(2)$ and has already been described in \cite{NR}.\\
\ \\
\noindent For case $(2)$, we only address the merging case (the argument for division is similar). Let $\pm_c$ and $\pm_f$ denote the decoration on the cleaved and free circles, respectively. Then:\\
\ \\

\noindent $\weigh {{E}_{T',c}}{\overline {D_c}}(+_c+_f) +\weigh {\overline {D_c}}{E_{T,c}}(+_c+_f)$\\

\indent $=\weigh {{E}_{T',c}}{({I_{\partial(r,s)} \otimes +_c})}+\weigh {\overline {D_c}}{\big [(\lefty{e_C}+\righty {w_C} \righty {e_C}) \otimes -_c+_f}$ 
\ \\
\ \\
\indent $\hspace {0.12in}$ $+(\righty {w_f} I_{\partial(r,s)} \otimes +_c-_f)\big]$
\ \\
\ \\
\indent $=(\mu_{\mathcal{B}\Gamma_{n}} \otimes \ensuremath{\mathbb{I}}_d) \big [I_{\partial(r,s)} \otimes  (\lefty {e_C}+ (\righty {w_C}+\righty {w_f})\righty {e_C})\otimes -_c \big ]+(\mu_{\mathcal{B}\Gamma_{n}} \otimes \ensuremath{\mathbb{I}}_d) \big[(\lefty {e_C}+\righty {w_C} \righty {e_C})\otimes$ 
\ \\
\ \\
\indent  $\hspace {0.12in}$ $I_{\partial(r,s)} \otimes-_c +\righty {w_f}I_{\partial(r,s)}  \otimes \righty {e_C}\otimes -_c\big ]$
\ \\
\ \\
\indent $=(\lefty {e_C}+ (\righty {w_C}+\righty {w_f})\righty {e_C})\otimes -_c+(\lefty {e_C}+\righty {w_C} \righty {e_C}) \otimes-_c +\righty {w_f} \righty {e_C}\otimes -_c$
\ \\
\ \\
\indent $=0$

\ \\
If either $s_c=-$ or $s_f=-$, then both terms will be $0$ and, thus, (\ref{eqn:ident4}) is true.
\ \\
\ \\
\noindent For case $(3)$, we will again only prove the identity for the merging case. Let $\pm_1$, $\pm_2$ be the decorations on cleaved circles $C_1$, $C_2$ respectively and $\pm_c$ be the decoration on a merged cleaved circle $C$. Let $e_{\gamma}$, $e_{\gamma_1}$ and $e_{\gamma_2}$ be the bridge elements representing the changes of the cleaved links: $+_1+_2 \longrightarrow +_c$, $-_1+_2 \longrightarrow -_c$ and $+_1-_2 \longrightarrow -_c$ respectively. Then\\
\ \\
$\weigh {{E}_{T',c}}{\overline {D_c}}(+_1+_2) +\weigh {\overline {D_c}}{E_{T,c}}(+_1+_2)$
\ \\
\ \\
\indent $=\weigh {{E}_{T',c}}{(e_{\gamma} \otimes +_c)} +\weigh {\overline {D_c}}{\big[ \big((\lefty {e_{C_1}}+\righty {w_{C_1}} \righty {e_{C_1}}) \otimes -_1+_2 \big)}$
\ \\
\ \\
\indent  $\hspace {0.12in}$ $+\big ((\lefty {e_{C_2}}+\righty {w_{C_2}} \righty {e_{C_2}}) \otimes +_1-_2 \big) \big]$
\ \\
\ \\
\indent $=e_{\gamma}(\lefty {e_{C}}+\righty {w_{C}} \righty {e_{C}}) \otimes -_c +  (\lefty {e_{C_1}}+\righty {w_{C_1}} \righty {e_{C_1}})e_{\gamma_1} \otimes -_c + (\lefty {e_{C_2}}+\righty {w_{C_2}} \righty {e_{C_2}})e_{\gamma_2} \otimes -_c$
\ \\
\ \\
\indent$=\big(e_{\gamma} \lefty {e_C} + \lefty {e_{C_1}} e_{\gamma_1}+  \lefty {e_{C_2}} e_{\gamma_2} \big) \otimes -_c + \righty {w_{C_1}}(e_{\gamma} \righty {e_C}+\righty {e_{C_1}} e_{\gamma_1}) \otimes -_c +\righty {w_{C_2}}(e_{\gamma} \righty {e_C}+\righty {e_{C_2}} e_{\gamma_2}) \otimes -_c $
\ \\
\ \\
\indent $=0$                                           
\ \\
since the first, second and third terms equal $0$ by relations \ref{rel:lefty1} and \ref {rel:4} defining $\mathcal{B}\Gamma_{n}$.
\ \\
Finally, we prove that (\ref{eqn:ident5}) is true. Again, since the proof for divisions is similar, we only present the proof for the merging case. We have the following subcases: 
\begin {enumerate}
\item All of circles are free: The proof of (\ref{eqn:ident5}) is then similar to the case $(2)$ below and already described in \cite {NR}.
\item There is at least one free circle $F$ and one cleaved circle $C$. Since the case where $C$ contains the marked point is the special case of what we are about to prove, we can assume $C$ does not contain the marked point.  Let $\pm_c$ and $\pm_f$ be the decorations on the cleaved and free circles respectively. We rewrite the map on the RHS of (\ref{eqn:ident5}) as
$$
\begin{array}{lcl}
+_c+_f  & \longrightarrow & {w.\big[(\righty {e_C}}\otimes -_c+_f )+({I_{\partial (r,s_C)}}\otimes +_c-_f)\big]\\
\ \\
-_c+_f  & \longrightarrow & {w.I_{\partial (r,s_C)}}\otimes-_c-_f\\
\ \\
+_c-_f &\longrightarrow & {w.\righty {e_C}}\otimes -_c-_f\\
\end{array}
$$ 
if $r(c)=1$ or $+_c \xrightarrow {0} -_c$ if $r(c)=0$. 
\ \\
\ \\
If $r(c)=0$, the first term of the LHS of (\ref{eqn:ident5}) equals $0$ because $\overline {D_c}$ is supported on states $(r',s')$ with $r'(c)=1$. The second term of the LHS maps $+_c$ to $(\righty {e_C}+\righty {e_C})\otimes -_c=0$. Therefore, (\ref{eqn:ident5}) is also true in this case.
\ \\
\ \\
If $r(c)=1$, the second term of the LHS of (\ref{eqn:ident5}) is $0$ and the map of the first term can be described as
$$
\begin{array}{lcl}
+_c+_f & \longrightarrow & (w.I_{\partial({r,s})}\righty {e_C}\otimes -
_c+_f)+(w.I_{\partial({r,s})}I_{\partial({r,s})}\otimes +_c-_f) \\ 
\ \\
-_c+_f  & \longrightarrow & w.I_{\partial({r,s_C})}I_{\partial({r,s_C})} \otimes -_c-_f \\
\ \\
+_c-_f & \longrightarrow & w.\righty {e_C} I_{\partial ({r,s_C})}\otimes -_c-_f \\
\end{array}
$$ 
\begin{center}
\begin{figure}
\includegraphics[scale=0.5]{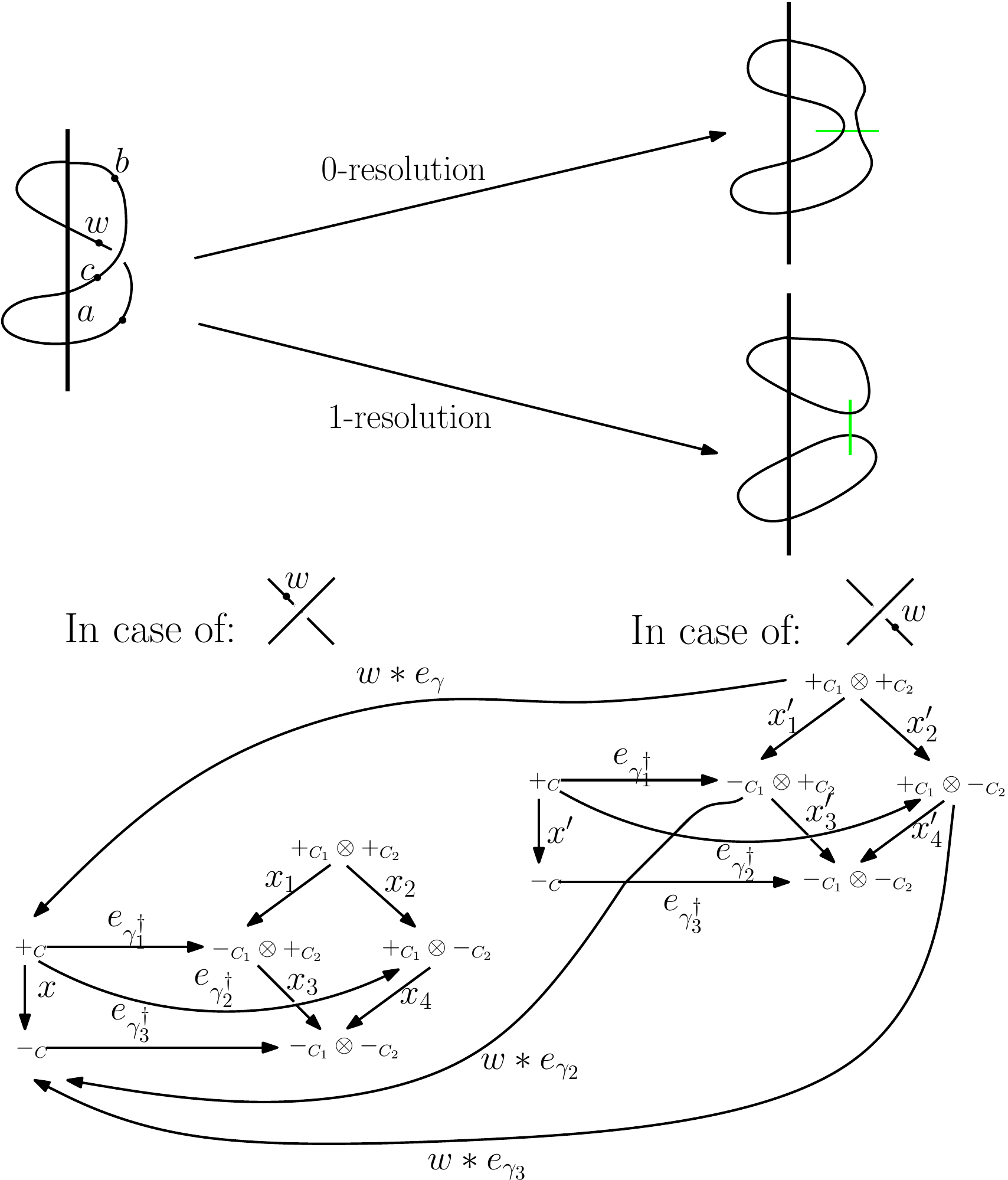}
\caption{An illustration of how to construct a type $D$ morphism between the two type $D$ structures before and after moving weight. The recorded algebra elements above the arrows correspond to the change in the cleaved links. The diagram does not cover all the terms of the type $D$ structures or the morphism. It only illustrates the last case of proposition \ref {prop:weigh}. However, this picture can be modified to give the pictures for the other cases by changing $e_{\gamma_i}$ and $e_{\gamma_i}^{\dagger}$ to suitable algebra elements.}
\label{fig:movingweight}
\end{figure}
\end{center}
By comparing the coefficients on each generator, the images of the LHS and the RHS agree. As a result, (\ref{eqn:ident5}) is true in this case.
\item All of the circles are cleaved circles: The proof will be similar to the above cases when we can prove that two sides of (\ref {eqn:ident5}) agree on every generator $\xi=(r,s)$ of $\rightcomplex T$ by using the relation $\righty e_{(\ga,\si,\si')} \righty e_{(\ga^\dagger,\si',\si_C)}=\righty {e_C}$. Figure \ref {fig:movingweight} illustrates the proof for this case. In this figure, $C$ is the cleaved circle of the $0$-resolution while $C_1$ and $C_2$ correspond to the top and the bottom cleaved circles of the $1$-resolution. In addition, $x_1=\lefty {e_{C_1}}+(w+b)\righty {e_{C_1}}$ and $x_2=\lefty {e_{C_2}}+(c+a)\righty {e_{C_2}}$ while $x'_1=\lefty {e_{C_1}}+b\righty {e_{C_1}}$ and $x'_2=\lefty {e_{C_2}}+(w+a+c)\righty {e_{C_2}}$. In the identity (\ref{eqn:ident5}):

$$
{\text {LHS}}(+_{c_1}+_{c_2})=w.e_{\gamma}e_{\gamma_1^{\dagger}}\otimes -_{c_1}+_{c_2}\oplus w.e_{\gamma}e_{\gamma_2^{\dagger}}\otimes +_{c_1}-_{c_2}
$$
and:
$$
{\text {RHS}}(+_{c_1}+_{c_2})=(x_1+x'_1)\otimes -_{c_1}+_{c_2}\oplus(x_2+x'_2)\otimes -_{c_2}+_{c_1}
$$
Using relation \ref {rel:special} in the set of other bridge relations, we see that ${\text {LHS}}(+_{c_1}+_{c_2})={\text {RHS}}(+_{c_1}+_{c_2})$. Similarly, we can prove that (\ref {eqn:ident5}) is true when applied for $-_{c_1}+_{c_2}$ or $+_{c_1}-_{c_2}$. 
Therefore, $\Psi$ is a type $D$ homomorphism.  $\Diamond$ 
\end {enumerate}

\begin {prop} $(\rightcomplex {T}, \rightTde T)$ is isomorphic to $(\rightcomplex {T'}, \rightTde {T'})$ as type $D$ structures.
\end {prop}
\noindent $\bf Proof$. Let $\Phi: \dcomplex {T'}{T}$ where $\Phi(r,s)= I_{\partial (r,s)} \otimes (r,s) + w.\underline {D_c}(r,s)$ where $\underline {D_c}$ is defined identically as $\overline {D_c}$ but from $\righty {T'}$ to $\righty {T}$. By proposition \ref {prop:weigh}, $\Phi$ is a $D$ homomorphism. We will prove: $\Phi \ast \Psi=I_{\rightcomplex {T}}$ and $\Psi \ast \Phi =I_{\rightcomplex {T'}}$ where $\ast$ stands for the composition of two type $D$ structures (described in section \ref {sect:stable}). After that, we can conclude $(\rightcomplex {T}, \rightTde T)$ is isomorphic to $(\rightcomplex {T'}, \rightTde {T'})$. Indeed,\\
\ \\
$\Phi \ast \Psi $
\ \\
\ \\
\indent $=\weigh {\Phi}{\Psi}$
\ \\
\ \\
\indent$=\weigh {\I_{T'}}{\I_T}+w.\big[\weigh {\underline{D_c}}{\I_T}+\weigh {\I_{T'}}{\overline {D_c}}\big]$ 
\ \\
\ \\
\indent $\hspace {0.12in}$ $+\weigh {\underline {D_c}}{\overline {D_c}}$ \ \\

\noindent We note that the last term is $0$ since $\underline {D_c}$ is supported on states $(r,s)$ where $r(c)=1$ while the image of $\overline {D_c}$ contains the states $(r_1,s_1)$ where $r_1(c)=0$. The sum of the two middle terms is also $0$ since the second term is equal to the third term. Moreover, $\weigh {\I_{T'}}{\I_T}=I_{\rightcomplex {T}}$. Combining all these facts, we obtain: $\Phi \ast \Psi=I_{\rightcomplex {T}}$. Similarly,  $\Psi \ast \Phi =I_{\rightcomplex {T'}}$ by an almost identical argument. $\Diamond$

\section{Graded differential algebra and stable homotopy} \label {sect:stable}
\subsection{Preliminary}
\noindent The different projections of the right tangle $\overrightarrow {\lk T}$ can have a different number of arcs and thus, the corresponding type $D$ structures will be vector spaces over different base fields. To relate these structures, we will need an appropriate algebraic tool: the stable homotopy, whose construction is based on the idea in section $4$ of \cite {LR3}. Let $\F$ be a field and $W$ be a vector space over $\F$. Let $\F_W$ be the field of rational function of $P_W$ where $P_W$ is the symmetric algebra of $W$. Recall: let $W$, $W'$ be two vector spaces over $\F$ and let $M$, $M'$ be two vector spaces over $\F_W$, $\F_{W'}$ respectively. A pair $(M, W)$ is stably isomorphic to $(M', W')$ if there is a triple $(W'',i,i')$ where $W''$ is a vector space over $\F$, $W \xrightarrow {i} W''$ and $W \xrightarrow {i'} W''$ are injective linear maps, such that $M \otimes_{\F_W} \F_{W''} \cong M' \otimes_{\F_{W'}} \F_{W''}$ as vector spaces over $\F_{W''}$. This relation is proved to be an equivalence relation in \cite {LR3}. We always assume that $\F$ is $\Z_2$ and $W$ is a vector space over $\Z_2$, unless otherwise stated. In this section, we give a modification of this definition which will allow us to relate two type $D$ (or $A$) structures over two different fields. Furthermore, the complexes obtained by gluing stable homotopy equivalences of type $A$ structures to stable homotopy equivalences of type $D$ structures will be stably homotopic in the sense of \cite {LR3}. 
\ \\
\subsection{$ {A_{W,\bf{\lk I}}}$ category and stable $ A_{\infty}$ homotopy equivalence} \label {sect:stable1}
\ \\
\noindent We will describe a way to construct $A_{\infty}$ modules over the ground field $\F_W$ from the one over $\Z_2$ by first upgrading our unital DGA $(\mathcal{B}\Gamma_{n},d_{\Gamma_{n}})$ (over $\Z_2$) to one over $\F_W$.\\ 
\ \\
\noindent We will let $(A,\lk I)$, $\mu_1$ and $\mu_2$ stand for $(\mathcal{B}\Gamma_{n},\lk I_n)$, $d_{\Gamma_{n}}$ and the product on $\mathcal{B}\Gamma_{n}$ respectively. \\ 
\begin {defn}
For each vector space $W$ over $\Z_2$, let $A_W:=A \otimes_{\Z_2} \F_W$, equipped with the following $\F_W-$linear maps
$$
\mu_{W,1}: A_W \longrightarrow A_W[-1]
$$
$$
\mu_{W,2}: A_W \otimes_{\F_W} A_W \longrightarrow A_W
$$
where $\mu_{W,1}=\mu_1 \otimes \I_{\F_W}$ and $\mu_{W,2}=\mu_2 \otimes \I_{\F_W}$, under the canonical isomorphism: $A_W \otimes_{\F_W} A_W \cong A^{\otimes 2} \otimes_{\Z_2}\F_W$. We also let $\lk I_W$ be $\lk I \otimes_{\Z_2} \F_W$ 
\end {defn}

\noindent Using the fact that tensor product is a functor, it is straightforward to verify the following proposition:\\
\begin {prop}
If $(A,\mu_1,\mu_2)$ is unital DGA over $\Z_2$ then $(A_W,\mu_{W,1},\mu_{W,2})$ is unital DGA over $\F_W$. 
\end {prop}
\noindent The proof of this proposition is left to the reader.\\

\noindent We recall the definition of an $A_{\infty}$ module, as in \cite {LOT}.\\

\begin {defn}\cite {LOT}
$(M,\{m_i\}_{i \in \N})$ is a right $A_{\infty}$ module over $(A_W,\lk I_W)$ if $M$ is a graded $\lk I_W$-module and for each $i \in \N$, $m_i: M\otimes_{\lk I_W} A_{W}^{\otimes (i-1)} \longrightarrow M[i-2]$ is an $\lk I_W$-linear map, which satisfies the following compatibility condition:
$$
0=\sum \limits_{i+j=n+1} m_i(m_j\otimes \I^{\otimes (i-1)})+\sum\limits_{i+j=n+1,j<3,k>0}m_i(\I^{\otimes k}\otimes \mu_{W,j} \otimes\I^{\otimes(i-k-1)})
$$
\noindent M is said to be strictly unital if for any $\xi \in M$, $m_2(\xi\otimes 1_{A_W})=\xi$, but for $n>2$, $m_n(\xi\otimes a_1\otimes...\otimes a_{n-1})=0$ if any $a_i=1_{A_W}$\\
\begin {defn} \cite {LOT}
Let $(M,\{m_i\})$ and $(M',\{m'_i\})$ be two $A_{\infty}$ modules over $(A_W,\lk I_W)$. An $A_{\infty}$ homomorphism $\Psi$ is a collection of maps:
$$
\psi_i: M\otimes_{\lk I_W} A_W^{\otimes(i-1)}\longrightarrow M'[i-1]
$$ 
indexed by $i \in \N$, satisfying the compatibility conditions:
$$
0=\sum\limits_{i+j=n+1}m'_i(\psi_j\otimes \I^{\otimes(i-1)})+\sum\limits_{i+j=n+1}\psi_i(m_j\otimes \I^{\otimes(i-1})
+\sum\limits_{i+j=n+1,k>0,j<3}\psi_i(\I^{\otimes k}\otimes \mu_{W,j}\otimes \I^{\otimes(i-k-1)})
$$
Additionally, a homotopy $H$ between two $A_{\infty}$ morphisms $\Psi$ and $\Phi$ is a set of maps $\{h_i\} $ with $h_i:M\otimes_{\lk I_W}A_W^{\otimes (i-1)}\longrightarrow M'[i]$ such that:
\begin {equation}
\psi_n+\phi_n=\sum\limits_{i+j=n+1}m'_i(h_j\otimes \I^{\otimes(i-1)})+\sum\limits_{i+j=n+1}h_i(m_j\otimes \I^{\otimes(i-1)})+\sum\limits_{i+j=n+1,k>0,j<3}h_i(\I^{\otimes k} \otimes \mu_{W,j} \otimes \I^{\otimes(i-k-1)})
\end {equation}
\end {defn}
\end {defn}
\noindent For an ease of notation, we sometimes let $M$ stand for $(M,\{m_i\})$ when it is clear from the context.\\
\begin {prop} \cite {LOT}
Let $A_{W,\lk I}$ be a collection of $A_{\infty}$ modules over $(A_W,\lk I_W)$. $Mor(M,M')$ is a collection of $A_{\infty}$ homomorphisms $(\Psi,\{\psi_i\})$ and the composition of $\Phi\ast \Psi$ is the set of maps:
$$
(\Phi\ast \Psi)_n=\sum \limits_{i+j=n+1}\phi_i(\psi_j \otimes \I^{\otimes (i-1)})
$$
for each $n \in \N$. Then $A_{W,\lk I}$ forms a category.    
\end {prop}
\begin {defn} \label {defn:functor}
Let $\varphi:W \hookrightarrow W'$ be a linear injection. Let $(M,\{m_i\}_{i \in \N})$  be an object of $A_{W,\lk I}$ and $(\Psi,\{\psi\}_{i\in \N}) \in Mor((M,\{m_i\}),(M',\{m'_i\}))$ be a morphism of $A_{W,\lk I}$. We define:
\begin{enumerate}
\item 
$$
\lk F_{\varphi}(M)=M \otimes_{\varphi} \F_{W'}
$$
\item For each $i\in \N$, $\lk F_{\vr}(m_i): \lk F_{\vr}(M)\otimes_{\lk I_{W'}}A_{W'}^{\otimes(i-1)} \longrightarrow \lk F_{\vr}(M)[i-2]$ is defined by:
$$
\lk F_{\varphi}(m_i)\big( (m\otimes r)\otimes (a_1\otimes r_1)\otimes....\otimes (a_{i-1}\otimes r_{i-1})\big)=m_i(m\otimes a_1\otimes...\otimes a_{i-1})\otimes rr_1...r_{i-1}
$$
where $m\in M$, $a_j \in A_W$ and $r$, $r_j\in \F_{W'}$ for $j=1,...,i-1$.
\item Similarly, for each $i\in \N$, $\lk F_{\vr}(\psi_i): \lk F_{\vr}(M)\otimes_{\lk I_{W'}}A_{W'}^{\otimes(i-1)} \longrightarrow \lk F_{\vr}(M')[i-1]$ is defined in the same manner as $\lk F_{\vr}(m_i)$  
\end{enumerate}
\end {defn}
\noindent Combining the facts that tensor product is a functor and there is the following canonical isomorphism:\\
$$
(M\otimes_{\varphi} \F_{W'})\otimes_{\lk I_{W'}} (A_W\otimes_{\varphi} \F_{W'})^{\otimes (i-1)}\cong (M\otimes_{\lk I_W}A_W^{\otimes(i-1)})\otimes_{\varphi}\F_{W'}
$$
for each $i \in \N$, it is straightforward to prove the following two propositions:
\begin {prop} \label {prop:12}
$\lk F_{\varphi}(M,\{m_i\})=(\lk F_{\varphi}(M),\{\lk F_{\varphi}(m_i)\})$ is well-defined and it is an object of $A_{W',\lk I}$. Likewise, $\lk F_{\varphi}(\Psi)$ is well-defined and belongs to $Mor(\lk F_{\varphi}(M,\{m_i\}), \lk F_{\varphi}(N,\{n_i\}))$
\end {prop}

\begin {prop} \label {prop:functor1}
For each injection $\varphi:W \hookrightarrow W'$, there is a functor $\lk F_{\varphi}:A_{W,\lk I} \longrightarrow A_{W',\lk I}$ defined by
$$
\lk F_{\varphi}(M,\{m_i\})=( \lk F_{\varphi}(M), \{\lk F_{\varphi}(m_i)\})
$$
$$
\lk F_{\varphi}(\Psi,\{\psi_i\})= (\lk F_{\varphi}(\Psi), \{\lk F_{\varphi}(\psi_i)\})
$$
Furthermore, if $\Psi$ and $\Phi$ are $A_{\infty}$ homotopic in $A_{W,\lk I}$ then $\lk F_{\varphi}(\Psi)$ is $A_{\infty}$ homotopic to $\lk F_{\varphi}(\Phi)$ in $A_{W',\lk I}$. 
\end {prop}

\noindent We now have enough tools to relate two $A_{\infty}$ structures over different fields by defining the stable homotopy equivalence of $A_\infty$ modules:\\
\begin {defn} \label {defn:AIrelate}
Let $W$ and $W'$ be two vector spaces over $\Z_2$. Let $(M,\{m_i\})$, $(M',\{m'_i\})$ be objects of $A_{W,\lk I}$ and $A_{W',\lk I}$ respectively. Then $(M,\{m_i\})$ is stably homotopy equivalent to $(M',\{m'_i\})$ if there is a triple $(\varphi,\varphi', W'')$ where $W''$ is a vector space over $\Z_2$ and $\varphi:W \hookrightarrow W''$ and $\varphi':W' \hookrightarrow W''$ are linear injections, such that $(\lk F_{\varphi}(M),\{\lk F_{\vr}(m_i)\})$ is homotopy equivalent to $(\lk F_{\varphi'}(M'),\{\lk F_{\vr}(m'_i)\})$ in the category $A_{W'',\lk I}$.
\end{defn}
\begin {prop} \label {prop:Aequivalent} 
Stable homotopy equivalence of $A_{\infty}$ modules is an equivalence relation.
\end {prop}
\noindent The proof of this proposition follows directly from lemma $4.3$ of \cite {LR3}, the fact that homotopy equivalence of $A_{\infty}$ modules is an equivalence relation, and the following lemma:\\
\begin {lemma} \label{lemma:Ainject}
Let $\widetilde W$ be a vector space over $\Z_2$ and $\widetilde {\varphi}: W'' \hookrightarrow \widetilde W$ be a linear injection. If $(M,\{m_i\})$ is stably homotopy equivalent to $(M',\{m'_i\})$ via $(\varphi,\varphi', W'')$, then $(M,\{m_i\})$ is stably homotopy equivalent to $(M',\{m'_i\})$ via $(\widetilde {\varphi} \circ \varphi, \widetilde {\varphi} \circ \varphi', \widetilde W)$ 
\end {lemma}
\noindent $\bf Proof.$ Since $(\lk F_{\varphi}(M),\{\lk F_{\vr}(m_i)\})$ is homotopy equivalent to $(\lk F_{\varphi'}(M'),\{\lk F_{\vr'}(m'_i)\})$, $(\lk F_{\widetilde {\varphi}}(\lk F_{\varphi}(M),\{\lk F_{\widetilde {\varphi}}(\lk F_{\vr}(m_i))\})$ is homotopy equivalent to $(\lk F_{\widetilde {\varphi}}(\lk F_{\varphi'}(M')),\{\lk F_{\widetilde {\varphi}}(\lk F_{\vr'}(m'_i))\})$ by proposition \ref {prop:functor1}. Thus, we will finish the proof of this lemma if we can prove 
$$
(\lk F_{\widetilde{\varphi} \circ \varphi}(M),\{\lk F_{\widetilde{\varphi} \circ \varphi}(m_i)\})  \cong_{A_{\widetilde W, \lk I}}  (\lk F_{\widetilde{\varphi}} \circ \lk F_{\varphi}(M),\{\lk F_{\widetilde {\varphi}}(\lk F_{\vr}(m_i))\})
$$ 
and
$$
(\lk F_{\widetilde{\varphi} \circ \varphi'}(M'),\{\lk F_{\widetilde{\varphi} \circ \varphi'}(m'_i)\})  \cong_{A_{\widetilde W, \lk I}}  (\lk F_{\widetilde{\varphi}} \circ \lk F_{\varphi'}(M'),\{\lk F_{\widetilde {\varphi}}(\lk F_{\vr'}(m'_i))\})
$$ 
\ \\
First of all, we have: 
$$
\lk F_{\widetilde{\varphi} \circ \varphi}(M)=M\otimes_{\widetilde{\varphi} \circ \varphi}\F_{\widetilde W}\cong (M\otimes_{\varphi}\F_{W''})\otimes_{\widetilde {\varphi}}\F_{\widetilde W}= \lk F_{\widetilde{\varphi}} \circ \lk F_{\varphi}(M)
$$ 
Secondly, under this identification of the underlying modules, we need to prove $\lk F_{\widetilde{\varphi} \circ \varphi}(m_i)= \lk F_{\widetilde{\varphi}} \circ \lk F_{\varphi}(m_i)$. Using the definition \ref {defn:functor}, we have:
$$
\lk F_{\widetilde{\varphi} \circ \varphi}(m_i) \big((m\otimes r)\otimes (a_1\otimes r_1)\otimes...\otimes (a_{i-1}\otimes r_{i-1})\big)=m_i(m \otimes a_1 \otimes...\otimes a_{i-1})\otimes _{\widetilde{\varphi} \circ \varphi} rr_1...r_{i-1}
$$ 
\noindent On the other hand, \\
\ \\
$
\lk F_{\widetilde{\varphi}} \circ \lk F_{\varphi}(m_i)\big(((m\otimes_{\varphi} 1)\otimes_{\widetilde {\varphi}} r)\otimes ((a_1\otimes_{\varphi} 1)\otimes_{\widetilde {\varphi}} r_1)\otimes...\otimes ((a_{i-1}\otimes_{\varphi} 1)\otimes_{\widetilde {\varphi}} r_{i-1})\big)
$
\ \\
\ \\
\indent $
=\lk F_{\varphi}(m_i)\big((m\otimes_{\varphi}1) \otimes (a_1\otimes_{\varphi}1)\otimes...\otimes (a_{i-1}\otimes 1)\big)\otimes_{\widetilde {\varphi}} rr_1...r_{i-1}
$
\ \\
\ \\
\indent $
=\big(m_i(m\otimes a_1\otimes...\otimes a_{i-1})\otimes_{\varphi}1\big)\otimes_{\widetilde {\varphi}}rr_1...r_{i-1}
$
\ \\
Therefore, $\lk F_{\widetilde{\varphi} \circ \varphi}(m_i)= \lk F_{\widetilde{\varphi}} \circ \lk F_{\varphi}(m_i)$.         $\Diamond$
\ \\

\subsection{${D_{W,\lk I}}$ Category and Stable $D$ homotopy equivalence}
\ \\
We first review the definition of the type $D$ structure and the $D_W$ category as in \cite {LOT}.
\begin {defn} \cite {LOT}
Let $N$ be a graded $\lk I_W$-module. A (left) $D$ structure on N is a linear map:
\begin {equation}
\delta: N \longrightarrow (A_W \otimes_{\lk I_W} N) [-1]
\end {equation}
satisfying:
$$
(\mu_{W,2}\otimes  \I_N)(\I_{A_W} \otimes \delta)\delta +(\mu_{W,1}\otimes \I_N)\delta=0
$$
      
\noindent A morphism of type $D$ structures $(N,\delta) \longrightarrow (N',\delta')$ is a map $\psi: N \longrightarrow A_W\otimes_{\lk I_W} N'$ such that:\\
$$
(\mu_{W,2} \otimes \I_N)(\I_{A_W}\otimes \delta')\psi+(\mu_{W,2}\otimes \I_N)(\I_{A_W}\otimes \psi)\delta+(\mu_{W,1}\otimes \I_N)\psi=0
$$ 
$H:N\longrightarrow (A_W\otimes_{\lk I_W} N')[1]$ is a homotopy of two type $D$ morphisms $\psi$ and $\phi$ if:
$$
\psi+\phi=(\mu_{W,2} \otimes \I_N)(\I_{A_W}\otimes \delta')H+(\mu_{W,2}\otimes \I_N)(\I_{A_W}\otimes H)\delta+(\mu_{W,1}\otimes \I_N)H
$$ 
\end {defn}
\begin {prop} \cite {LOT}
Let $D_{W,\lk I}$ be the collection of type $D$ structures $(N,\delta)$ over $(A_W,\lk I_W)$. Let $Mor((N,\delta),(N',\delta'))$ be the set of type $D$ morphisms and the composition $\phi \ast \psi$ of $\psi: N \longrightarrow A_W\otimes N'$ and $\phi: N' \longrightarrow A_W\otimes_{\lk I_W} N''$ is defined to be $(\mu_{W,2}\otimes \I_{N''})(\I_{A_W} \otimes \phi)\psi$. Additionally, the identity morphism at $(N,\delta)$ is $I_N:N \longrightarrow A_W\otimes_{\lk I_W} N$ given by $x\longrightarrow 1_{A_W}\otimes x$. Then $D_{W,\lk I}$ forms a category. 
\end {prop}

\noindent Using the same technique as in \ref {sect:stable1}, for each $\vr: W\hookrightarrow W'$ linear injection, we will construct a functor from $D_{W,\lk I}$ to $D_{W',\lk I}$. Before doing that, we need the following definition: \\
\begin {defn} \label {defn:34}
Let $\vr: W\hookrightarrow W'$ be a linear injection. Let $(N,\delta)$, $(N',\delta')$ be two objects of $D_{W,\lk I}$ and $\psi \in Mor((N,\delta),(N',\delta'))$. We define:
\begin {enumerate}
\item $$\lk G_{\vr}(N):=N\otimes_{\vr} \F_{W'}$$
\item Furthermore,
$$
\lk G_{\vr}(\delta): N\otimes_{\vr}\F_{W'} \longrightarrow A_{W'}\otimes_{\lk I_{W'}} (N\otimes_{\vr}\F_{W'})[-1]
$$
is defined to be $\delta \otimes \I_{\F_{W'}}$ under the canonical isomorphism $A_{W'}\otimes_{\lk I_{W'}} (N\otimes_{\vr}\F_{W'}) \cong (A_W \otimes_{\lk I_W} N)\otimes_{\vr}\F_{W'}$, and 
\item $$
\lk G_{\vr}(\psi): N\otimes_{\vr}\F_{W'} \longrightarrow A_{W'} \otimes_{\lk I_{W'}} (N'\otimes_{\vr}\F_{W'})
$$
is defined to be $\psi \otimes \I_{\F_{W'}}$ under the canonical isomorphism $A_{W'} \otimes_{\lk I_{W'}} (N'\otimes_{\vr}\F_{W'}) \cong (A_W\otimes_{\lk I_{W}} N')\otimes_{\vr}\F_{W'}$
\end {enumerate}
\end {defn}
\begin {prop} \label {prop:21}
Using the same notation as in the above definition, there exists a functor $\lk G_{\vr}: D_{W,\lk I}\longrightarrow D_{W',\lk I}$ defined by:
$$
\lk G_{\vr}(N,\delta)=(\lk G_{\vr}(N), \lk G_{\vr}(\delta))
$$
and $\lk G_{\vr}(\psi)$ is defined as above. Furthermore, let $\psi$ and $\phi$ belong to $Mor((N,\delta),(N',\delta'))$ in $D_{W,\lk I}$ and if $H$ is a homotopy from $\psi$ to $\phi$, then $\lk G_{\vr}(H)$ is a homotopy from $\lk G_{\vr}(\psi)$ to $\lk G_{\vr}(\phi)$ in $D_{W',\lk I}$ where:
$$
\lk G_{\vr}(H): N\otimes_{\vr}\F_{W'} \longrightarrow A_{W'}\otimes_{\F_{W'}} (N'\otimes_{\vr}\F_{W'})[1]
$$
is defined to be $H\otimes \I_{F_{W'}}$ under the isomorphism $A_{W'}\otimes_{\lk I_{W'}} (N'\otimes_{\vr}\F_{W'})\cong (A_W\otimes_{\lk I_W}N')\otimes_{\vr} \F_{W'}$
\ \\
\ \\
\noindent Therefore, $\lk G_{\vr}$ induces a functor from the homotopy category of $D_{W,\lk I}$ to the homotopy category of $D_{W',\lk I}$.\\
\end {prop}
\noindent We now can give a definition of stable $D$ homotopy.\\
\begin {defn}
Let $W$ and $W'$ be vector spaces over $\Z_2$. Let $(N,\delta)$, $(N',\delta')$ be objects of $D_{W,\lk I}$ and $D_{W',\lk I}$, respectively. Then $(N,\delta)$ is stably homotopy equivalent to $(N,\delta')$ if there is a triple $(\vr,\vr',W'')$ where $W''$ is a vector space over $\Z_2$, $\vr:W\hookrightarrow W''$ and $\vr:W'\hookrightarrow W''$ are linear injections, such that $\lk G_{\vr}(N,\delta)$ is homotopy equivalent to $\lk G_{\vr'}(N',\delta')$ in the category $D_{W'',\lk I}$
\end {defn}
\begin {prop} \label {prop:Dequivalent}
Stable homotopy of type $D$ structures is an equivalence relation.
\end {prop}
\noindent  Since the proofs of the propositions \ref {prop:21} and \ref {prop:Dequivalent} are similar to the proofs of propositions \ref {prop:functor1} and \ref {prop:Aequivalent} respectively, we leave them to the readers. We have the following remark about the property of the composition of two functors, which is useful for the next subsection. \\
\begin {remark} \label {remark:Dfunctor}
Let $\vr_1:W\hookrightarrow W_1$ and $\vr_2: W_1 \hookrightarrow W_2 $ be injective linear maps. Let $(N,\delta)$ be an object of $A_W$. Then $\lk G_{\vr_2 \circ \vr_1}(N,\delta) \cong_D \lk G_{\vr_2}(\lk G_{\vr_1})(N,\delta)$ 
\end {remark}   
\subsection{Pairing an $ A_{\infty}$ module and a type $ D$ structure over different DGAs} \label {sect:pairing} 
In \cite {LOT}, there is the result that we can pair an object $(M,\{m_i\})$ of $A_{W,\lk I}$ and an object $(N,\delta)$ of $D_{W,\lk I}$ to form a chain complex $(M\bo N, \by)$ over $\F_W$. Additionally, if $(M,\{m_i\})\simeq_{A_{\infty}}(M',\{m'_i\})$ in $A_{W,\lk I}$ category and $(N,\delta) \simeq_D (N',\delta')$ in $D_{W,\lk I}$ category, then $(M\bo N, \by)$ is chain homotopic to $(M'\bo N', \by)$. 
\ \\
\ \\
For our purpose, since we need to pair a type $A$ and a type $D$ structures over distinct differential graded algebras, we will modify the way to pair them to get a chain complex. Furthermore, we will prove that under the change of either type $A$ or type $D$ by a stable homotopy equivalence, the glued chain complexes are stably chain homotopic.
\ \\
\ \\
Let $W$ and $W'$ be two vector spaces over $\Z_2$. Let $(M,\{m_i\})$, $(N,\delta)$ be objects of $A_{W,\lk I}$ and $D_{W',\lk I}$ respectively. Let $W_1:=W\oplus W'$ and let $\vr:W\hookrightarrow W_1 $, $\vr':W'\hookrightarrow W_1 $ be two canonical linear injections. 
\begin {defn}
Define $M\bos N$ to be the graded vector space $\lk F_{\vr}(M)\otimes_{\lk I_{W_1}}\lk G_{\vr'}(N)$ over $\F_{W_1}$ and $\bys:M\bos N \longrightarrow (M\bos N)[-1]$ to be the map :
$$
\bys=\sum \limits_{k=0}^{\infty}\big (\lk F_{\vr}(m_{k+1})\otimes \I^{\otimes {(k+1)}}_{\lk G_{\vr'}(N)}\big )\circ \big (\I_{\lk F_{\vr}(M)} \otimes \Delta_k\big )
$$
where $\Delta_k: \lk G_{\vr'}(N) \longrightarrow A^{\otimes k}_{W_1} \otimes_{\lk I_{W_1}} \lk G_{\vr'}(N)[-k] $ is defined by induction $\Delta_0=\I_{\lk G_{\vr'}(N)}$, $\Delta_1=\lk G_{\vr'}(\delta)$ and the relation: $\Delta_n=(\I^{\otimes(n-1)} \otimes \lk G_{\vr'}(\delta))\Delta_{n-1}$.
\end {defn}
\noindent $\bf Note.$ $(M\bos N,\bys)$ is defined exactly the same as the definition of $(\lk F_{\vr}(M) \bo  \lk G_{\vr'}(N), \by)$ in \cite {LOT} and thus, it is a chain complex. \\
\begin {prop} \label {prop:lhomotopy}
Let $W$, $W'$ and $\ww$ be three vector spaces over $\Z_2$. Let $(N,\delta)$ be an object of $D_{\ww,\lk I}$. Let $(M,\{m_i\})$, $(M',\{m'_i\})$ be objects of $A_{W,\lk I}$ and $A_{W',\lk I}$ respectively, such that $(M,\{m_i\})$ is homotopy equivalent to $(M',\{m'_i\})$ via the triple $(\vr,\vr',W'')$. Then $(M\bos N,\bys)$ is stably chain homotopic to $(M'\bos N,\bys)$.
\end {prop}
\noindent Before proving this proposition, we state the following lemma which can be proved in a similar manner as lemma \ref {lemma:Ainject}.  \\
\begin {lemma} \label {lemma:chainrelate}
Let $(M,\{m_i\})$ and $(N,\delta)$ be objects of $A_{W,\lk I}$ and $D_{W,\lk I}$ categories respectively. Let $\vr: W \hookrightarrow W'$ be a linear injection. Then $((M\bo N)\otimes_{\vr} \F_{W'}, \by \otimes \I_{\F_{W'}})$ is chain isomorphic to $(\lk F_{\vr}(M) \bo \lk G_{\vr}(N), \by)$.
\end {lemma}

\noindent $\bf Proof$ $\bf of$ $\bf Proposition$ \ref {prop:lhomotopy} $\bf .$ \\
Define $V:=W \oplus \ww$, $V':=W' \oplus \ww$ and $V'':=W'' \oplus \ww$. Let $\w {\vr}: V \longrightarrow V''$ and $\w {\vr}': V' \longrightarrow V''$ denotes the maps $\vr \oplus \I_{\ww}$ and $\vr' \oplus \I_{\ww}$. Let $p:\ww \longrightarrow V$, $p':\ww \longrightarrow V'$, $\pi: W \longrightarrow V$, $\pi':W' \longrightarrow V'$ and $\pi'':W'' \longrightarrow V''$ be the natural injections. We immediately have the following relations:
\begin {enumerate}
\item $\w{\vr} \circ p= \w{\vr}'\circ p'$
\item $\w{\vr} \circ \pi =\pi'' \circ \vr$
\item $\w{\vr}' \circ \pi' =\pi'' \circ \vr'$
\end {enumerate}
\ \\
We will prove that $(M\bos N, \bys)$ is stably chain homotopic to $(M' \bos N,\bys)$ via the triple $(\w {\vr},\w {\vr}', V'')$. Indeed, using lemma \ref {lemma:chainrelate} and the fact that the underlying complex of $(M\bos N, \bys)$ is $\lk F_{\pi}(M)\otimes_{\lk I_V}\lk G_{p}(N)$  , we have:
$$
((M \bos N)\otimes_{\w {\vr}} F_{V''},\bys \otimes \I_{F_{V''}}) \simeq (\lk F_{\w {\vr}} (\lk F_{\pi}(M))\bo \lk G_{\w {\vr}}(\lk G_p(N), \by) 
$$
Additionally, by the proof of lemma \ref {lemma:Ainject}, we have: $\lk F_{\w {\vr}} (\lk F_{\pi}(M)) \simeq_A \lk F_{\w {\vr}\circ \pi}(M)=\lk F_{\pi''\circ \vr}(M) \simeq_A \lk F_{\pi''}(\lk F_{\vr}(M))$. Note that the second identity comes from the relation ${\w {\vr}}\circ \pi={\pi''\circ \vr}$. Similarly, $\lk G_{\w {\vr}}(\lk G_p(N))$ is $D$-homotopy equivalent to $\lk G_{\w {\vr}\circ p}(N)$. Therefore, following \cite {LOT}, $((M \bos N)\otimes_{\w {\vr}} F_{V''},\bys \otimes \I_{F_{V''}})$ is chain homotopic to $(\lk F_{\pi''}(\lk F_{\vr}(M)) \bo \lk G_{\w {\vr}\circ p}(N),\by) $. Likewise,
$$
((M' \bos N)\otimes_{\w {\vr}'} F_{V''},\bys \otimes \I_{F_{V''}}) \simeq (\lk F_{\pi''}(\lk F_{\vr'}(M')) \bo \lk G_{\w {\vr}'\circ p'}(N),\by)
$$
Since $(M,\{m_i\})$ is stable $A$-homotopy equivalent to $(M',\{m'_i\})$ via the triple $(\vr,\vr',W'')$, we have $\lk F_{\vr}(M) \simeq_A \lk F_{\vr'}(M')$ in the category $A_{W'',\lk I}$. By proposition \ref {prop:functor1} , $\lk F_{\pi''}(\lk F_{\vr}(M)) \simeq_A  \lk F_{\pi''}(\lk F_{\vr'}(M'))$. Furthermore, since $\w {\vr}\circ p=\w {\vr}'\circ p'$, we have $\lk G_{\w {\vr}\circ p}(N) \equiv \lk G_{\w {\vr}'\circ p'}(N)$. Therefore, $ (\lk F_{\pi''}(\lk F_{\vr}(M)) \bo \lk G_{\w {\vr}\circ p}(N),\by) \simeq (\lk F_{\pi''}(\lk F_{\vr'}(M')) \bo \lk G_{\w {\vr}'\circ p'}(N),\by)$ and thus, $((M \bos N)\otimes_{\w {\vr}} F_{V''},\bys \otimes \I_{F_{V''}}) \simeq ((M' \bos N)\otimes_{\w {\vr}'} F_{V''},\bys \otimes \I_{F_{V''}})$. As a result, $(M\bos N, \bys)$ is stably chain homotopic to $(M' \bos N,\bys)$ via the triple $(\w {\vr},\w {\vr}', V'')$. $\Diamond$          
\ \\
\ \\
The same method can be applied to prove the following theorem:
\begin {prop}
Let $W$, $W'$ and $\ww$ be vector spaces over $\Z_2$. Let $(M,\{m_i\})$ be an object of $A_{\ww,\lk I}$. Let $(N,\delta)$ and $(N',\delta')$ be objects of $D_{W,\lk I}$ and $D_{W',\lk I}$ respectively such that $(N,\delta)$ is stably homotopy equivalent to $(N',\delta')$ via the triple $(\vr,\vr',W'')$. Then $(M\bos N,\bys)$ is stably chain homotopic to $(M\bos N',\bys)$.  
\end {prop}
      
\section{Invariance under Reidemeister moves} \label {sect:IVR}
\noindent In this section, we will prove that the homotopy type of $(\rightcomplex {T}, \rightTde T)$ is invariant under the Reidemeister moves. Let $\righty {T}$ be a diagram of $\righty {\lk T}$ before the Reidemeister moves and we let $\righty {T'}$ be the diagrams of $\righty {\lk T}$ after the Reidemeister moves. First of all, we will use the ``weight moves" trick, described in section \ref {sect:movingweight}, to move the weights off the free circles in the figures \ref{fig:RI}, \ref {fig:RIIShift},  \ref {fig:RIII_LShift} and \ref {fig:RIII_RShift}, then using the cancellation method to get a deformation retraction of $(\rightcomplex {T},\rightTde {T})$. As the next step, we will implement the stable homotopy of type D structures in section \ref {sect:stable} to show that $(\rightcomplex {T},\rightTde {T})$  is (stably) homotopy equivalent to $(\rightcomplex {T'},\rightTde {T'})$. 
\ \\
\ \\
\noindent {\bf {Attention}}: In the proof of invariance under Reidemeister moves, we sometimes use the index $T$ in $\mathcal{B}\Gamma_{T,n}$ to emphasize the dependence of the ground field $\F_{\righty {T}}$ on the diagram $\righty {T}$. \\
\diagram[0.7]{RI}{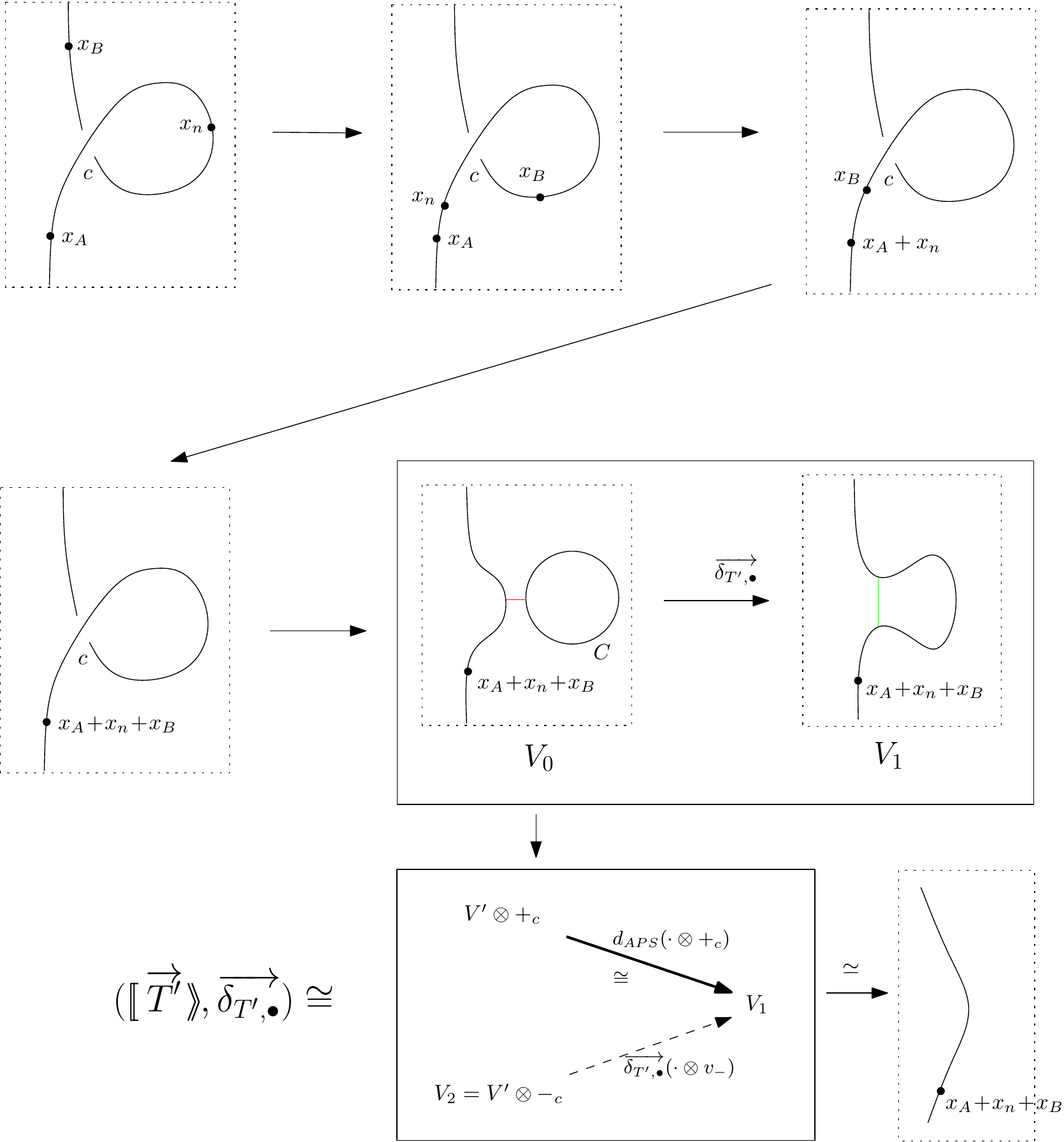}{In the top row, we use the weight shift isomorphisms to move all the weights to the bottom of the diagram. Surgering the crossing $c$ in both ways gives a finer view into the complex. Regardless of whether the local arc is on a free circle or a cleaved circle, the recorded algebra element of thickened arrow is always an invertible element of $\mathcal{B}\Gamma_{n}$ (however, the recorded algebra element of the dashed arrow depends upon the type of the local arc). When the complex is reduced along the thickened arrow, we obtain the complex for the diagram before the Reidemeister I move with the weight $x_{A} + x_{n} + x_{B}$ on the local arc.}

\subsection{Invariance under the first Reidemeister move} Figure \ref{fig:RI} shows the complex for a diagram prior to and after an Reidemeister I move applied to a right-handed crossing.  
\ \\
As usual, we can decompose $\rightcomplex {T'}=V_0\oplus V_1$ corresponding to states $(r,s)$ where $r(c)=0$ or $r(c)=1$. Furthermore, since each state generating $V_0$ always has a free circle $C$ as in the local diagram, we can continue decomposing $V_0=(V'\otimes +_c) \oplus (V'\otimes -_c)$ where $\pm_c$ are the decorations on $C$.
\ \\
 Let $(r,s)$ be a state generating $V' \otimes +_c$. As we can see, the only state $(r',s')$, which generates $V_1$ and is in the image of $\rightTde {T'}(r,s)$, is the one obtained by applying $d_{APS}$. In this case, the coefficient of $(r',s')$ in $\mathcal{B}\Gamma_{n}$ is $I_{\partial {(r,s)}}=I_{\partial {(r',s')}}$. Using this fact, we can use the cancellation lemma \ref {lemma:cancel} to cancel out $V'\otimes +_c$ and $V_1$. What we have left is $V_2=V' \otimes -_c$ with the new type $D$ structure $\overrightarrow{\delta_1}: V_2 \rightarrow\mathcal{B}\Gamma_{n}\otimes V_2$ where $\delta_1$ is the sum of $\overrightarrow{\delta'}=\rightTde{T}|_{V_2}$ and the perturbation term. 
\ \\
Taking a deeper look into the cancellation process, we see that the perturbation term arises from the following diagram: $\xi \rightarrow \xi_1 \rightarrow \xi_2 \rightarrow \xi_3$ where $\xi \in V_2$, $\xi_1$ is a generator of $V_1$ which is in the image of $\rightTde{T'}(\xi)$, $\xi_2 \in V' \otimes +_c$ which has the image $\xi_1$ under $d_{APS}$ as above (the coefficient in this case is an idempotent) and $\xi_3$ generates $\rightcomplex {T'}$ which belongs to $\rightTde {T'}(\xi_2)$. We note that $\xi_3$ has to be a generator of $(V' \otimes +_c) \oplus V_1$. Indeed, the only possibility for $\xi_3 \in V'\otimes-_c=V_2$ is that $\xi_3$ is obtained from $\xi_2$ by applying the vertical map $\rightvde {V}$ to change the decoration on $C$. But since we already moved the weight out of $C$, $\righty {w_C}=0$. Therefore, $\xi_3 \in (V' \otimes +_c) \oplus V_1$ and as a result, $\xi_3$ is canceled out by another term in either $V' \otimes +_c$ or $V_1$. Consequently, the perturbation term will be $0$ and $(V_2,\overrightarrow{\delta'}) \simeq (\rightcomplex {T'}, \rightTde {T'})$ .  
\ \\
\ \\
Although there is one-to-one corresponding between the generators of $\rightcomplex {T}$ and of $V_2$, we are still working over different fields $\field {\righty {T}}$ and $\field {\righty {T'}}$. Additionally, the local arc in $\righty {T}$ is labeled by $y_i$ while the one of $V_2$ is labeled by $x_A+x_B+x_n$. To relate them, let $\widetilde \vr:\mathcal{B}\Gamma_{T,n} \rightarrow \mathcal{B}\Gamma_{T',n}$ be the map induced by the inclusion $\vr:\field {\righty {T}} \rightarrow \field {\righty {T'}}$ defined by: $y_i \rightarrow x_A+x_B+x_n$ and $y_j \rightarrow x_j$ for $i \not= j$. Then, we define:
$$
\rightTde {T_{\vr}}: \rightcomplex {T}\otimes_{\vr} \field {\righty {T'}} \longrightarrow \mathcal{B}\Gamma_{T',n} \otimes_{\lk I_{T',n}} \big[ \rightcomplex {T}\otimes_{\vr} \field {\righty {T'}} \big]
$$
by specifying $\pairing {\rightTde {T_{\vr}}\big( (r,s) \otimes1 \big)}{(r',s')\otimes 1}=\widetilde \vr \pairing {\rightTde {T}(r,s) }{(r',s')}$. Following the proposition \ref {prop:21}, $\rightTde {T_{\vr}}$ is a type $D$ structure on $\rightcomplex {T}\otimes_{\vr} \field {\righty {T'}}$ over $(\mathcal{B}\Gamma_{T',n},\lk I_{T',n})$. 
 \ \\
If we further identify the generators $(r,s)$ of $V_2$ with $(r,s)\otimes 1$ of $\rightcomplex {T}\otimes_{\vr} \field {\righty {T'}}$, we have $\pairing {\overrightarrow {\delta'}(r,s)}{(r',s')}=\pairing {\rightTde {T_{\vr}}\big( (r,s) \otimes1 \big)}{(r',s')\otimes 1}$. Therefore, 
$$
( \rightcomplex {T}\otimes_{\vr} \field {\righty {T'}},\rightTde {T_{\vr}}) \cong (V_2,\overrightarrow {\delta'})\cong (\rightcomplex {T'}, \rightTde {T'}) 
$$ 
over $(\mathcal{B}\Gamma_{T',n},\lk I_{T',n})$. Consequently, $(\rightcomplex {T}, \rightTde {T})$ is (stably) homotopy equivalent to $(\rightcomplex {T'}, \rightTde {T'})$.  $\Diamond$
\diagram[0.75]{RIIweights}{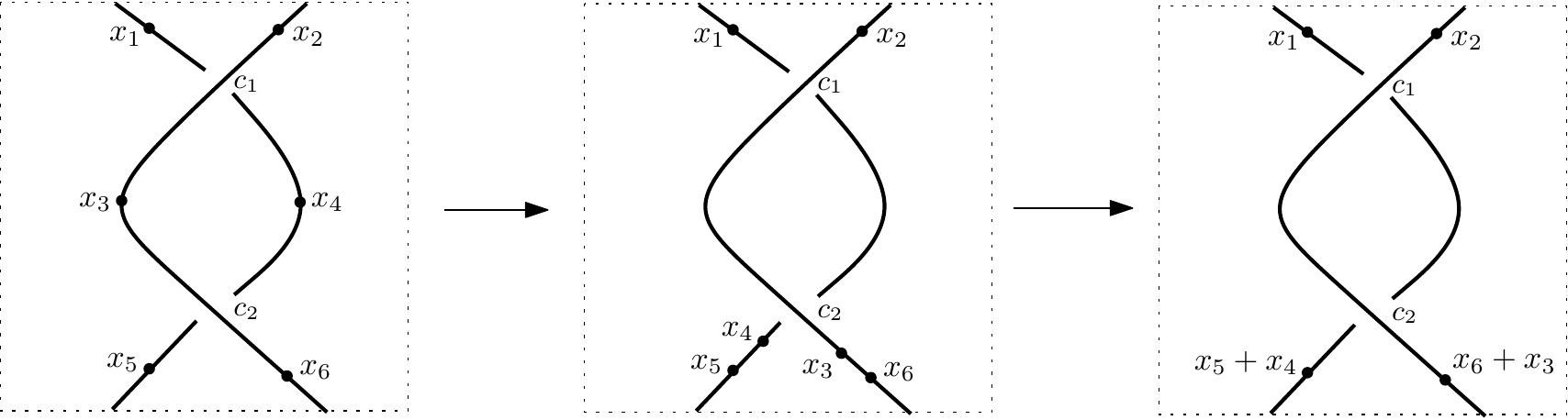}{}
\diagram[0.75]{RIIShift}{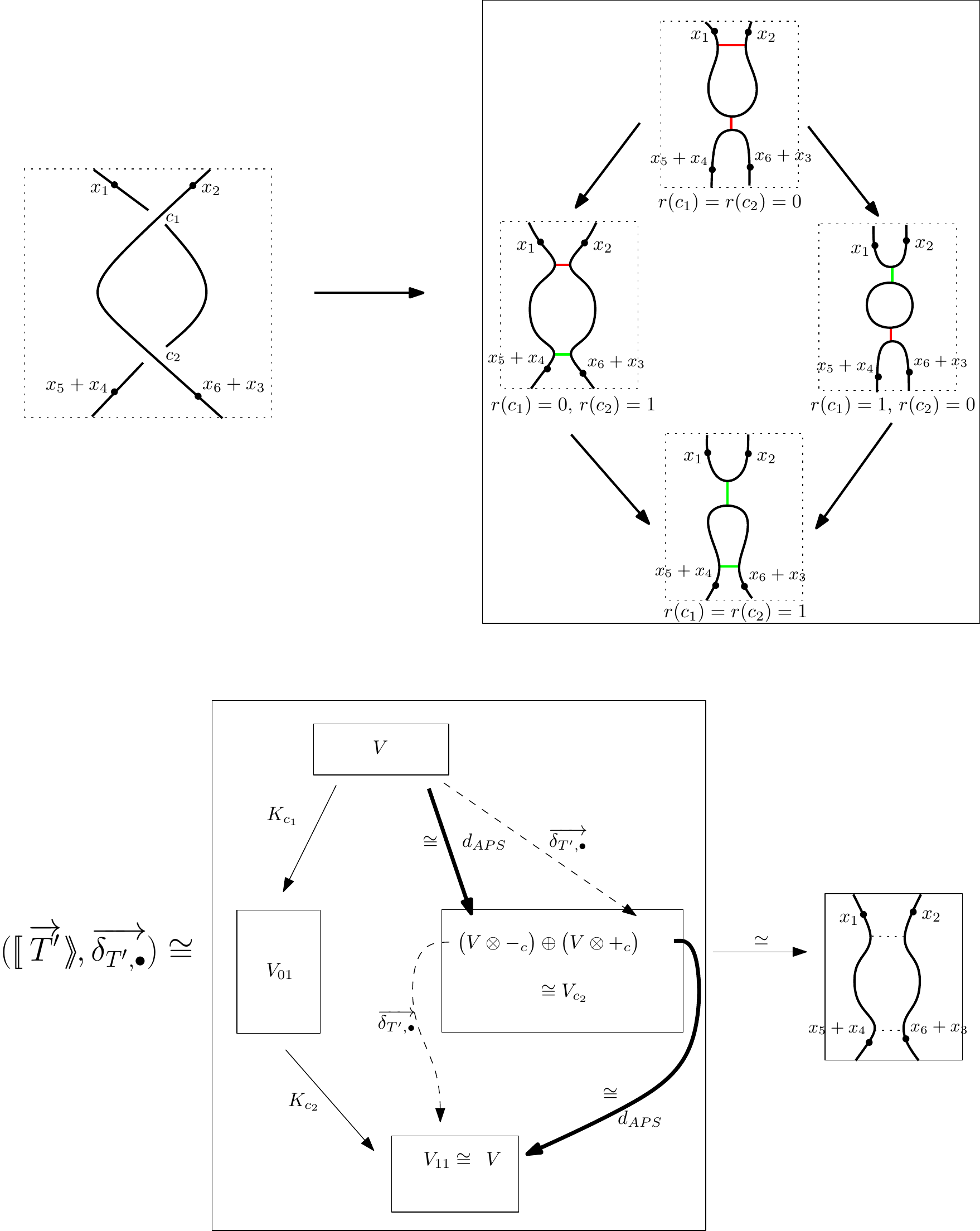}{This figure illustrates the proof of invariance under the second Reidemeister move. As we can see, regardless of whether the local arcs lie on free circles or cleaved circles, the recorded algebra elements of the thicker arrows are always idempotents. Additionally, if we cancel the bottom thicker arrow first, and then the top thicker one, we introduce no new perturbation terms since the weight on $C$ is $0$ . These cancellations produce the deformation retraction of the type $D$ structure.}  
\noindent $\bf Note.$ The invariance under a left-handed Reidemeister I move can be obtained by combining the invariance of the type $D$ structure under the right-handed Reidemeister I move and under the Reidemeister II move. \\  
\subsection{Invariance under the second Reidemeister move} Similar to the proof of the Reidemeister I invariance, we first shift all of the middle weights to the bottom in figure \ref {fig:RIIweights}. Let $c_1$, $c_2$ be two crossings in the local diagrams. For $i=0,1$, $j=0,1$, we let $V_{ij}$ be the set of states where $c_0$ is resolved by $i$ and $c_1$ is resolved by $j$. We can decompose $\rightcomplex {T'}=V_{00} \oplus V_{01}\oplus V_{10} \oplus V_{11}$. We further decompose $V_{10}=(V\otimes +_c) \oplus (V\otimes -_c)$ where $\pm_c$ are the decorations on the free circle $C$. 
Basically, the same sort of argument as in the proof of the Reidemeister I invariance can be used to explain why we can cancel out $V \otimes +_c$ and $V_{11}$ without creating any new perturbation. The reason again comes from the fact that $\righty {w_C}=0$, and thus, we do not have any maps from $V \otimes +_c$ to $V \otimes -_c$. This is illustrated in Figure \ref {fig:RIIShift}. 
\ \\

\noindent The next step is to cancel out $V_{00}$ and $V \otimes -_c$ by applying an isomorphism $I\otimes d_{APS}: V_{00} \rightarrow V\otimes -_c$. Again, no perturbation term appears since there is no map from $V_{01}$ to $V_{00}\oplus (V \otimes -_c)$. At the end, we will be left solely with $V_{01}$ and a new type $D$ structure $\delta': V_{01} \rightarrow \mathcal{B}\Gamma_{T',n} \otimes V_{01}$ where $\delta'=\rightTde {T'}|_{V_{01}}$. Note that the generators of $V_{01}$ corresponds $1$-$1$ with the generators of $\rightcomplex {T}$ but they are two possibly distinct type $D$ structures over the different ground fields. Let the weights on two local arcs of $\righty {T}$ be $y_L$ and $y_R$. To relate $\mathcal{B}\Gamma_{T,n}$ and $\mathcal{B}\Gamma_{T',n}$, we construct the $\widetilde {\vr}:\mathcal{B}\Gamma_{T,n} \longrightarrow \mathcal{B}\Gamma_{T',n}$ which is induced by the inclusion $\vr:\field {T} \longrightarrow \field {T'}$, defined by: $y_L \rightarrow x_1+x_4+x_5$, $y_R \rightarrow x_2+x_3+x_6$ and $y_i \rightarrow x_i$ for $i \not= L$ or $R$. The rest can be done exactly as in the proof of the invariance under the Reidemeister I. $\Diamond$ \\

\diagram[0.5]{RIII_LShift}{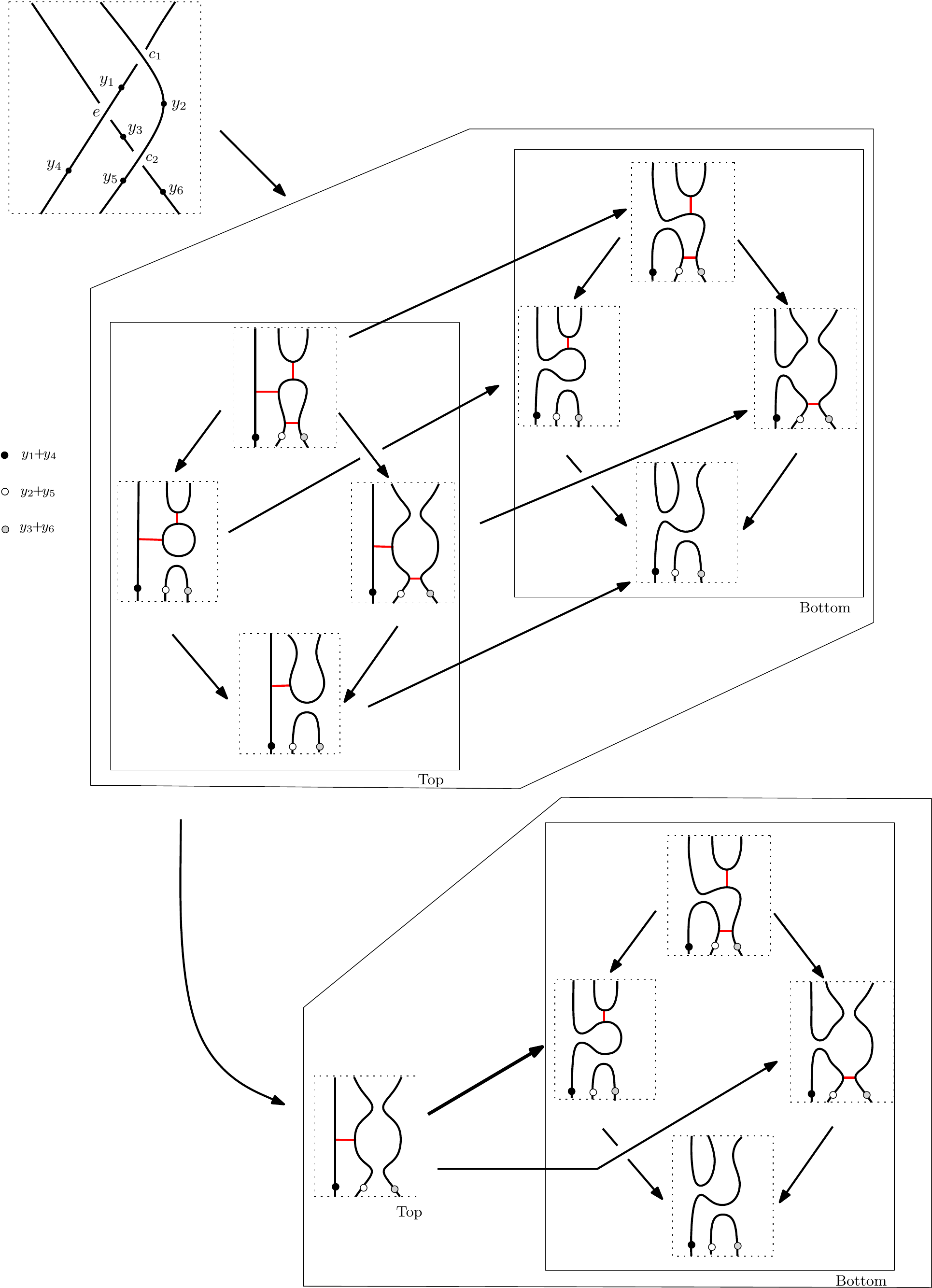}{The local picture for a diagram before the Reidemeister III move. We decompose the module along the
eight possible ways of resolving the local crossings. The four resolutions with the crossing $c$ resolved by a $0$-resolution
replicate the diagrams in the proof of Reidemeister II invariance. Using the cancellation process in the top diagram gives the
diagram at the bottom. A new perturbation map may occur from the thicker arrow in the bottom figure; however, under the identification of the generators of the lower diagrams in \ref {fig:RIII_LShift} and \ref {fig:RIII_RShift} , it will be the same as the map of the lower diagram of \ref {fig:RIII_RShift} which is obtained by surgering a bridge at the crossing $d$ .}

\diagram[0.5]{RIII_RShift}{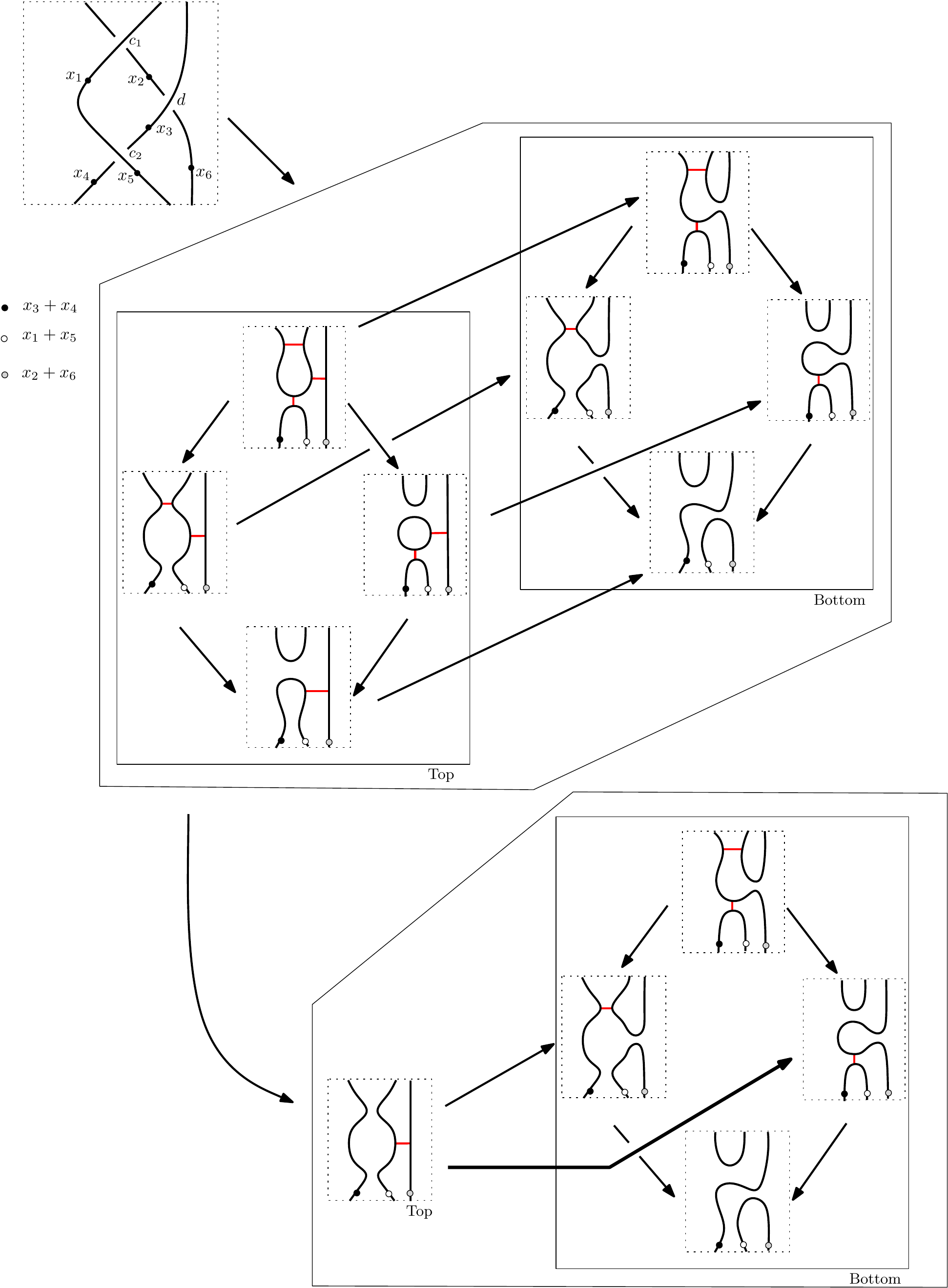}{The local picture for a diagram after the Reidemeister III move. Once again there is a new perturbation map, shown by the thicker arrow in the bottom figure.}
\subsection{Invariance under the third Reidemeister move} Since the proof of the invariance under the third Reidemeister move is similar to the proof of the invariance under the first and the second Reidemeister moves, we only mention the strategy of the proof and the figures to illustrate it. However, we will construct a homomorphism to relate $\mathcal{B}\Gamma_{T,n}$ and $\mathcal{B}\Gamma_{T',n}$ since it is somewhat different from the proof of the invariance under the Reidemeister I and II moves. The strategy is to follow these steps: 1) shift the weights to the bottom as in figure \ref {fig:RIII_LShift}, \ref {fig:RIII_RShift}. 2) use the cancellation method exactly as in the Reidemeister II move for the top faces of the top diagrams in figure \ref {fig:RIII_LShift}, \ref {fig:RIII_RShift}, we will get the lower diagrams with the two new perturbation maps. The important point is that if we canonically identify the generators of the lower diagrams in figures \ref {fig:RIII_LShift}, \ref {fig:RIII_RShift}, the maps from tops to the bottoms of the lower diagrams will be the same  (see \cite{LR1} for more detail).
\ \\
\ \\
Assuming the weights for the arcs not shown in the figures are $z_7,....,z_t$. The type $D$ structures of the first and second  lower diagrams will be denoted by $(R_1,\delta_1)$ and $(R_2,\delta_2)$. Let $(S,\delta_S)$ be a type $D$ structure with the identical maps as in the lower figure but with the black circle representing $u_1$, the white circle representing $u_2$ and the gray circle representing $u_3$. We note that $(S,\delta_S)$ is a type $D$ structure over $\mathcal{B}\Gamma_{S,n}=\mathcal{B}\Gamma_n \otimes_{\Z_2} \field {S}$ where $\field {S}$ is the field of fractions of $\Z_2 [u_1,u_2,u_3,z_7,...,z_t]$. We define the map $\vr_1: \mathcal{B}\Gamma_{S,n} \rightarrow \mathcal{B}\Gamma_{T,n}$ which is induced by the inclusion $\widetilde {\vr_1}: \field {S} \rightarrow \field {T}$: $u_1 \rightarrow y_1+y_4$, $u_2 \rightarrow y_2+y_5$ and $u_3 \rightarrow y_3+y_6$. As in the proof of the Reidemeister I invariance, this map will give us a way to relate $(S,\delta_S)$ and $(R_1,\delta_1)$ and as a result, $(S,\delta_S)$ is (stably) homotopy equivalent to $(R_1,\delta_1)$. 
\ \\
\ \\
Similarly, the map $\widetilde {\vr_2}: \mathcal{B}\Gamma_{S,n} \rightarrow \mathcal{B}\Gamma_{T',n}$, which is induced by the inclusion $\vr_2: \field {S} \rightarrow \field {T'}$: $u_1 \rightarrow x_3+y_4$, $u_2 \rightarrow y_1+y_5$ and $u_3 \rightarrow y_2+y_6$, will result in $(S,\delta_S)$ being (stably) homotopy equivalent to $(R_2,\delta_2)$. Then by Proposition \ref {prop:Dequivalent} , $(R_1,\delta_1)$ is (stably) homotopy equivalent to $(R_2,\delta_2)$. As a result, $(\rightcomplex {T},\rightTde{T})$ is (stably) homotopy equivalent to $(\rightcomplex {T'},\rightTde{T'})$. $\Diamond$

\section{A type $A$ structure in totally twisted Khovanov homology} \label {sect:Adef}
\noindent In \cite{LR2}, L. Roberts describes a type $A$ structure on the underlying module $\leftcomplex {T}$ over $ \mathcal{B}\Gamma_{n}$. It is characterized by two maps: 1) the differential $m_1$ which increases the bigrading by $(1,0)$ and thus, is an order $1$ map relative to the $\ze-$grading  and 2) the action $m_2$ which preserves bigrading and thus, $\ze-$grading preserving. For the twisting bordered homology, we will describe another type $A$ structure on $\leftcomplex {T}$ over $\mathcal{B}\Gamma_{n}$ by the following maps:\\
$$
\m{1}: \leftcomplex {T} \rightarrow \leftcomplex {T}[-1]
$$
$$
\m{2}: \acomplex {T}{T}
$$
\ \\
Let $\xi=(r,s)$ be a generator of $\leftcomplex {T}$ and $e$ be a generator of $\mathcal{B}\Gamma_{n}$. 
\ \\
\ \\
For $\m{1}$: we denote $\m{1}(\xi):=d_{APS}(\xi)+\partial_{\lk V}(\xi)$. Since the vertical map $\partial_{\lk V}$ decreases the bigrading by (0,2), $\m{1}$ is $\ze-$grading preserving into $\leftcomplex {T}[-1]$.  
 \ \\
\ \\
For $\m{2}$: we first define the action of a generator of $\mathcal{B}\Gamma_{n}$ on $\leftcomplex {T}$ 
\ \\
\ \\
$\bullet$ If $e\not= \lefty {e_C}$, then $\m{2}(\xi\otimes e):=m_2(\xi \otimes e)$ where $e$ is either an idempotent, or a bridge element, or a right decoration element and $m_2$ is defined in \cite {LR2}.
\ \\
\ \\
$\bullet$ If $e=\lefty {e_C}$, then $\m{2}(\xi\otimes \lefty {e_C}):=m_2(\xi \otimes \lefty {e_C})+\lefty {w_C} m_2(\xi \otimes \righty {e_C})=m_2(\xi \otimes \lefty {e_C})+\lefty {w_C}(r,s_C)$. The fact that $m_2$ preserves $\ze-$grading immediately result in the fact that $\m{2}$ also preserves $\ze-$grading.
\ \\
\ \\
To define the action of a general element of $\mathcal{B}\Gamma_{n}$ on $\leftcomplex {T}$, for $p_1, p_2 \in\mathcal{B}\Gamma_{n}$, we impose the relation:
\begin {equation} \label{eqn:m2}
\m{2}(\xi \otimes p_1p_2)=\m{2}(\m{2}(\xi \otimes p_1)\otimes p_2)
\end {equation}
For this definition to be well defined, we need to prove the following proposition:
\begin{prop} If two products of the generators $p_1$ and $p_2$ define equal elements in $\mathcal{B}\Gamma_{n}$ then $\m{2}(\xi\otimes p_1)=\m{2}(\xi\otimes p_2)$.
 \end{prop}
\noindent $\bf Proof.$ It suffices to prove $\m{2}(\xi \otimes p)=0$ if $p$ is a relation defining $\mathcal{B}\Gamma_{n}$. First of all, we recall the two following facts:\\
\begin {enumerate}
\item $m_2(\xi \otimes p)=0$ as in \cite{LR2}

\item $\m{2} (\xi \otimes e)=m_2(\xi \otimes e)$ for every generator of $\mathcal{B}\Gamma_{n}$ unless $e=\lefty {e_C}$ 
\end {enumerate}
\ \\
Combining these two facts and using the relation (\ref {eqn:m2}), we have if $p$ does not involve $\lefty {e_C}$, then $\m{2}(\xi \otimes p)= m_2(\xi \otimes p)=0$. If $p$ involves $\lefty {e_C}$, let $p_1$ be an element of $ \mathcal{B}\Gamma_{n}$ obtained from $p$ by substituting $\righty {e_C}$ for each term $\lefty {e_C}$ appeared in $p$. In this situation, we have the following cases:
\ \\
\ \\
$\bf {Case}$ $\bf {I}$. If $p_1$ is a relation defining $\mathcal{B}\Gamma_{n}$, then we have two possibilities:
\ \\
$\bullet$ If $p=\lefty {e_C} \lefty {e_D} +\lefty {e_D}\lefty {e_C}$: Since $\m{2} (\xi \otimes \lefty {e_C})=m_2(\xi \otimes \lefty {e_C})+\lefty {w_C} m_2(\xi\otimes \righty {e_C})$ and $\m{2} (\xi \otimes \lefty {e_D})=m_2(\xi \otimes \lefty {e_D})+\lefty {w_D}m_2(\xi\otimes \righty {e_D})$, we have:
\ \\
$\m{2} (\xi \otimes p )$
\ \\
\ \\
\indent $=\m{2} \big(\xi \otimes(\lefty {e_C}\lefty {e_D} +\lefty {e_D}\lefty {e_C})\big)$
\ \\
\ \\
\indent $= \m{2}(\m2(\xi \otimes \lefty {e_C}) \otimes \lefty {e_D}) +\m{2}(\m2(\xi \otimes \lefty {e_D}) \otimes \lefty {e_C})$
\ \\
\ \\
\indent $=\m{2}\big(\big[ m_2(\xi \otimes \lefty {e_C})+\lefty {w_C} m_2(\xi\otimes \righty {e_C}) \big]\otimes \lefty {e_D}\big)+\m{2}\big(\big[ m_2(\xi \otimes \lefty {e_D})+\lefty {w_D} m_2(\xi\otimes \righty {e_D}) \big]\otimes \lefty {e_C}\big)$ 
\ \\
\ \\
\indent $=\big[m_2(\xi \otimes \lefty {e_C}\lefty{e_D})+m_2(\xi \otimes \lefty {e_D}\lefty{e_C})\big] + \lefty {w_D}\big[m_2(\xi \otimes \lefty {e_C}\righty{e_D})+m_2(\xi \otimes \righty {e_D}\lefty{e_C})\big]+$
\ \\
\ \\
\indent $\hspace {0.12 in}$$\lefty {w_C}\big[m_2(\xi \otimes \lefty {e_D}\righty{e_C})+m_2(\xi \otimes \righty {e_C}\lefty{e_D})\big]+ \lefty {w_C}\lefty {w_D}\big[ m_2(\xi \otimes \righty {e_C}\righty{e_D})+m_2(\xi \otimes \righty {e_D}\righty{e_C})\big]$ 
\ \\
\ \\
\indent $=0$ since $m_2(\xi \otimes q)=0$ for every $q$ which is a relation defining $\mathcal{B}\Gamma_{n}$.
\ \\
\ \\
$\bullet$ If $\lefty {e_C}$ is the only left decoration element involved in the relation $p$. We have: $\m{2}(\xi\otimes p)=m_2(\xi \otimes p)+\lefty {w_C}m_2(\xi \otimes p_1)=0$. 
\ \\
\ \\
$\bf {Case}$ $\bf {II}$. If $p_1$ is not a relation defining $\mathcal{B}\Gamma_{n}$, we see that the only relations which involve left decoration elements and satisfies that $p_1$ is not a relation defining $\mathcal{B}\Gamma_{n}$, come from either merging two $+$ cleaved circles or dividing a $+$ cleaved circle by surgering along a bridge $\gamma$. Since the proof of the case of division is similar, we just present the proof of the case of merging. In this case, $p=\lefty {e_{C_1}}m_{\ga_1}+\lefty {e_{C_2}}m_{\ga_2}+m_\ga\lefty {e_C}$. We need to prove:
\begin {equation} \label {eqn:checkm2}
\m{2}(\xi \otimes\lefty {e_{C_1}}m_{\ga_1} )+\m{2}(\xi \otimes\lefty {e_{C_2}} m_{\ga_2} )+\m{2}(\xi \otimes m_\ga  \lefty {e_C} )=0
\end {equation}   
Rewriting the left side:
\ \\
$\Leftrightarrow \m{2}(\m{2}(\xi \otimes\lefty {e_{C_1}}) \otimes m_{\ga_1} )+\m{2}(\m{2}(\xi \otimes\lefty {e_{C_2}}) \otimes m_{\ga_2} )+\m{2}(\m{2}(\xi \otimes m_\ga) \otimes \lefty {e_C} )=0$
\ \\
\ \\
$\Leftrightarrow \big[m_2(m_2(\xi \otimes\lefty {e_{C_1}}) \otimes m_{\ga_1} )+m_2(m_2(\xi \otimes\lefty {e_{C_2}})\otimes m_{\ga_2} )+m_{2}(m_{2}(\xi \otimes m_\ga)\otimes \lefty {e_C} )\big]+$
\ \\
\ \\
$\big[  m_2(\lefty {w_{C_1}}m_2(\xi \otimes\righty {e_{C_1}}) \otimes m_{\ga_1} )+m_2(\lefty {w_{C_2}}m_2(\xi \otimes\righty {e_{C_2}})\otimes m_{\ga_2} )+\lefty {w_C}m_{2}(m_{2}(\xi \otimes m_\ga)\otimes \righty {e_C} )  \big]=0$
\ \\
\ \\
The sum of the first bracket will be $0$ since $m_2(\xi\otimes p)=0$ and thus it suffices to prove:
$  m_2(\lefty {w_{C_1}}m_2(\xi \otimes\righty {e_{C_1}}) \otimes m_{\ga_1} )+m_2(\lefty {w_{C_2}}m_2(\xi \otimes\righty {e_{C_2}})\otimes m_{\ga_2} )+\lefty {w_C}m_{2}(m_{2}(\xi \otimes m_\ga)\otimes \righty {e_C} )  =0$ 
\ \\
Since $\lefty {w_C}=\lefty {w_{C_1}}+\lefty {w_{C_2}}$, we can rewrite the left side:
\ \\
$\Leftrightarrow \lefty {w_{C_1}}\big[  m_2(\xi \otimes  \righty {e_{C_1}} m_{\ga_1} )+m_2(\xi \otimes m_\ga \righty {e_{C}} )\big]+\righty {e_{C_2}}\big[  m_2(\xi \otimes  \righty {e_{C_2}} m_{\ga_2}) +m_2(\xi \otimes m_\ga \righty {e_{C}} )\big]=0$
\ \\
The sums of the first and second brackets equal $0$ since  $\righty {e_{C_1}} m_{\ga_1} + m_\ga \righty {e_{C}}$ and $\righty {e_{C_2}} m_{\ga_2} +m_\ga \righty {e_{C}}$ are the relations defining $\mathcal{B}\Gamma_{n}$. Identity \ref {eqn:checkm2} is therefore true, so, $\m{2}$ is well-defined. $\Diamond$
\section{The properties of the type A structure} \label {sect:Aprp}
\noindent In this section, we will prove that $(\leftcomplex {T},\m{1},\m{2})$ is an $A_\infty$ module over the differential graded algebra $\mathcal{B}\Gamma_{n}$ with $\m {n}=0$ for $n \geq 3$.\\
\begin {prop} Let $\xi=(r,s)$ be a generator of $\leftcomplex {T}$ and $p_1,p_2 \in \mathcal{B}\Gamma_{n}$. The maps $\m{1}$ and $\m {2}$ satisfy the following relations:
\ \\
\begin{enumerate}
\item $\m{1}(\m {1}(\xi))=0$
\item $\m{2}(\m{1} (\xi) \otimes p_1)+\m{2}(\xi \otimes d_{\Gamma_{n}}(p_1))+ \m{1} (\m {2}(\xi \otimes p_1))=0 $
\item $\m{2}(\xi \otimes p_1p_2)=\m{2}(\m{2}(\xi \otimes p_1)\otimes p_2)$
\end {enumerate}
\end {prop}
\noindent $\bf Proof.$ The first identity comes from the fact that $\m{1}$ is the differential of the complex $\leftcomplex {T}$. The third identity comes from the fact that the construction of $\m{2}$ is well defined, as shown in section \ref {sect:Adef}. Therefore, we only need to verify the second identity. By bootstrapping the relation of longer words described in proposition 21 of \cite {LR2}, it suffices to prove the second identity when $p_1$ is a generator of  $\mathcal{B}\Gamma_{n}$. We have the following two cases:\\
\ \\
$\bf Case$ $\bf I.$ If $p_1 \not= \lefty {e_C}$, using the fact that $d_{\Gamma_{n}}(p_1)=0$ and the relations between $\m{1}$, $\m{2}$ and $m_1$, $m_2$, we rewrite the second identity as
$$
\big [m_{2}(m_{1} (\xi) \otimes p_1)+m_{2}(\xi \otimes d_{\Gamma_{n}}(p_1))+ m_{1} (m_ {2}(\xi \otimes p_1))\big]+\big[ m_2(\partial_{\lk V}(\xi) \otimes p_1)+\partial_{\lk V}(m_2(\xi \otimes p_1)) \big]=0
$$
Due to the fact that $(\leftcomplex {T},m_1,m_2)$ is an $A_\infty$ module over the differential graded algebra $\mathcal{B}\Gamma_{n}$, the sum in the first bracket is $0$. Therefore, we only need to prove that: 
\begin{equation} \label {eqn:aprop}
m_2(\partial_{\lk V}(\xi) \otimes p_1)+\partial_{\lk V}(m_2(\xi \otimes p_1))=0
\end{equation}
There are three possibilities for $p_1$:
\ \\
\begin {enumerate}
\item If $p_1$ is an idempotent $I_{(L,\sigma)}$, then both terms will equal $\partial_{\lk V}(\xi)$ if $\partial(r,s)=I_{(L,\sigma)}$ and they are both $0$ if $\partial(r,s)\not= I_{(L,\sigma)}$. As the result, identity \ref {eqn:aprop} is true in this case.

\item If $p_1=\righty {e_C}$, then $m_2(\partial_{\lk V}(\xi) \otimes p_1)=\partial_{\lk V}(m_2(\xi \otimes p_1))=\sum\limits_{\text {D is + free circle}} \lefty {w_D}(r,s_{C,D})$.

\item  If $p_1$ is a bridge element corresponding to surgering along a bridge $\gamma$, then
$$
m_2(\partial_{\lk V}(\xi) \otimes p_1)=\partial_{\lk V}(m_2(\xi \otimes p_1))=\sum_{\gamma, \text {D is + free circle}}(r_{\gamma,D},s_{\gamma,D})
$$ 
\noindent Therefore, identity \ref {eqn:aprop} is true in this case.\\
\end {enumerate}
\ \\
$\bf Case$ $\bf II.$ If $p_1=\lefty {e_C}$, we need to prove the following:
\ \\ 
$\m{2}(\m{1} (\xi) \otimes \lefty {e_C})+\m{2}(\xi \otimes d_{\Gamma_{n}}(\lefty {e_C}))+ \m{1} (\m {2}(\xi \otimes \lefty {e_C}))=0 $
\ \\
\ \\
$\Leftrightarrow \m{2}(m_1 (\xi) \otimes \lefty {e_C})+\m{2}(\partial_{\lk V}(\xi) \otimes \lefty {e_C})+m_2(\xi \otimes d_{\Gamma_{n}}(\lefty {e_C}))+\m{1} (m_2(\xi \otimes \lefty {e_C}))+$
\indent $\lefty {w_{C}}\m{1} (m_{2}(\xi \otimes \righty {e_C}))=0$
\ \\
\ \\
$\Leftrightarrow \big [ m_{2}(m_1 (\xi) \otimes \lefty {e_C}) +  \sum \limits_{\ga} \lefty {w_{C_\ga}} m_{2}(\xi_{\ga} \otimes \righty {e_C}) \big]+\big[ m_{2}(\partial_{\lk V}(\xi) \otimes \lefty {e_C})+\lefty {w_{C}}m_{2}(\partial_{\lk V}(\xi) \otimes \righty {e_C})\big]+$
\indent $m_2(\xi \otimes d_{\Gamma_{n}}(\lefty {e_C}))+\big[ m_{1} (m_2(\xi \otimes \lefty {e_C}))+\partial_{\lk V}(m_2(\xi \otimes \lefty {e_C}))\big] +\big[\lefty {w_{C}}m_{1} (m_{2}(\xi \otimes \righty {e_C}))+$
\ \\
\ \\
\indent $\lefty {w_{C}}\partial_{\lk V}(m_{2}(\xi \otimes \righty {e_C}))\big]=0$ 
\ \\
where $\xi_{\gamma}$ is in the image of $m_1(\xi)$, $C_{\ga}$ corresponds to $C$ in $\partial (\xi_{\ga})$ and $\gamma$ is taken over the resolution bridges which contributes to $m_1$. Rewriting the left hand side:
\ \\
$\Leftrightarrow \big[ m_{2}(m_{1} (\xi) \otimes \lefty {e_C})+m_{2}(\xi \otimes d_{\Gamma_{n}}(\lefty {e_C}))+ m_{1} (m_{2}(\xi \otimes \lefty {e_C}))\big]+$
\ \\ 
\ \\
\indent $ \lefty {w_{C}} \big[ m_{2}(\partial_{\lk V}(\xi) \otimes \righty {e_C})+v(m_{2}(\xi \otimes \righty {e_C}))\big] +$
\ \\
\ \\
\indent $\big[  \sum \limits_{\ga} \lefty {w_{C_\ga}} m_{2}(\xi_{\ga} \otimes \righty {e_C}) +m_{2}(\partial_{\lk V}(\xi) \otimes \lefty {e_C})+ \partial_{\lk V}(m_2(\xi \otimes \lefty {e_C}))+ \lefty {w_{C}}m_{1} (m_{2}(\xi \otimes \righty {e_C})) \big]=0 $
\ \\
The sum inside the first bracket is $0$ since $(\leftcomplex {T},m_1,m_2)$ is an $A_\infty$ module over the differential graded algebra $\mathcal{B}\Gamma_{n}$. The sum inside the second bracket is also $0$ by identity (\ref {eqn:aprop}). Therefore, we only need to make sure the third sum is also $0$:
\begin {equation}\label {eqn:lastidentity}
\sum \limits_{\ga} \lefty {w_{C_\ga}} m_{2}(\xi_{\ga} \otimes \righty {e_C}) +m_{2}(\partial_{\lk V}(\xi) \otimes \lefty {e_C})+ \partial_{\lk V}(m_2(\xi \otimes \lefty {e_C}))+ \lefty {w_{C}}m_{1} (m_{2}(\xi \otimes \righty {e_C})) =0 
\end {equation}
Since we have the following identity:
$$
m_2(m_1(\xi)\otimes \righty {e_C})=m_2(m_1(\xi)\otimes d_{\Gamma_{n}}(\righty {e_C}))+m_1(m_2(\xi\otimes \righty {e_C}))=m_1(m_2(\xi\otimes \righty {e_C})),
$$
we can rewrite (\ref  {eqn:lastidentity}) as following:
$$
\sum \limits_{\ga} m_{2}((\lefty {w_{C_\ga}}+\lefty {w_{C}})\xi_{\ga} \otimes \righty {e_C})+m_{2}(\partial_{\lk V}(\xi) \otimes \lefty {e_C})+ \partial_{\lk V}(m_2(\xi \otimes \lefty {e_C}))=0
$$
The second term contains pairs of $(D,\ga)$ from changing the decoration on a $+$ free circle $D$ first, then resolving an active resolution bridge $\ga$, which changes the decoration on $C$ from $+$ to $-$. On the other hand, the third term will be the sum of generators $(\ga_1,D_1)$ coming from surgering an active resolution bridge $\ga_1$ to change the decoration on $C$  first, then changing a $+$ free circle $D_1$ of $\xi_{\ga_1}$ to $-$. Taking the sum of the second and the third terms, the pairs $(D,\ga)$ in the second term will be canceled out by the reverse pair $(\ga,D)$ if $(\ga,D)$ belongs to the third term and vice versa. However, there are two exceptional cases when reversing a pair of the second (third) term does not belong to the third (second) term: 1) $(D,\ga)$ where $D$ is a $+$ free circle of $\xi$ and $\ga$ has one foot on $D$ and another on $C$ or 2) $(\ga,D_\ga)$ where $\ga$ is an active resolution bridge whose feet belongs $C$ and $D_\ga$ is new $+$ free circle which is created by surgering $\ga$. The first term contains the generators whose coefficients $\lefty {w_{C_\ga}}+ \lefty {w_{C}}\ne 0$ if and only if the active bridge $\ga$ has at least one foot on $C$. By comparing the ``weight" coefficients, those generators will be canceled out by generators in the above two exceptional cases. As a result, the sum of the three terms is $0$.  $\Diamond$ 
  
\section {Invariance of the type $A$ structure under the weight moves and Reidemeister moves } \label {sect: Ainvariant}
\subsection {Invariance under the weight moves} \label {sect:amove}
\ \\
\noindent In section \ref {sect:movingweight}, we proved that the type $D$ structure in the twisted tangle homology described in section \ref {sect:3} is invariant under the weight moves by using the trick in \cite {NR} or \cite {TJ}. Similarly, in this section, we will prove that under the weight moves, the type $A$ structure in the twisted tangle homology is an invariant. Additionally, based on the construction of the type $A$ and the type $D$ structures, there exists a type $DA$ bimodule version in twisted tangle homology for an $(m,n)$-tangle (This work is joint with L. Roberts and is in the preparation to submit). Then thanks to the gluing process in \cite {LOT1} which pairs the type $A$ and type $DA$ structures, it suffices to prove the invariance for the local case (in the following figure), and the invariance for the global case follows as a consequence.\\
$$
\inlinediag{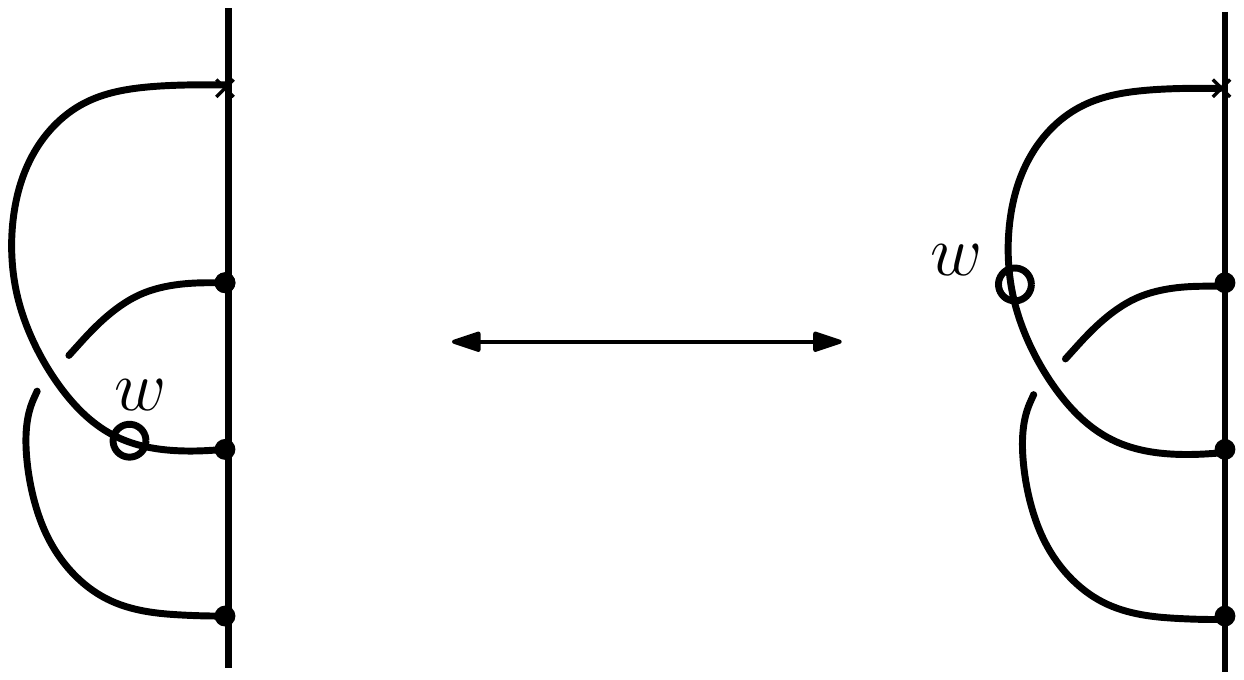} \hspace{0.75in} \inlinediag{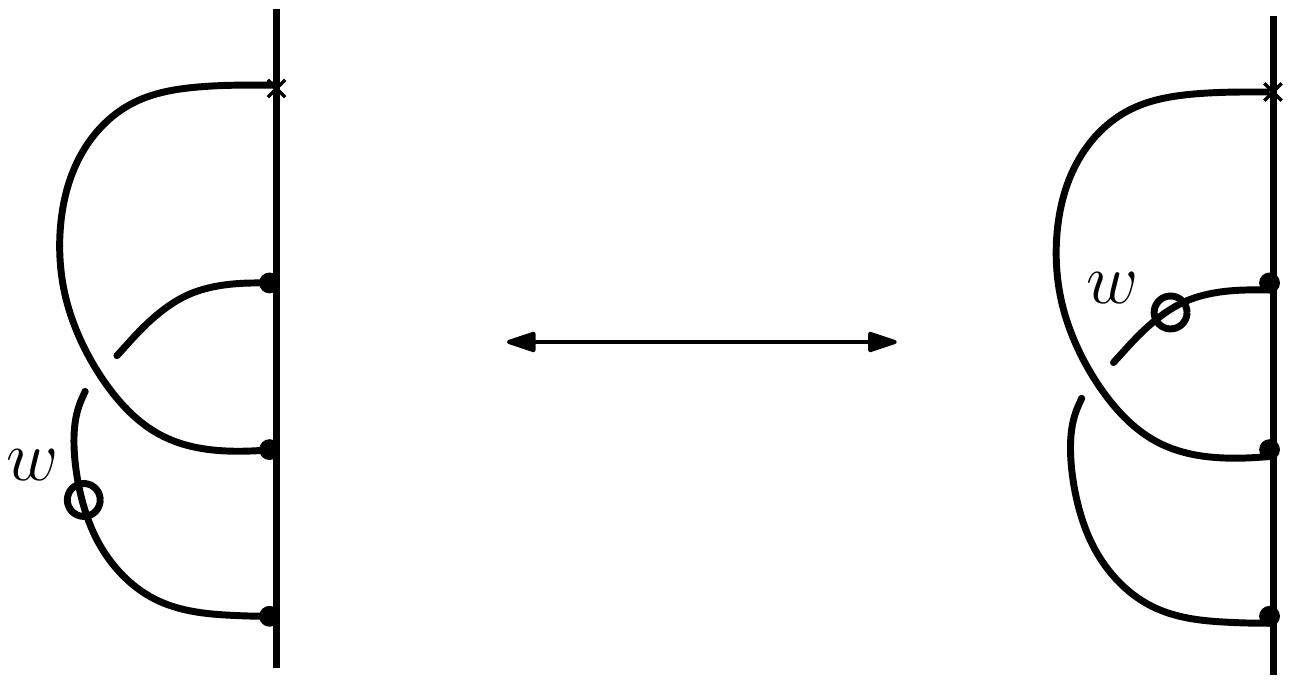} 
$$

\noindent We will only give a proof for transition on the left in the above figure, since the proof of the right case is similar. Let $\lefty {T}$ be the weighted tangle before the weight $w$ moved along the crossing $c$ and let $\lefty {T'}$ be the weighted tangle after the movement of $w$. Let $\{\m {1}, \m {2}\}$ and $\{\mm {1},\mm {2}\}$  be the maps defining the type $A$ structures for $\leftcomplex {T}$ and $\leftcomplex {T'}$ respectively. Note that in this case, both $\m {1}$ and $\mm {1}$ are $0$-maps. We will construct an $A_{\infty}$ morphism $\Psi=\{\psi_1,\psi_2\}$ from $\leftcomplex {T}$ to $\leftcomplex {T'}$ as following:\\
\ \\

\begin{center}
\begin{figure}
\includegraphics[scale=0.5]{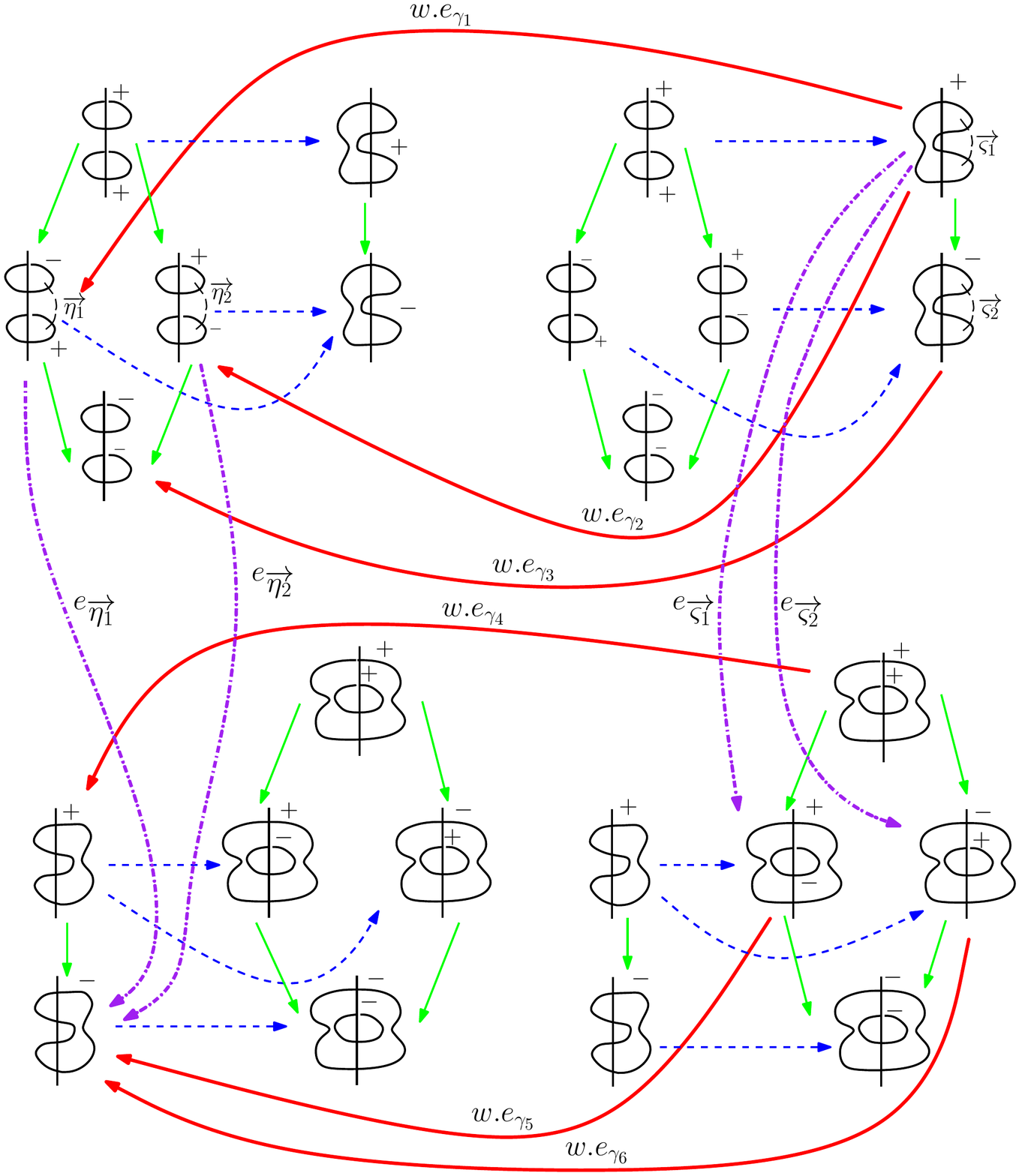}
\caption{This figure illustrates the definition of $\psi_2$. In this figure, the thick red arrows define the map $\psi_2$. Additionally, the dashed dotted purple, the blue and the green arrows stand for the action of the right, left bridge elements and the right decoration elements respectively. The symbol above each arrow specifies the element in $\mathcal{B}\Gamma_{4}$ acting on the complex.}
\label{fig:achainmove}
\end{figure}
\end{center}
$\bullet$ $\psi_1:\leftcomplex {T} \longrightarrow \leftcomplex {T'}$ is the identity map since the generators of $\leftcomplex {T}$ can be canonically identified with the generators of $\leftcomplex {T'}$.
\ \\
$\bullet$ $\psi_2: \awcomplex {T}{T'} [1]$ is defined based on its image on $\xi \otimes e$ where $\xi=(r,s)$ is a generator of $\leftcomplex {T}$ and $e$ is a generator of $\mathcal{B}\Gamma_{2}$. We define $\psi_2(\xi \otimes e)=w.(r_{\gamma},s_{\gamma})$ if $\xi$ and $e$ satisfy the following conditions:
\begin {enumerate}
\item $\partial (r,s)=s(e)$, $r(c)=1$ and $e$ is a left bridge element.
\item $r_{\gamma}$ is obtained by surgering along the inactive bridge resolution $\gamma$ at crossing $c$ and $s_\ga$ is computed by using the Khovanov sign rules (Note that: $r_{\gamma}$ is a resolution of $\lefty {T'}$ under the identification $\lefty {T}$ and $\lefty {T'}$) and  $(r_{\gamma},s_{\gamma})=t(e)$.
\end {enumerate}  
Otherwise, we define $\psi_2(\xi \otimes e)=0$. We also describe the definition of $\psi_2$ as the thick red arrows of the figure \ref {fig:achainmove}. Then we will define $\psi_2$ when $e$ is any element in $\mathcal{B}\Gamma_{2}$ by imposing the following relation for each $\xi \in \leftcomplex {T}$ and $e_1,e_2 \in \mathcal{B}\Gamma_{2}$:
\begin {equation} \label {eqn:impose}
\psi_2(\xi \otimes e_1e_2)=\mm {2} (\psi_2(\xi \otimes e_1)\otimes e_2)+ \psi_2(\m {2}(\xi \otimes e_1)\otimes e_2)
\end {equation}
For $\psi_2$ to be well-defined, we need to verify that with the relation we just imposed, $\psi_2(\xi\otimes p)=0$ for each relation $p$ defining $\mathcal{B}\Gamma_{2}$. Since for each $e$ generating $\mathcal{B}\Gamma_{2}$, $\psi_2(\xi\otimes e)=0$ unless $e$ is a left bridge element, $\psi_2(\xi \otimes p)$ is trivially $0$ if $p$ does not involve left bridge element(s). Therefore, there are two cases to verify: 1) if $p$ is the relation (\ref {rel:3}) of {\bf Disjoint support and squared bridge relations} and 2) if $p$ is the relation (\ref {rel:7}) of {\bf Relations for decoration edges}. With the aid of Figure \ref {fig:achainmove}, the proof is straightforward. We will do one example to illustrate the method: Since we have relation: $p=e_{\gamma_1}e_{\righty{\eta_1}}+e_{\gamma_2}e_{\righty {\eta_2}}=0$, which is the square relation, we need to verify $ \psi_2(\xi \otimes e_{\gamma_1}e_{\righty {\eta_1}})=\psi_2(\xi \otimes e_{\gamma_2}e_{\righty {\eta_2}} )$ where $\xi$ is the right top corner state in Figure \ref {fig:achainmove}. Indeed, since $\xi(c)=r(c)=1$, using relation (\ref {eqn:impose}) and the fact that $\m {2}(\xi \otimes e_{\ga_1})=0$ (because its resolution bridge is inactive), we have:
$$
\psi_2(\xi \otimes e_{\gamma_1}e_{\righty {\eta_1}})=\mm {2}(\psi_2(\xi \otimes e_{\gamma_1})\otimes e_{\righty {\eta_1}})
$$
Thus, in Figure \ref {fig:achainmove}, we go along the thick red arrow to get the result of $ e_{\gamma_1}$ acting on $\xi$ and follow the dashed dotted purple curve due to the action of $e_{\righty {\eta_1}}$. Similarly, $\psi_2(\xi \otimes e_{\gamma_2}e_{\righty {\eta_2}} )$ is calculated by first going along the red arrow under the action of ${e_{\gamma_2}}$ and then, following by the dashed dotted purple arrow $e_{\righty {\eta_2}}$. Since the result will be the same if we follow either of those two paths, we finish the proof of $\psi_2$ to be well-defined in this case. The proof of other cases can be handled by the same method and we leave the verification to the readers.  
\ \\
\begin {prop}
$\Psi=\{\psi_1,\psi_2\}$ is an $A_{\infty}$ morphism from $\leftcomplex {T}$ to $\leftcomplex {T'}$
\end {prop}
\noindent {\bf Proof}. We need to verify the three following conditions:\\
\begin {enumerate}
\item $\psi_1$ is the chain map from $(\leftcomplex {T},\m {1})$ to $(\leftcomplex {T'},\mm {1})$
\item $\mm {2}(\psi_1(\xi)\otimes e)+\mm {1}(\psi_2(\xi \otimes e))+\psi_1(\m {2} (\xi \otimes e))+\psi_2(\m {1}(\xi)\otimes e)+\psi_2(\xi \otimes d_{\Gamma_2}(e) )=0 $
\item $\psi_2(\xi \otimes e_1e_2)=\psi_2 (\m {2}(\xi \otimes e_1)\otimes e_2) +\mm {2}(\psi_2(\xi \otimes e_1)\otimes e_2)$
\end {enumerate}
The first condition is trivially true since both $\m {1}$ and $\mm {1}$ are zero maps. The third one comes from relation (\ref {eqn:impose}) that we impose on $\psi_2$. For the second condition to be verified, it suffices to prove that:
\begin {equation} \label {equation:10}
\psi_2(\xi \otimes d_{\Gamma_2}(e) )=\mm {2}(\psi_1(\xi)\otimes e)+\psi_1(\m {2} (\xi \otimes e))
\end {equation}
\noindent First of all, we will prove that (\ref {equation:10}) is true when $e$ is a (length 0 or 1) generator of $\mathcal{B}\Gamma_{2}$. If $e$ is not a left decoration element, the left hand side is $0$ by the definition of $d_{\Gamma_2}$. Similarly, the right hand side is also $0$ because the action of $e$ on $\leftcomplex {T'}$ does not change after we move the weight and thus, $\mm {2}(\psi_1(\xi)\otimes e)=\psi_1(\m {2} (\xi \otimes e))$. If $e$ is a left decoration element $\lefty {e_C}$ where $C$ is $+$ cleaved circle of $\partial (\xi)$, there are two possibilities:\\
\ \\
\ \\
\indent 1) $C$ is the only cleaved circle of $\partial (\xi)$ then:\\
$$
\mm {2}(\psi_1(\xi)\otimes e)+\psi_1(\m {2}(\xi \otimes e))=0
$$
because both terms in the identity are $\lefty {w_C} (r,s_C)$ (the moved weight $w$ is still on $C$ in this case)
\ \\
\ \\
\indent 2) $C$ is not the only one, then:\\
$$
\mm {2}(\psi_1(\xi)\otimes e)+\psi_1(\m {2}(\xi \otimes e))=w. (r,s_C)
$$
On the other hand, by the definition of $d_{\Gamma_2}$ and relation (\ref {eqn:impose}), $\psi_2(\xi \otimes d_{\Gamma_2}(\lefty {e_C}))$ is calculated as the sum of paths starting at $\xi$, then following either the dashed blue arrow (the action of a left bridge element) or the thick red arrow (the action of $\psi_2$), and then following either the thick red arrow or the dashed blue arrow (see figure \ref {fig:achainmove}). Note that if $C$ is the only cleaved circle of $\xi$ then there are two such paths and their sum will be canceled out. Otherwise, we have only one such path and the end point of this path is $(r,s_C)$. The reason why we have only one path in this case is that according to the definition:
$$
\psi_2(\xi \otimes d_{\Gamma_2}(\lefty {e_C}))=\psi_2(\xi \otimes \lefty {e_{\ga}}\lefty {e_{\ga^{\dagger}}} )=\psi_2 (\m {2}(\xi \otimes \lefty {e_\ga})\otimes \lefty {e_{\ga^{\dagger}}}) +\mm {2}(\psi_2(\xi \otimes \lefty {e_\ga})\otimes \lefty {e_{\ga^{\dagger}}})
$$
and depending on $\xi(c)=0$ or $\xi(c)=1$, the second term or the first term in the latest sum will disappear. Furthermore, the weight $w$ comes from the chain map $\psi_2$. Therefore, we finish the proof when $e$ is of length $0$ or $1$.
\ \\
\ \\
\noindent Next we prove that identity (\ref {equation:10}) is true when $e$ is a general element in $\mathcal{B}\Gamma_{2}$ by bootstrapping the relation of longer words. Let $e_1$ and $e_2$ be two elements of $\mathcal{B}\Gamma_{2}$. For an ease of notation, we let $m_2$, $m'_2$ and $d$ stand for $\m {2}$, $\mm {2}$ and $d_{\Gamma_2}$ respectively. Then
\ \\
\ \\
$
\psi_2(\xi \otimes d(e_1e_2) )
$
\ \\
\ \\
\indent$=\psi_2(\xi \otimes de_1e_2 )+\psi_2(\xi \otimes e_1de_2 )$
\ \\
\ \\
\indent $=\psi_2 (m_{2}(\xi \otimes de_1)\otimes e_2) +m'_{2}(\psi_2(\xi \otimes de_1)\otimes e_2)+\psi_2 (m_2(\xi \otimes e_1)\otimes de_2)+$
\ \\
\ \\
\indent $\hspace {0.12in}$ $m'_{2}(\psi_2(\xi \otimes e_1)\otimes de_2)$
\ \\
\ \\
\indent $=\big[\psi_2 (m_{2}(\xi \otimes de_1)\otimes e_2) +m'_{2}(\psi_2(\xi \otimes e_1)\otimes de_2)\big]+m'_{2}(\psi_2(\xi \otimes de_1)\otimes e_2)+$
\ \\
\ \\
\indent $\hspace {0.12in}$ $\psi_2 (m_2(\xi \otimes e_1)\otimes de_2) $
\ \\
\ \\
\indent $=\big[\psi_2 (m_{2}(\xi \otimes de_1)\otimes e_2) +m'_{2}(\psi_2(\xi \otimes e_1)\otimes de_2)\big]+\big[m'_2(m'_2(\psi_1(\xi)\otimes e_1)\otimes e_2)+$
\ \\
\ \\
\indent $\hspace {0.12 in}$ $m'_2((\psi_1(\xi)\otimes e_1)\otimes e_2) \big]+\big[\psi_1(m_2(m_2(\xi\otimes e_1)\otimes e_2))+m'_2((\psi_1(\xi)\otimes e_1)\otimes e_2) \big]$
\ \\
\ \\
\indent $=m'_2(m'_2(\psi_1(\xi)\otimes e_1)\otimes e_2)+\psi_1(m_2(m_2(\xi\otimes e_1)\otimes e_2))$
\ \\
\ \\
\noindent We have the fourth equality is true because equation (\ref {equation:10}) is true for $e_1$ and $e_2$. The first bracket of the last equality is $0$ because $\psi_2 (m_{2}(\xi \otimes de_1)\otimes e_2)=m'_{2}(\psi_2(\xi \otimes e_1)\otimes de_2)=0$. Indeed, $de_i$ ($i=1,2$) is either $0$ or a sum of product(s) which contains a factor of the form $\lefty {e_{\gamma }} \lefty {e_{\gamma^{\dagger}}}$. As a result, the actions of $de_i$ on our complexes $\leftcomplex {T}$ and $\leftcomplex {T'}$ are trivial since $\lefty {T}$ and $\lefty {T'}$ have only one crossing. Therefore,\\ 
\ \\
\noindent $\psi_2(\xi \otimes d(e_1e_2) )$
\ \\
$=m'_2(m'_2(\psi_1(\xi)\otimes e_1)\otimes e_2)+\psi_1(m_2(m_2(\xi\otimes e_1)\otimes e_2))$
\ \\
$=m'_{2}(\psi_1(\xi)\otimes e_1e_2)+\psi_1(m_{2} (\xi \otimes e_1e_2))$
\ \\
\noindent As a result, we have proved that $\Psi=\{(\psi_1,\psi_2)\}$ is an $A_{\infty}$ morphism from $\leftcomplex {T}$ to $\leftcomplex {T'}$ $\Diamond$.\\
\ \\
\ \\
If we define $\Phi$ identically as $\Psi$ but reversing the roles of $\leftcomplex {T'}$ and $\leftcomplex {T}$, we immediately have the following:
\begin {enumerate}
\item $(\Psi \ast \Phi)_1=I_{\leftcomplex {T}}$ since $\psi_1$ and $\phi_1$ are identity maps.
\item $(\Psi \ast \Phi)_2=\psi_1\circ \phi_2+\psi_2\circ (\phi_1\otimes \I)=0$ since both terms are $\phi_2$ under the canonical identification of generators of $\leftcomplex {T}$ and $\leftcomplex {T'}$  
\item $(\Psi \ast \Phi)_3=\psi_2(\phi_2\otimes \I)=0$ since both $\psi_2$ and $\phi_2$ are supported on states $\xi$ where $\xi(c)=1$ and their images contain states whose crossing $c$ is resolved by $0$-resolution.
\end {enumerate}  
Therefore, $\Psi \ast \Phi=1_{\leftcomplex {T}}$. Similarly, we have $\Phi \ast \Psi=1_{\leftcomplex {T'}}$. As a result, $(\leftcomplex {T}, \m {1}, \m {2})$ is isomorphic to $(\leftcomplex {T'},\mm {1}, \mm {2})$ as type $A$ structures.
\subsection {Invariance under Reidemeister Moves}
Due to a bimodule structure introduced in the beginning of this section, we only need to prove the invariance under the Reidemeister moves in the following figures:
$$
\inlinediag{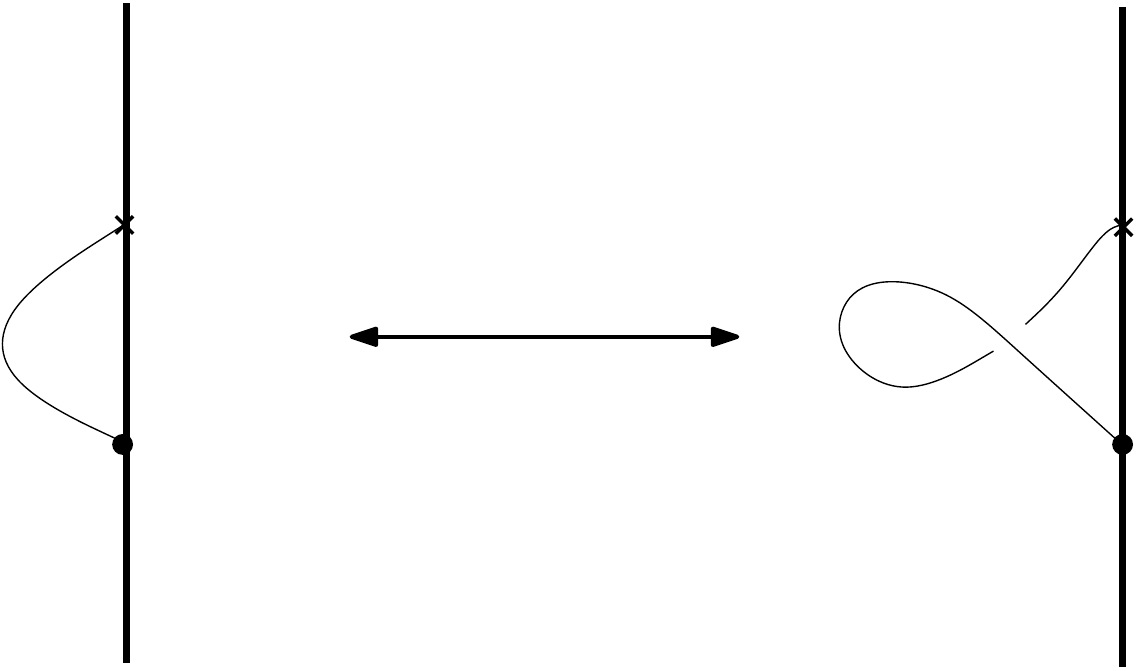} \hspace{0.3in} \inlinediag{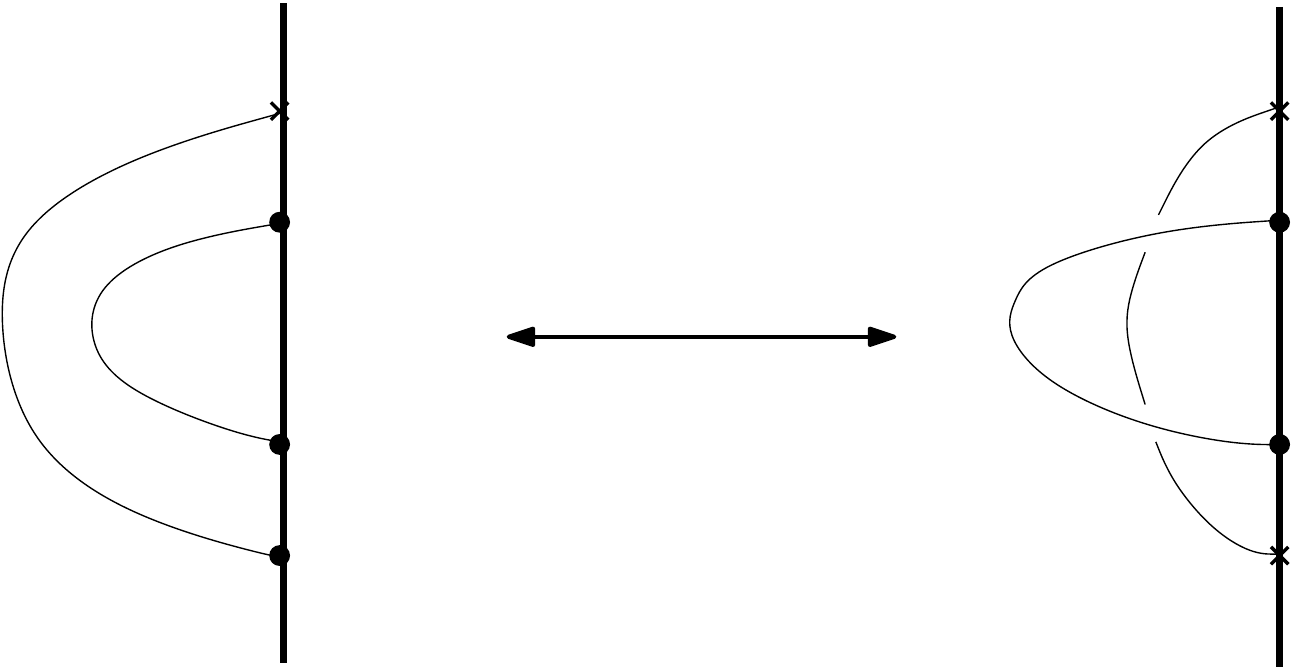} \hspace{0.3in}  \inlinediag{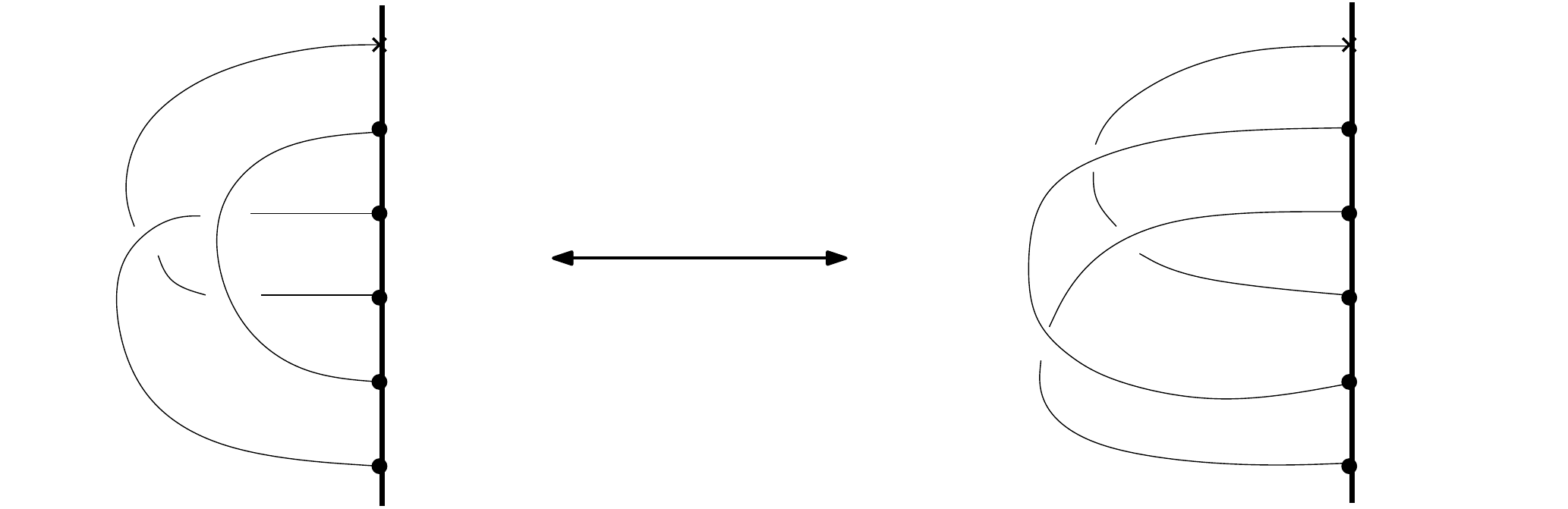}
$$

\noindent The strategy is to use the isomorphism in \ref {sect:amove}  to move the weights close to the $y$-axis. We then can modify Roberts's proofs of the invariance in the untwisted case (see section $5$ and $6$ of \cite {LR2}). Let $\lefty {T}$ and $\lefty {T'}$ be the tangles before and after the Reidemeister moves. We will briefly sketch his arguments for the untwisted case and how we can apply those arguments for our case:\\
\begin {enumerate}
\item Since the weights are close to the $y$-axis, we have $\m {1}=m_1$.
\item Decomposing the complex of tangle $\leftcomplex {T'}$ as a direct sum where each summand corresponds to the resolution of the crossing(s) of the tangle $\lefty {T'}$, there always exists a summand $V$ consisting of states which have a free circle in the generating resolutions. We then decompose $V=V_+ \oplus V_-$, based on the decoration on the free circle.
\item For the differential $m_1$, there are two types of isomorphisms from merging a $+$ free circle to a decorated circle or dividing out a decorated circle to get a new $-$ free circle. We can use this isomorphism to cancel out the summand $V_+$ and its image $d_{APS}(V_+)$  in the case of the Reidemeister move I, or $V_+$ and $d_{APS}(V_+)$ first and then cancel out $V_-$ and its preimage under the dividing isomorphism in the case of the Reidemeister moves II and III. After the cancellation, we get exactly the same chain module as before the Reidemeister move with the possibility that the higher order actions might appear.
\item Roberts proves that the higher order actions actually do not show up because of two properties. The first property is that the image of $m_1$ on $V_+$ is another summand of complex, which is canceled out by the cancellation process. The second property is that the images of higher order actions always lie on $V_+$ and this will be canceled out at the end.
\item For our case, the same technique can be used to prove there is no higher order actions. Since we have $\m {1}=m_1$, we definitely have the first property ($\bf Note.$ If we do not move weights close to the $y$-axis, the image of $\m {1}$ on $V_+$ also contains $V_-$). Additionally, since $\m {2}$ is different from $m_2$ only on the action of left decoration elements, the image of $\m {2}(V_+)$ lies on $V_+$ and therefore the image of the higher order actions lie on $V_+$. At the end, they all disappear when we cancel out $V_+$. The same argument can be applied for $V_-$ if needed (as in Reidemeister move II or III). 
\item After canceling out those terms, we obtain almost the same type $A$ structure associated to the tangles before the Reidemeister moves. The only difference is that we are working over different ground fields. This issue can be handled exactly like the case of the type $D$ structure (see section \ref {sect:IVR}) by using the stable equivalence relation, described in section \ref {sect:stable1}.
\end {enumerate}  
Therefore, pursuant to the bimodule version in twisted tangle homology, we have the following theorem:
\begin {thm}
Let $\lefty {\lk T}$ be a left tangle with a diagram $\lefty {T}$. The homotopy class of $(\leftcomplex {T},\m {1},\m {2})$, defined as in \ref {sect:Adef}, is an invariant of the left tangle $\lefty {\lk T}$.
\end {thm}   

\section{Relation to the totally twisted Khovanov homology by gluing left and right tangles} \label {sect:gluing}
\noindent As described in section \ref {sect:1}, let $T$ be a link diagram for a link $\lk T$ which is divided by the $y$-axis into two parts: a left tangle $\lefty {T}$ and a right one $\righty {T}$. Using the pairing technique in section \ref {sect:pairing}, we will prove that the chain complex $(\leftcomplex {T} \bos \rightcomplex {T},\bys )$ obtained by gluing the type $A$ structure $(\leftcomplex {T},\m {1},\m {2})$ and the type $D$ structure $(\rightcomplex {T},\rightTde {T})$ is chain isomorphic to the totally twisted Khovanov homology of $T$.\\
\ \\
Let $\field {T}$ be the field of fractions of $\ring{T}=\Z_2[x_f|f\in \arcs{T}]$ where $\arcs {T}$ is the set of segments whose endpoints are either the crossings of $T$ or the intersection points of $T$ and the $y$-axis.
\noindent Let $\complex {T}$ be the graded Khovanov complex over $\field {T}$, equipped with a totally twisted differential $\widetilde {\partial}=\partial_{KH}+\partial_{\lk V,T}$ where $\partial_{KH}$ is the regular Khovanov map and $\partial_{\lk V,T}$ is a Koszul map defined similar to $\partial_{\lk V}$ (see \cite {LR3} for more detail). Recall that $(\complex{T},\widetilde {\partial})$ is a link invariant. 

\ \\
\noindent There are natural injections $\phi_l: \field {\lefty {T}}\hookrightarrow \field {T} $ and $\phi_r: \field {\righty {T}}\hookrightarrow \field {T} $, which come from the fact that $\arcs{T}$ is the disjoint union of $\arcs {\lefty {T}}$ and $\arcs {\righty {T}}$. To describe the glued complex, we first need to describe the type $A$ structure $(\lk F_{\phi_l}(\leftcomplex {T}), \lk F_{\phi_l}(\m {1}),\lk F_{\phi_l}(\m{2}))$ and the type $D$ structure $(\lk G_{\phi_r}(\rightcomplex {T}),\lk G_{\phi_r}(\rightTde{T}))$. By using the formulas in Definitions \ref {defn:functor} and \ref {defn:34}, over $\field {T}$, the generators of $\lk F_{\phi_l}(\leftcomplex {T})$ and $\lk G_{\phi_r}(\rightcomplex {T})$ are identified with the generators of $\leftcomplex {T}$ (as a vector space over $\field {\lefty {T}}$) and $\rightcomplex {T}$ (as a vector space over $\field {\righty {T}}$), respectively. Additionally, under this identification, the maps $ \lk F_{\phi_l}(\m {1})$, $ \lk F_{\phi_l}(\m {2})$ and $\lk G_{\phi_r}(\rightTde{T})$ are identical to $\m {1}$, $\m {2}$ and $\rightTde {T}$ respectively.\\
\ \\
\noindent Therefore, we can use $(\leftcomplex {T},\m {1},\m {2})$ (respectively $(\rightcomplex {T},\rightTde {T}$)) to stand for $(\lk F_{\phi_l}(\leftcomplex {T}), \lk F_{\phi_l}(\m {1}),\lk F_{\phi_l}(\m{2}))$ (respectively $(\lk G_{\phi_r}(\rightcomplex {T}),\lk G_{\phi_r}(\rightTde{T}))$). Keep in mind that from now to the end of this section, $\leftcomplex {T}$ and $\rightcomplex {T}$ are vector spaces over $\field {T}$ while $\m {1}$, $\m {2}$ and $\rightTde {T}$ are defined in \ref {sect:Adef} and \ref {sect:3}. The weight of a cleaved circle in type A (respectively type D) is calculated as the sum of weights on the left-side (respectively right-side) arcs which belong to this circle. The module structure of the glued complex then can be described as:
$$
\leftcomplex {T} \bos \rightcomplex {T}=\leftcomplex {T}\otimes_{\lk I_n} \rightcomplex {T}
$$
\noindent Additionally, the differential of this complex is given by the following formula (see section \ref {sect:pairing}):\\
$$
\bys(x\otimes y)=\m {1}(x)\otimes y + (\m {2}\otimes\I)(x\otimes \rightTde {T}(y)) 
$$
\ \\
\noindent We now prove the gluing theorem:\\
\begin {thm}
$(\leftcomplex {T} \bos \rightcomplex {T},\bys )$ is isomorphic to $(\complex {T},\widetilde {\partial})$.
\end {thm} 
\noindent $\bf Proof.$ Due to the module structures of $\leftcomplex {T}$, $\rightcomplex {T}$ as vector spaces over $\field {T}$ and the action of $\lk I_n$ on them, $\leftcomplex {T} \bos \rightcomplex {T}$ is a vector space over $\field {T}$ whose generators are identified with the generators of $\complex {T}$ (see section $7$ of \cite {LR2} for more detail). Furthermore, this identification was proved to preserve the bigrading, and thus, it is $\zeta$-grading preserving. Therefore, it suffices to prove that, under this identification, $\bys=\widetilde {\partial}$. Since $\m {1}=d_{APS}+\partial_{\lk V}$ and $\rightTde {T}=\rightDelta {T}+\rightvde{V}$, the differential $\bys$ can be decomposed as following:\\
$$
\bys (x\otimes y)=\big[d_{APS}(x)\otimes y+(m_2\otimes \I)(x\otimes \rightDelta {T}(y))\big ]+\big[\partial_{\lk V}(x)\otimes y+((\m {2}-m_2)\otimes \I)(x\otimes \rightDelta {T}(y))
$$
 \ \\
\indent $\hspace {0.7in}$$+(\m {2}\otimes \I)(x\otimes \rightvde{V}(y))\big]$
\ \\
\ \\
\noindent From \cite {LR2}, we know that:\\
$$
d_{APS}\otimes \I+(m_2\otimes \I)(\I \otimes \rightDelta {T})=\partial_{KH}
$$
Therefore, we will complete the proof of the theorem if we can show:
\begin {equation} \label {eqn:1}
\partial_{\lk V}(x)\otimes y+((\m {2}-m_2)\otimes \I)(x\otimes \rightDelta {T}(y))+(\m {2}\otimes \I)(x\otimes \rightvde{V}(y))=\partial_{\lk V,T}(x\otimes y)
\end {equation}
\noindent for $x$ and $y$ are generators of $\leftcomplex {T}$ and $\rightcomplex {T}$ such that $I_{\partial (x)}=I_{\partial (y)}$. \\
\ \\
We note that $x\otimes y$ is a resolution $(r,s)$ of $\complex {T}$, consisting of a collection of circles which can be divided into four groups:
\begin {enumerate}
\item Circles decorated by $-$
\item Left free circles decorated by $+$
\item Right free circles decorated by $+$
\item Cleaved circles decorated by $+$
\end {enumerate} 
Therefore, the right hand side of identity \ref {eqn:1} can be written as following:
$$
\partial_{ {\lk V,T}}(x\otimes y)=\sum\limits_{C \in {(2)\dot \bigcup(3)\dot \bigcup(4)}} w_C.(r,s_C)
$$
where $w_C=\sum \limits_{f\in \arcs {C}} x_f$. 
\noindent On the other hand,\\
$$
\partial_{\lk V}(x)\otimes y=\sum \limits_{C \in (2)} w_C.(r,s_C)
$$ 
Additionally, according to the construction of $\m {2}$, the action $\m {2}-m_2$ is supported only on left decoration elements. As a result:
$$
((\m {2}-m_2)\otimes \I)(x\otimes \rightDelta {T}(y))=\sum \limits_{C\in (4)}((\m {2}-m_2)\otimes \I)(x\otimes \lefty {e_C}\otimes y_C)
$$
$$
\hspace*{1.7cm} =\sum\limits_{C \in (4)} \lefty {w_C}x_C \otimes y_C
$$
$$
\hspace*{1.3cm}=\sum\limits_{C \in (4)} \lefty {w_C}(r,s_C)
$$
From the definition of $\rightvde {V}$, we can calculate the last term of the left hand side of \ref {eqn:1} as:
$$
(\m {2}\otimes \I)(x\otimes \rightvde{V}(y))=\sum \limits_{C \in  (4)}\righty {w_C}.(\m {2}\otimes \I)(x\otimes \righty {e_C}\otimes y_C)+\sum \limits_{C \in  (3)}w_C.(\m {2}\otimes \I)(x\otimes I_{\partial (y_C)}\otimes y_C)
$$
$$
\hspace*{-2cm}=\sum \limits_{C \in  (4)} \righty {w_C}(r,s_C)+\sum \limits_{C \in  (3)} w_C(r,s_C)
$$
\noindent Rewriting the left hand side of identity \ref {eqn:1}, we have:\\
$$
LHS= \sum \limits_{C \in  (2)}w_C.(r,s_C)+\sum \limits_{C \in  (3)}w_C.(r,s_C)+\sum \limits_{C \in  (4)}(\lefty {w_C}+\righty {w_C}).(r,s_C)
$$
Since $w_C=\lefty {w_C}+\righty {w_C}$ for each cleaved circle, we can conclude that identity \ref {eqn:1} is true and thus, $\bys=\widetilde {\partial}$. Thus, $(\leftcomplex {T} \bos \rightcomplex {T},\bys )$ is chain isomorphic to $(\complex {T},\widetilde {\partial})$. $\Diamond$.
\section {Examples of the type D and the type A structures} \label {sect: example}
\noindent In this section, we will give examples calculating the totally twisted Khovanov type $D$ and type $A$ structures for several knots and links.\\
\noindent 
\ \\
$\bf Example$ $\bf 1:$ {Hopf link}
\ \\
Consider the Hopf link $T$ which is transverse to the $y$-axis at 4 points. It divides the link $T$ into two parts: the left tangle $\lefty {T}$ and the right one $\righty {T}$. Label each arc of $T$ as in the following picture.    \\
$$
\inlinediag[0.6]{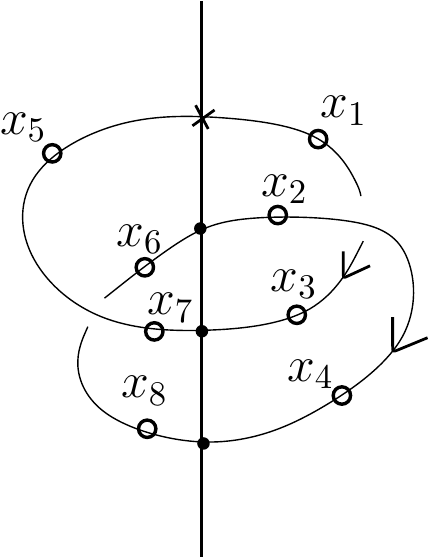}
$$ 
We first describe the type $D$ structure on $\rightcomplex {T}$. As shown in the following figure, $\rightcomplex {T}$ can be thought as a vector space over $\field {T}$ (by the same argument as in section \ref {sect:gluing}), generated by six elements: $\righty {\xi_1}$,..., $\righty {\xi_6}$ corresponding to the bottom row of the figure from left to right. Their $\ze$-gradings are $\frac {-1}{4}$, $\frac {1}{4}$, $0$, $\frac {1}{4}$, $\frac {3}{4}$ and $\frac {1}{2}$ respectively. For example: since the bigrading of $\righty {\xi_4}$ is $(-1,-5/2)$ due to $n_+(\righty {T})=0$ and $n_-(\righty {T})=2$, the $\ze$-grading of $\righty {\xi_4}$ is $-1+5/4=1/4$. We also denote $\{\lefty {\eta_i}\}_{i=1,..,4}$ and $\{\righty {\gamma_j}\}_{j=1,2}$ the left bridges and right active resolution bridges as in the figure. Let $C$ and $D$ be the $+$ circles of $\righty {\xi_1}$ and $\righty {\xi_4}$ respectively. 
$$
\inlinediag[1]{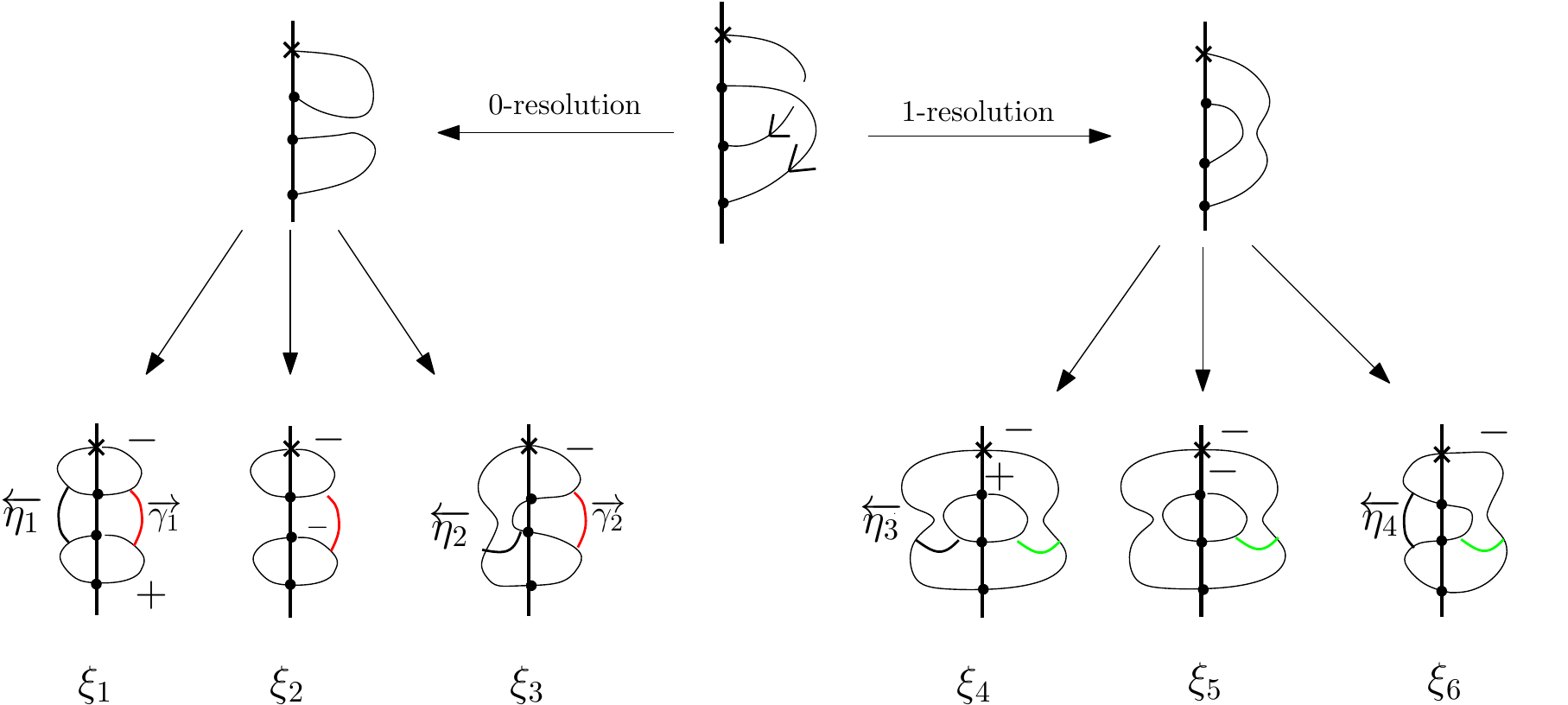}
$$ 
\ \\
We now can describe $\rightTde {T}$:
\begin{equation} \label {eqn:hopfd}
\begin {split}
\rightTde {T}(\righty {\xi_1})&=e_{\lefty {\eta_1}}\otimes \righty {\xi_3}+e_{\righty {\gamma_1}}\otimes \righty {\xi_6}+[\lefty {e_C}+(x_3+x_4)\righty {e_C}]\otimes \righty {\xi_2} \\
\rightTde {T}(\righty {\xi_3})&=e_{\righty {\gamma_2}}\otimes \righty {\xi_5}+e_{\lefty {\eta_2}}\otimes  \righty {\xi_2}\\
\rightTde {T} (\righty {\xi_4})&=e_{\lefty {\eta_3}}\otimes  \righty {\xi_6}+[\lefty {e_D}+(x_2+x_3)\righty {e_D}]\otimes \righty {\xi_5}\\
\rightTde {T} (\righty {\xi_6})&=e_{\lefty {\eta_4}}\otimes  \righty {\xi_5}
\end {split}
\end{equation}
We next describe the type $A$ structure on $\leftcomplex {T}$. Similarly to $\rightcomplex {T}$, $\leftcomplex {T}$ is a six dimensional vector space over $\field {T}$, generated by $\{\lefty {\xi_i}\}_{i=\overline {{1,6}}}$ where $\lefty {\xi_i}$ has the same boundary as $\righty {\xi_i}$. We label the bridges of each generator $\lefty {\xi_i}$ exactly as we labeled the ones for $\righty {\xi_i}$ in the type $D$ structure. For example: $\lefty {\xi_1}$ has a left active resolution bridge $\lefty {\eta_1}$ and right bridge $\righty {\gamma_1}$. We also let $\righty {\gamma_3}$ and $\righty {\gamma_4}$ be the right bridges of $\lefty {\xi_4}$ and $\lefty {\xi_6}$ respectively. Additionally, the $\ze$-grading of $\lefty {\xi_1}$,...,$\lefty {\xi_6}$ are $\frac {-1}{4}$, $\frac {1}{4}$, $\frac {1}{2}$, $\frac {1}{4}$, $\frac {3}{4}$ and $0$, respectively. Since none of generators has a free circle, the type $A$ structure on $\leftcomplex {T}$ can be described as follows:
\begin{equation} \label {eqn:hopfa}
\begin {split}
\m {1}&=0\\
\m {2} (\lefty {\xi_1}\otimes e_{\lefty {\eta_1}})&=\lefty {\xi_3}\\
\m {2} (\lefty {\xi_1}\otimes e_{\righty {\gamma_1}})&=\lefty {\xi_6}\\
\m {2} (\lefty {\xi_1}\otimes \righty {e_C}) &= \lefty {\xi_2}\\
\m {2} (\lefty {\xi_1}\otimes \lefty {e_C}) &= (x_7+x_8)\lefty {\xi_2}\\
\m {2} (\lefty {\xi_3}\otimes e_{\righty {\gamma_2}})&= \lefty {\xi_5}\\
\m {2} (\lefty {\xi_4}\otimes e_{\righty {\gamma_3}})&=\lefty {\xi_3}\\
\m {2} (\lefty {\xi_4}\otimes \righty {e_D}) &= \lefty {\xi_5}\\
\m {2} (\lefty {\xi_4}\otimes \lefty {e_D}) &= (x_6+x_7)\lefty {\xi_5}\\
\m {2} (\lefty {\xi_6}\otimes e_{\righty {\gamma_4}})&=\lefty {\xi_2}\\
\m {2} (\lefty {\xi_6}\otimes e_{\lefty {\eta_4}})&=\lefty {\xi_5}
\end {split}
\end{equation}
Then $\leftcomplex {T} \boxtimes \rightcomplex {T}$ is a graded vector space over $\field {T}$, generated by $\{\xi_i\}_{i=\overline {1,6}}$ where $\xi_i=\lefty {\xi_i}\otimes_{\lk I_2}\righty {\xi_i}$. Since $\ze(\xi_i)=\ze(\lefty {\xi_i})+\ze(\righty {\xi_i})$, the $\ze$-grading of $\xi_1$,...,$\xi_6$ are $\frac {-1}{2}$, $\frac {1}{2}$, $\frac {1}{2}$, $\frac {1}{2}$, $\frac {3}{2}$ and $\frac {1}{2}$ respectively. Since $\m {1}=0$, the differential of  $\leftcomplex {T} \boxtimes \rightcomplex {T}$ is:
$$
\partial^{\boxtimes}(x\otimes y)=(\m {2}\otimes \I)(x\otimes \rightTde {T} (y))
$$
Therefore, using \ref {eqn:hopfd} and \ref {eqn:hopfa}, we have:
\begin{equation}
\begin {split}
\partial^{\boxtimes}(\xi_1)&=\xi_3+\xi_6+(x_3+x_4+x_7+x_8)\xi_2\\
\partial^{\boxtimes}(\xi_3)&=\xi_5\\
\partial^{\boxtimes}(\xi_4)&=(x_2+x_3+x_6+x_7)\xi_5 \\
\partial^{\boxtimes}(\xi_6)&=\xi_5\\
\end {split}
\end{equation}
The above description of $\ze$-complex $(\leftcomplex {T} \boxtimes \rightcomplex {T},\partial^{\boxtimes})$ is exactly the same as the totally twisted Khovanov homology of Hopf link whose homology is $\field {T}\oplus \field {T}$, occurring in the $\ze$-grading $\frac {1}{2}$.
\ \\
\ \\
\ \\
$\bf Example$ $\bf 2:$ Consider the following tangle diagram $\righty {T}$:
\ \\
\ \\
$$
\inlinediag[0.5]{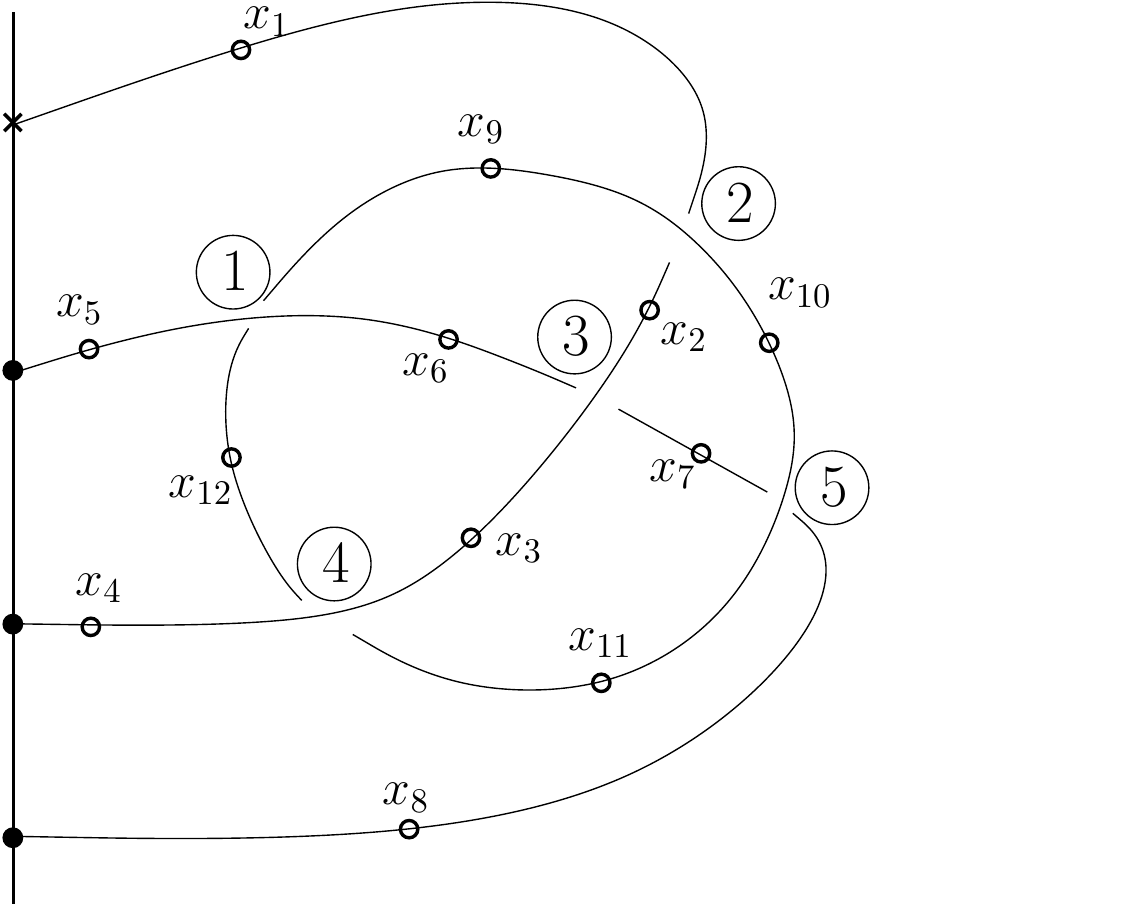}
$$ 
We label five crossings of $\righty {T}$ as in the above figure. In this example, we will investigate the image of $\rightDelta {n,T}$ on some generators of $\rightcomplex {CT}$. Recall: $(\rightcomplex {CT},\rightDelta {n,T})$ is the deformation retraction of $(\rightcomplex {T},\rightTde {T})$, defined in Section \ref  {sect:1}. For this tangle $\righty {T}$, there are two left planar matchings $m_1$ (containing an arc which connects the mark point and the bottommost point on the $y$-axis) and $m_2$. A generator of $\rightcomplex {CT}$ is obtained by resolving the crossings of $\righty {T}$ by either a $0$-resolution or a $1$-resolution in such a way that it has no free circle, then capping it off by either $m_1$ or $m_2$ and finally, decorating the unmarked cleaved circle (it might not exist) by $\pm$. Therefore, we can encode a generator of $\rightcomplex {CT}$ in the form of either $\xi_{\ast\ast\ast\ast\ast,i}$ or $\xi_{\ast\ast\ast\ast\ast,i,\pm}$ where each $\ast$ receives the value in $\{0,1\}$, depending on the resolution of the corresponding crossing. $i$ receives the value of either $1$ or $2$, depending on the way we choose either left planar matching $m_1$ or $m_2$. Additionally, if the generator has an unmarked cleaved circle, we use the second form to specify the decoration on the unmarked circle. Otherwise, we use the first form to represent this generator. For example, the following figure represents $\xi _{00000,1}$. 
$$
\inlinediag[0.5]{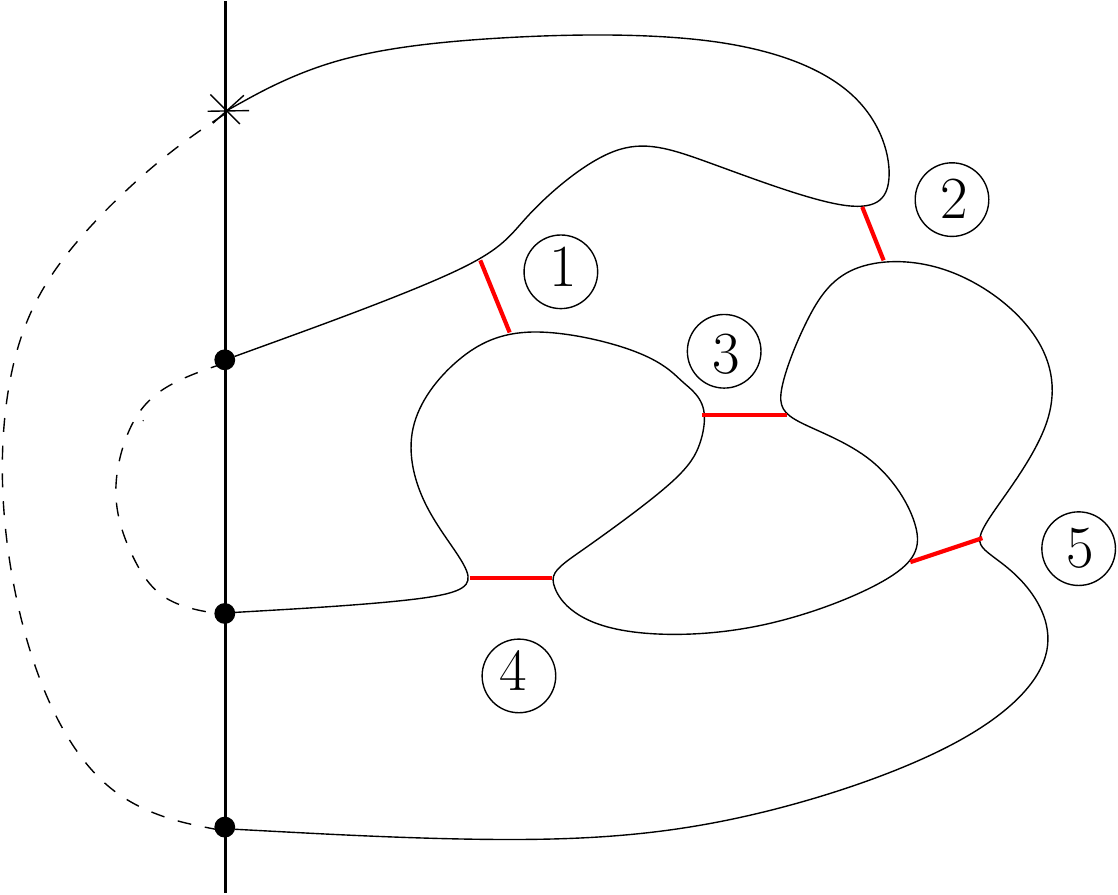}
$$ 
We have:
\begin {equation} \label {eqn:ex}
\begin{split}
\rightDelta {n,T} (\xi_{00000,1})=[x_3+x_6+&x_{12}]^{-1} I_{\partial (\xi_1)} \otimes \xi_{10010,1}+ [x_2+x_7+x_{10}]^{-1}I_{\partial (\xi_1)} \otimes {\xi_{01001,1}}\\
&+([x_3+x_7+x_{11}]^{-1}+[x_2+x_7+x_{10}]^{-1}) I_{\partial (\xi_1)} \otimes \xi_{00101,1}\\
&+([x_3+x_6+x_{12}]^{-1}+[x_3+x_7+x_{11}]^{-1}) I_{\partial (\xi_1)}\otimes \xi_{00110,1} \\
&+e_{\righty {\gamma_1}}\otimes \xi_{10000,1,-}+e_{\righty {\gamma_2}}\otimes \xi_{01000,1,-}\\
&+e_{\lefty {\eta_1}} \otimes \xi_{00000,2,-}
\end{split}
\end {equation}
where $I_{\partial (\xi_1)}:=I_{\partial (\xi_{00000,1})}$, $\{e_{\righty {\gamma_i}}\}_{i=1,2}$ are the right bridge elements of $\mathcal{B}\Gamma_{2}$, corresponding to the change of the cleaved links after surgering $\xi_{00000,1}$ along the active resolution bridge at crossing $i$. $\lefty {\eta_1}$ is the unique left bridge of the left planar matching $m_1$.
\ \\
\ \\
\noindent As we can see, in the right hand side of equation \ref {eqn:ex}, the first (or the second term) comes from a generator, obtained from $\xi_{00000,1}$ by surgering along two active resolution bridges. Depending on which crossing we resolve first, there are two ways to make this surgery. However, there is only one way which creates a free circle and the second way will change the boundary. On the other hand, both ways to surger $\xi_{00000,1}$ to obtain a generator of  the third (or fourth) term contain free circles, and, as the result, we count the coefficients for both paths. 
\ \\
\ \\
Similarly, we compute $\rightDelta {n,T} (\xi_{00000,2,+})$:
\begin {equation} \label {eqn:ex1}
\begin{split}
\rightDelta {n,T} (\xi_{00000,2,+}&)=[x_3+x_6+x_{12}]^{-1} I_{\partial (\xi_2)} \otimes \xi_{10010,2,+}+ [x_2+x_7+x_{10}]^{-1}I_{\partial (\xi_2)} \otimes {\xi_{01001,2,+}}\\
&+([x_3+x_7+x_{11}]^{-1}+[x_2+x_7+x_{10}]^{-1}) I_{\partial (\xi_2)} \otimes \xi_{00101,2,+}\\
&+([x_3+x_6+x_{12}]^{-1}+[x_3+x_7+x_{11}]^{-1}) I_{\partial (\xi_2)}\otimes \xi_{00110,2,+} \\
&+e_{\righty {\gamma_1}}\otimes \xi_{10000,2}+e_{\righty {\gamma_2}}\otimes \xi_{01000,2}\\
&+e_{\lefty {\eta_2}} \otimes \xi_{00000,2}+[\lefty {e_C}+(x_2+x_3+x_4+x_6+x_7+x_8+x_{10}+x_{11}+x_{12})\righty {e_C}] \otimes \xi_{00000,2,-}
\end{split}
\end {equation}
where $I_{\partial ({\xi_2})}=I_{\partial (\xi_{00000,2,+})}$. $\lefty {\eta_2}$ is a unique left bridge of $m_2$ and $C$ stands for the unmarked cleaved circle of $\xi_{00000,2,+}$. We note that the main difference between equation \ref {eqn:ex} and equation \ref {eqn:ex1} is that the image of $\rightDelta {n,T}$ on $\xi_{00000,2,+}$ contains a term coming from a decoration element.

\end{document}